\documentclass[11pt]{article}

\usepackage[a4paper,margin=1in]{geometry}
\usepackage[T1]{fontenc}
\usepackage[utf8]{inputenc}
\usepackage{lmodern}
\usepackage{amsmath,amssymb,amsthm,mathtools}
\usepackage{mathrsfs}
\usepackage{enumitem}
\usepackage{needspace}
\usepackage{microtype}
\usepackage{graphicx,booktabs}
\usepackage[
  backend=biber,
  style=numeric,
  sorting=nyt,
  giveninits=true,
  maxbibnames=99,
  doi=true,
  eprint=true,
  url=true,
  isbn=false
]{biblatex}
\AtBeginBibliography{\small}
\usepackage[hidelinks]{hyperref}
\hypersetup{pdftitle={Brown Measures and Spectral Convergence for Arbitrary GUE Polynomials},
  pdfauthor={Benoit Collins and Yu Kitagawa and Akihiro Miyagawa and Wangjun Yuan and Ping Zhong}}
\usepackage[nameinlink,capitalize,noabbrev]{cleveref}

\allowdisplaybreaks[4]
\numberwithin{equation}{section}

\newtheorem{theorem}{Theorem}[section]
\newtheorem{proposition}[theorem]{Proposition}
\newtheorem{lemma}[theorem]{Lemma}
\newtheorem{corollary}[theorem]{Corollary}
\theoremstyle{definition}
\newtheorem{definition}[theorem]{Definition}
\newtheorem{example}[theorem]{Example}
\theoremstyle{remark}
\newtheorem{remark}[theorem]{Remark}
\crefname{theorem}{Theorem}{Theorems}
\Crefname{theorem}{Theorem}{Theorems}
\crefname{proposition}{Proposition}{Propositions}
\Crefname{proposition}{Proposition}{Propositions}
\crefname{lemma}{Lemma}{Lemmas}
\Crefname{lemma}{Lemma}{Lemmas}
\crefname{corollary}{Corollary}{Corollaries}
\Crefname{corollary}{Corollary}{Corollaries}
\crefname{definition}{Definition}{Definitions}
\Crefname{definition}{Definition}{Definitions}
\crefname{example}{Example}{Examples}
\Crefname{example}{Example}{Examples}
\crefname{remark}{Remark}{Remarks}
\Crefname{remark}{Remark}{Remarks}

\newcommand{\C}{\mathbb C}
\newcommand{\R}{\mathbb R}
\newcommand{\E}{\mathbb E}
\newcommand{\Pp}{\mathbb P}
\newcommand{\tr}{\operatorname{tr}}
\newcommand{\Tr}{\operatorname{Tr}}
\newcommand{\sg}{\sigma_{\min}}
\newcommand{\HS}{\mathrm{HS}}
\newcommand{\End}{\operatorname{End}}
\newcommand{\Span}{\operatorname{span}}
\newcommand{\ran}{\operatorname{Ran}}
\newcommand{\supp}{\operatorname{supp}}
\newcommand{\dd}{\,\mathrm{d}}
\newcommand{\Lap}{\nabla_z^2}
\newcommand{\Gr}{\mathrm{Gr}}
\newcommand{\ip}[2]{\langle #1,#2\rangle}

\newcommand{\spec}{\operatorname{spec}}
\newcommand{\id}{\operatorname{id}}
\title{\bfseries Brown Measures and Spectral Convergence\\
for Arbitrary GUE Polynomials}
\author{Beno{\^\i}t Collins\\[0.4em]
  \normalsize Department of Mathematics\\
  \normalsize Graduate School of Science\\
  \normalsize Kyoto University\\
  \normalsize Kyoto 606-8502, Japan\\
  \normalsize \href{mailto:collins@math.kyoto-u.ac.jp}{\texttt{collins@math.kyoto-u.ac.jp}}
  \and
  Yu Kitagawa\\[0.4em]
  \normalsize Department of Mathematics\\
  \normalsize Graduate School of Science\\
  \normalsize Kyoto University\\
  \normalsize Kyoto 606-8502, Japan\\
  \normalsize \href{mailto:kitagawa.yu.57z@st.kyoto-u.ac.jp}{\texttt{kitagawa.yu.57z@st.kyoto-u.ac.jp}}
  \and
  Akihiro Miyagawa\\[0.4em]
  \normalsize Department of Mathematics\\
  \normalsize Graduate School of Science\\
  \normalsize Kyoto University\\
  \normalsize Kyoto 606-8502, Japan\\
  \normalsize \href{mailto:miyagawa.akihiro.8p@kyoto-u.ac.jp}{\texttt{miyagawa.akihiro.8p@kyoto-u.ac.jp}}
  \and
  Wangjun Yuan\\[0.4em]
  \normalsize Department of Mathematics\\
  \normalsize Southern University of\\
  \normalsize Science and Technology\\
  \normalsize Shenzhen, China\\
  \normalsize \href{mailto:ywangjun@connect.hku.hk}{\texttt{ywangjun@connect.hku.hk}}
  \and
  Ping Zhong\\[0.4em]
  \normalsize Department of Mathematics\\
  \normalsize University of Houston\\
  \normalsize Houston, TX 77204, USA\\
  \normalsize \href{mailto:pzhong@central.uh.edu}{\texttt{pzhong@central.uh.edu}}}
\date{}

\begin{document}
\maketitle
\begin{abstract}
We prove that, for every noncommutative polynomial in independent GUE
matrices, the empirical eigenvalue distribution converges almost surely
to the Brown measure of the same polynomial in free semicircular
variables. This is the Gaussian case of a conjecture of Belinschi,
\'Sniady and Speicher. We also prove that the spectrum converges almost
surely in Hausdorff distance to the spectrum of the limiting operator,
so that the spectral radius converges as well; deviations have
probability smaller than every inverse power of the matrix size. For
this we show that the support of the Brown measure of every polynomial
in free semicircular variables is its whole spectrum, a statement about
operators that involves no random matrices. We also give an algorithm,
with a convergence guarantee, that computes the limiting Brown measure
from matrices whose size depends only on the polynomial. A holomorphic
rotation introduces an independent Ginibre component and yields the
almost-sure limit without a least-singular-value bound for the original
GUE polynomial. We also give a second-moment criterion for weak
convergence of Brown measures. On the free probability side, we use the Dyson
equation of a linear pencil to show that harmonicity of its logarithmic
potential forces invertibility.
\end{abstract}
\setcounter{tocdepth}{2}
\clearpage
{\small\tableofcontents}
\clearpage

\section{Introduction}\label{sec:intro}

\subsection{The problem and the results}
Let $G_1^{(n)},\ldots,G_k^{(n)}$ be independent normalized GUE matrices
of size $n$, and let $P$ be a polynomial in $k$ noncommuting variables
with complex coefficients. When $P(G^{(n)})$ is selfadjoint, the asymptotic behaviour of its eigenvalues is well understood. Voiculescu's
asymptotic freeness theorem \cite{Voi} identifies the limit of the joint
$*$-moments of the matrices with those of a free semicircular family,
and Haagerup and Thorbj\o rnsen \cite{HT} proved that the norms of all
polynomials converge almost surely. For selfadjoint polynomials, these
theorems imply both weak convergence of eigenvalue measures and
Hausdorff convergence of spectra.

\begin{figure}[htbp]
 \centering
 \includegraphics[width=\textwidth]{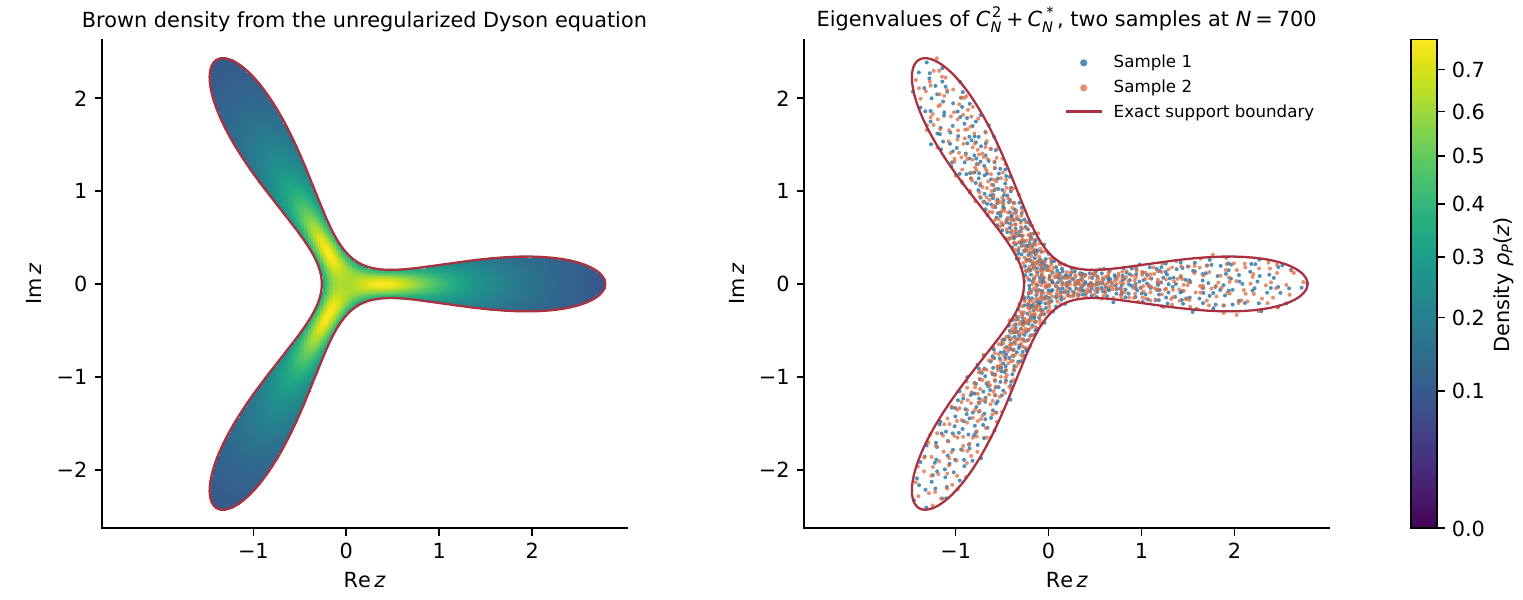}
 \caption{The Brown measure of $c^2+c^*$ and its Ginibre approximation.
 Left: the density computed from the unregularized Dyson equation
 on the interior branch selected by the limit $\eta\downarrow0$.
 Right: two independent samples of the eigenvalues of $C_N^2+C_N^*$,
 where $C_N$ denotes a complex Ginibre matrix of size $N=700$.
 The red curve in both panels is the boundary derived in
 \cref{alg:adjoint-measure}.}
 \label{fig:adjoint-example}
\end{figure}

For a nonselfadjoint polynomial, the eigenvalues can be complex.
Consider, for example, $C_N^2+C_N^*$ for a normalized complex Ginibre
matrix $C_N$ of size $N$. Its limiting eigenvalue distribution is the Brown measure
of $P=c^2+c^*$ for a standard circular element $c$
(\cref{fig:adjoint-example}). The threefold symmetry follows from
$\omega c\overset{*\text{-law}}=c$ and
$(\omega c)^2+(\omega c)^*=\omega^2P$ for $\omega=e^{2\pi i/3}$.
Even the convergence of eigenvalue measures for a nonnormal polynomial
does not follow from convergence of $*$-moments and norms.

In free probability theory, the random variables considered
here are operators in a finite von Neumann algebra with a faithful
normalized trace $\tau$, which plays the role of expectation.
Unital subalgebras are freely independent if a product of centered
factors has trace zero whenever consecutive factors come from different
subalgebras \cite{Voi,VDN}. Elements are freely independent if the
unital subalgebras they generate are free.
Let $s_1,\ldots,s_k$ be freely independent semicircular variables of
variance one, and put
\[
 A_n=P(G_1^{(n)},\ldots,G_k^{(n)}),\qquad A=P(s_1,\ldots,s_k),
 \qquad \mu_{A_n}=\frac1n\sum_{j=1}^n\delta_{\lambda_j(A_n)}.
\]
Write $\mu_A$ for the Brown measure of $A$, the probability measure
whose logarithmic potential is $U_A(z)=\log\Delta_\tau(A-z1)$,
where $\Delta_\tau$ is the Fuglede--Kadison determinant. The potential
determines the measure by $\Lap U_A=2\pi\mu_A$ in the distributional
sense, with $\Lap$ the planar Laplacian. For normal $A$, the Brown
measure is its spectral distribution, and for finite matrices it
is the eigenvalue counting measure.

Belinschi, \'Sniady, and Speicher conjectured that $\mu_{A_n}$ converges
to $\mu_A$ \cite[Section 7]{BSS}; their conjecture also covers Wigner
matrices. The convergence of complex eigenvalue distributions
and the question of a matching lower bound for
the spectral radius are highlighted in \cite[Section 6.9]{vH2026}. There are two difficulties:
(i) the logarithm in Hermitization can amplify a very small singular
value, so moment convergence alone does not identify the limiting
Brown measure, and (ii) strong convergence, that is, convergence of $*$-moments and norms,
excludes spectral outliers,
but the support of the Brown measure can be strictly smaller than the operator spectrum.
Our first theorem overcomes (i); our third theorem overcomes (ii).

\begin{theorem}[Convergence of the eigenvalue distribution]\label[theorem]{thm:main}
Let $P$ be a noncommutative polynomial with complex coefficients.
Almost surely,
\begin{align}
 \mu_{A_n}&\longrightarrow\mu_A\quad\text{weakly},\label{eq:main-weak}\\
 U_{A_n}&\longrightarrow U_A\quad\text{locally in }L^1(\C).\label{eq:main-potential-as}
\end{align}
Moreover, with $b(B)=\int|z|^2\dd\mu_B(z)$:
\begin{enumerate}[label=(\alph*)]
\item for every compact $K\subset\C$,
\begin{equation}\label{eq:main-potential}
 \E\int_K|U_{A_n}(z)-U_A(z)|\dd^2z\longrightarrow0,
\end{equation}
where $\dd^2z$ denotes planar Lebesgue measure;
\item $b(A_n)\to b(A)$ almost surely;
\item for every compact $K\subset\C$ and every fixed $t,q>0$,
\begin{equation}\label{eq:main-potential-tail}
 \Pp\{\|U_{A_n}-U_A\|_{L^1(K)}>t\}
 +\Pp\{|b(A_n)-b(A)|>t\}=O_{P,K,t,q}(n^{-q}).
\end{equation}
\end{enumerate}
\end{theorem}

Here $P$ is fixed, and matrices of different sizes may be coupled
arbitrarily on one probability space. The theorem settles the GUE case
of the Belinschi--\'Sniady--Speicher conjecture and answers, for GUE
polynomials, the question on complex eigenvalue distributions in
\cite[Section 6.9]{vH2026}. By \cref{subsec:free-inputs}, it applies
equivalently to arbitrary $*$-polynomials in independent complex
Ginibre matrices, answering the Gaussian case of the question on
adjoints in \cite[Question 1.8]{Han}.
For Ginibre polynomials without adjoints, Cook, Guionnet, and Husson
\cite{CGH} proved weak convergence to the Brown measure in degree at most two,
and Han \cite{Han} proved it in arbitrary degree. Neither result applies to the polynomial in \cref{fig:adjoint-example}, which involves an adjoint.

Shortly after the first version of this paper appeared, Han
\cite[Theorems 1.1 and 1.4]{HanWigner} extended weak convergence of
the eigenvalue measures to polynomials in independent Wigner matrices
and to $*$-polynomials in independent matrices with i.i.d.\ entries,
under density and moment assumptions. His convergence, including local
$L^1$ convergence of the logarithmic potentials, holds in probability.
For Gaussian inputs, \cref{thm:main} gives almost-sure convergence
and the deviation bounds \eqref{eq:main-potential-tail}. The results of \cref{thm:hausdorff,sp:main,thm:algorithm} on Hausdorff convergence, coincidence of the
support of the Brown measure with the spectrum, and computation of the Brown measure are not
covered by \cite{HanWigner}.

Our second theorem concerns the spectrum. For nonempty compact sets $K,L\subset\C$ we use
\[
 d_H(K,L)=\max\left\{\sup_{z\in K}\operatorname{dist}(z,L),
                         \sup_{z\in L}\operatorname{dist}(z,K)\right\},
\]
and we write $r(B)=\max_{z\in\spec B}|z|$ for the spectral radius.

\begin{theorem}[Convergence of the spectrum]\label[theorem]{thm:hausdorff}
Let $P$ be a noncommutative polynomial with complex coefficients.
Almost surely,
\begin{equation}\label{eq:main-hausdorff}
 d_H(\spec A_n,\spec A)\longrightarrow0,\qquad
 r(A_n)\longrightarrow r(A).
\end{equation}
\begin{enumerate}[label=(\alph*)]
\item For every $\epsilon,q>0$,
\begin{equation}\label{eq:main-hausdorff-tail}
 \Pp\{d_H(\spec A_n,\spec A)>\epsilon\}
 =O_{P,\epsilon,q}(n^{-q}).
\end{equation}
\item For every $\epsilon>0$ there are $c_\epsilon,C_\epsilon>0$
and $n_\epsilon$ such that
\begin{equation}\label{eq:main-outer}
 \Pp\left\{\sup_{\lambda\in\spec A_n}
       \operatorname{dist}(\lambda,\spec A)\ge\epsilon\right\}
 \le C_\epsilon e^{-c_\epsilon n}\qquad(n\ge n_\epsilon).
\end{equation}
\end{enumerate}
\end{theorem}

For GUE polynomials this gives the matching lower bound for the
spectral radius discussed in \cite[Section 6.9]{vH2026}.
The proof in \cref{sp:matrix} first shows that spectral outliers have exponentially small probability. It then combines the potential convergence of \cref{thm:main} with \cref{sp:main} to prove Hausdorff convergence.

\begin{theorem}[Support of the Brown measure]\label[theorem]{sp:main}
For every scalar noncommutative polynomial $P$ with complex coefficients
in a free standard semicircular family,
\begin{equation}\label{sp:main-equality}
 \supp\mu_{P(s)}=\spec P(s).
\end{equation}
\end{theorem}

\begin{remark}\label[remark]{sp:connected}
The $C^*$-algebra generated by a free semicircular
family has no nontrivial projections
\cite{PimsnerVoiculescu,Voiculescu1993KGroups};
see also \cite[Theorem 7.1]{HST}.
Consequently, every element of this algebra has
connected spectrum. By \cref{sp:main}, the support
of the Brown measure of $P(s)$ is therefore connected.
\end{remark}

\Cref{sp:main} is purely operator-theoretic and does not rely on random-matrix approximations.
It answers, for polynomials, a question of Miyagawa \cite{Miyagawa}, who asks it for
polynomials and rational functions in free semicircular variables.
The polynomial hypothesis is crucial here: the support of the Brown measure is
always contained in the operator spectrum, but the reverse inclusion
fails even for some elements of the $C^*$-algebra generated by a
semicircular family; see \cref{sp:nonpolynomial-example}.
For nonconstant $P$, the results of Shlyakhtenko--Skoufranis
\cite{SS} imply that $U_{P(s)}$ is finite everywhere;
see \cref{sp:finite-potential}. In particular $\mu_{P(s)}$ has no
atoms, since an atom forces the logarithmic potential to be
$-\infty$ there.
An isolated support point would be an atom, so the equality above
also shows that $\spec P(s)$ has no isolated points when $P$ is nonconstant.

We also give an algorithm for computing the limiting Brown measure using a Dyson equation and a formula for the Fuglede--Kadison determinant.

\begin{theorem}[Computing the Brown measure]\label[theorem]{thm:algorithm}
Let $P$ be a noncommutative polynomial with complex coefficients.
Its coefficients explicitly determine an affine pencil $L_z$,
$z\in\C$, of a size $m$ depending only on $P$, such that
$U_A(z)=m\log\Delta_{\tau_m}(L_z)$. For $\eta>0$ put
\[
 U_\eta(z)=\tfrac12(\Tr\otimes\tau)\log(L_z^*L_z+\eta^2I),\qquad
 \rho_\eta(z)=(2\pi)^{-1}\Lap U_\eta(z)\ge0.
\]
\begin{enumerate}[label=(\alph*)]
\item For each $z\in\C$ and $\eta>0$, a finite-dimensional matrix
Dyson equation has a unique solution with negative imaginary part,
and a convergent matrix iteration computes it. An explicit formula
gives $U_\eta(z)$ from this solution, and two linear systems give
$\rho_\eta(z)$.
\item As $\eta\downarrow0$,
\[
 U_\eta(z)\downarrow U_A(z)\ \text{for every }z\in\C,\qquad
 U_\eta\longrightarrow U_A\ \text{in }L^1_{\mathrm{loc}}(\C),
 \qquad \rho_\eta\dd^2z\longrightarrow\mu_A
 \ \text{vaguely},
\]
the pointwise limit being understood in the extended real sense.
\end{enumerate}
The construction applies at every spectral parameter, including $z=P(0)$.
More generally, any admissible pencil in the sense of
\eqref{alg:direct-pencil}--\eqref{alg:det-normalization} may be used
in place of $L_z$.
\end{theorem}

The construction is given in \cref{alg:algorithm}: the pencil in
\eqref{alg:direct-pencil}--\eqref{alg:det-normalization}, its Dyson equation in \eqref{alg:dyson},
the iteration in \cref{alg:iteration}, the formula for $U_\eta$ in
\eqref{alg:action}, and the linear systems in \eqref{alg:derivative-system}.
Thus the measure computed by the algorithm is the weak limit in
\cref{thm:main}, and its support is the Hausdorff limit in
\cref{thm:hausdorff}.
For $c^2+c^*$, \cref{alg:adjoint-example} obtains the density and
an exact parametrization of the boundary in \cref{fig:adjoint-example}.
For nonconstant polynomials in free circular variables without adjoints,
\cref{sp:circular-main} also gives absolute continuity and regularity
of the potential. The final example, $c^2+cd+d$, has disk support of
radius $\sqrt3$ (\cref{alg:quadratic-example}).

\subsection{Main ideas}
\paragraph{A second-moment criterion.}
Suppose the $*$-moments of random bounded operators $B_n$ converge in
probability to those of $B$, with the norm bounds of \cref{prop:criterion}.
Then weak convergence in probability of their Brown measures, convergence
in mean of their potentials locally in $L^1$, and convergence of the
expected second moments $\E b(B_n)\to b(B)$ are equivalent. The disk
identity bounds the potential error by the second-moment deficit
in \eqref{eq:deficit}. For GUE polynomials this reduces the weak
convergence problem to
\begin{equation}\label{eq:main-second}
 \E b(A_n)\longrightarrow b(A).
\end{equation}
This route gives a short proof of convergence in probability and in
mean. The almost-sure proof of \cref{sec:almost-sure} does not depend
on it, and gives these conclusions as well (\cref{rem:two-proofs}).

\paragraph{A holomorphic rotation between two independent copies.}
Let $H^{(n)}$ be an independent copy of $G^{(n)}$ and put, for $w\in\C$,
\begin{equation}\label{eq:rotation-intro}
 Z_j^{(n)}(w)=\cos(w)G_j^{(n)}+\sin(w)H_j^{(n)},\qquad
 Z^{(n)}(iy)\overset{d}=G'^{(n)}+\lambda_y C^{(n)},\quad
 \lambda_y=\sqrt{2}\sinh y\quad(y\ge0),
\end{equation}
where $G'$ is a normalized GUE family and $C$ is an independent
standard Ginibre family. Two elementary facts make this family useful.
For real $x$, $Z^{(n)}(x+w)$ and $Z^{(n)}(w)$ have the same law; and
for $w=iy$ with $y\ne0$ a Ginibre component is present. A quantity that
is subharmonic in $w$ and depends only on the law is therefore a convex
function of $y=\Im w$ (\cref{lem:rotation-strip}).

We apply this first to the expected second moment of the eigenvalues,
$y\mapsto\E b(P(Z^{(n)}(iy)))$. It is even and convex, and bounded at
$y=1$ uniformly in $n$, hence Lipschitz at $y=0$ uniformly in $n$
(\cref{thm:convex}). Convergence for each fixed $y>0$, where a Ginibre
component is present, therefore gives \eqref{eq:main-second}. This argument does not require an estimate on the least singular value of $P(G^{(n)})-z$.

For almost-sure convergence, we apply the rotation to
$\frac1n\log|\det(P(Z^{(n)}(w))-z)|$,
which is subharmonic in $w$, with
$y_n=(\log n)^{-L}$. A Cartan-type estimate (\cref{as:cartan}) replaces
convexity, and the invariance under real translations turns an estimate
valid for most real $x$ into a probability bound at $x=0$
(\cref{sec:almost-sure}). Applied to a free semicircular family $(s,t)$
in place of $(G,H)$, the same rotation relates the elliptic and
semicircular cases in the proof of \cref{sp:boundary}.

\paragraph{The Dyson equation and the Fuglede--Kadison determinant.}
To prove \cref{sp:main} and to compute the Brown measure,
we reduce the polynomial to a linear pencil with matrix
coefficients. Hermitization
replaces such a pencil by a selfadjoint block operator, whose
operator-valued Cauchy transform is the unique solution with negative
imaginary part of a quadratic matrix equation, the Dyson equation of
the pencil (\cref{sp:mde}); see \cite{VoiculescuOV,SpeicherOV} for the
general framework. We use an identity expressing the regularized
Fuglede--Kadison determinant of the pencil in terms of the Dyson-equation solution (\cref{sp:action}). The identity also appears in \cite[Theorem 3.1]{Kargin} and \cite[version 1, equation (2.8)]{YinFK}.
Together with the Dyson equation, this identity provides a finite-dimensional procedure for
computing the Brown measure of $P(s)$, even when the measure is not absolutely continuous
with respect to planar Lebesgue measure (\cref{thm:algorithm}).

\paragraph{Reduction to the elliptic case.}
For $y\ne0$, the inputs of the rotation are elliptic:
GUE plus Ginibre for the matrices, and $as+ibt$
with $b\ne0$ for their limits.
For the matrices, conditioning on $G'$ leaves
Ginibre matrices with deterministic means,
for which least-singular-value estimates
are available (\cref{struct:thm:main}).
For the free limits, every nonzero irreducible block of the linearized
pencil satisfies the no-intertwiner condition (\cref{sp:no-intertwiner})
by \cref{sp:elliptic-reduction}. This yields
equality of the support of the Brown measure and the spectrum for polynomials in elliptic elements
(\cref{sp:elliptic-theorem}). At the semicircular endpoint, passage to
$y=0$ shows that $\partial\spec P(s)\subseteq\supp\mu_{P(s)}$
(\cref{sp:boundary}). The harmonicity argument below shows that the
interior of the spectrum also lies in the support.

\paragraph{Harmonicity forces invertibility.}
Outside $\supp\mu_A$ the potential $U_A$ is harmonic. We prove in
\cref{sp:main} that if $U_A$ is harmonic on an open set, then $A-z$
is invertible there. After linearization this becomes a
statement about matrix pencils $I+\sum_jA_j(z)\otimes s_j$ with
holomorphic coefficients (\cref{sp:harmonic-disk}). The proof is based on the behaviour of the solution of the Dyson equation as the regularization parameter tends to zero. 
For a nonzero irreducible block satisfying the no-intertwiner
condition on a disk of harmonicity, noninvertibility at one point forces
the traces of all words in its coefficients to be constant throughout the
connected parameter domain. The block pencil would then be noninvertible
throughout that domain, contradicting invertibility at a known resolvent
point. At the semicircular endpoint the no-intertwiner condition is not
automatic: a block either is similar to a constant tuple, and its pencil
is then invertible throughout the domain, or satisfies the condition on
a smaller disk (\cref{sp:exceptional}).

\subsection{Related literature}
\paragraph{Convergence to Brown measures.}
For complex Gaussian matrices, \cref{thm:main} recovers the circular
and elliptic laws \cite{Ginibre,Mehta,Girko,Bai,TV,NO}, as well as
the product law for independent square complex Ginibre matrices in
\cite[Theorem 8.19]{GKT}. Related product laws for real elliptic matrices and for two
real symmetric Wigner matrices were established in \cite{ORSV}.
These works treat particular models with general entry distributions;
our theorem treats arbitrary polynomials in independent GUE matrices.
Driver--Hall--Ho--Kemp--Nemish--Nikitopoulos--Parraud
\cite[Theorem 1.21]{LimaBean} prove Brown-measure convergence for
matrix random walks with independent bi-unitarily invariant steps,
under their singular-value hypotheses. Random regularization, initiated
by \'Sniady \cite{Sniady} and developed for other nonnormal models in
\cite{GWZ,BPZ}, provides another route to convergence to Brown measures
in particular cases. In our approach, a singular-value estimate is
needed only for an auxiliary ensemble with a Ginibre component, and not
for the GUE polynomial itself.

The complex Gaussian rotation allows us to write the deformed GUE
family as $G'+\lambda C$, where $G'$ and $C$ are independent
GUE and Ginibre families, respectively.
After linearization and conditioning on $G'$, we obtain a
tensorized Ginibre matrix with a deterministic mean.
Han established the corresponding least-singular-value estimate
for the centered case \cite[Theorem 1.5]{Han}.
In \cref{sec:means}, we extend this estimate to means of the form
$\sum_j R_j\otimes D_j$ and obtain bounds with polynomial
dependence on $\lambda^{-1}$ (\cref{struct:thm:main}).
This allows $\lambda$ to tend to zero with $n$, as required
in the proof of almost-sure convergence.

The main random-matrix inputs also include asymptotic freeness \cite{Voi}, bounds near zero for the spectral distributions of polynomials in free semicircular variables and quantitative estimates for the expected spectral distribution functions \cite{BM}, Gaussian concentration and norm bounds \cite{AGZ}, and strong convergence \cite{HT}.

\paragraph{Convergence of spectra.}
Extensions of \cite{HT} include strong asymptotic freeness
with deterministic matrices and Haar unitaries
\cite{MaleStrong,CollinsMaleStrong}; see \cite{vH2026}
for a survey. For nonnormal models, we refer to the results on the support
in the single-ring model by Guionnet and Zeitouni \cite{GZ}
and to the analysis of outliers for holomorphic polynomials
with deterministic matrices in \cite{BBCC}. Brailovskaya, Cook, Kemp, and Parraud
\cite[Theorem 7.2]{BCKP} prove Hausdorff convergence for Brownian
motions on $\mathrm{GL}(n,\C)$ with strongly convergent unitary
initial data. Their general criterion \cite[Proposition 7.3(3)]{BCKP}
combines strong convergence, weak convergence of eigenvalue measures,
and equality of the support of the Brown measure with the spectrum.
After establishing the latter two properties for arbitrary GUE
polynomials, we use their finite-cover argument with quantitative
potential estimates to obtain summable Hausdorff deviation bounds.

\paragraph{The support of the Brown measure and the spectrum.}
Haagerup and Larsen \cite{HL} computed Brown measures
and spectra of $R$-diagonal operators; their formulas also exhibit cases
where the support of the Brown measure is strictly smaller than the spectrum.
Earlier positive results include
circular and elliptic additions under support assumptions, due to
Hall and Ho \cite[Corollary 4.8 and Theorem 4.15]{HH} (see also \cite{HH2022, HZ, ZhongElliptic}), and bounded
freely infinitely divisible $R$-diagonal elements \cite[Corollary 5.3]{KPZ},
with related Brown-measure formulas in \cite{BercoviciZhong}. Miyagawa
\cite{Miyagawa} establishes equality for homogeneous circular
polynomials, computes spectra in the quadratic circular case, and
also treats a nonhomogeneous example \cite[Example 4.5]{Miyagawa}.
Hall--Ho additionally prove
equality for free multiplicative Brownian motion with unitary initial
data \cite[Theorem 3.7]{HH}; earlier Brown-measure calculations and
support inclusions appear in \cite{DHK,HallHoMultiplicative,HallKempSupport}.

In the proof of our invertibility criterion, we construct an $L^2$ inverse using a Wick
expansion and then apply \cite[Theorem 1.1]{Miyagawa}
to show that the inverse is bounded. Related analyses of Dyson equations can be found in
\cite{AEK,CHNR,AK,AKSpectrum}; Hall--Ho use a different
analytic-extension criterion \cite[Theorem 2.6]{HH}.
The algorithm builds on the Brown-measure algorithm of
Belinschi--\'Sniady--Speicher \cite[Sections 3--4]{BSS} and the
linearization and subordination methods of \cite{BMS,HMS}.
For centered matrix-valued semicircular elements, Mai and Speicher
express the Fuglede--Kadison determinant in terms of the capacity of
the covariance map \cite[Theorem 1.2]{MaiSpeicher}.

The analysis of the Dyson equation also shows that the Brown measure of every nonconstant
polynomial in free circular variables without adjoints is absolutely continuous with respect
to two-dimensional Lebesgue measure (\cref{sp:circular-main}).

\subsection{Organization}
\Cref{sec:prelim} collects the preliminaries. \Cref{sec:weak} proves
\cref{thm:main}, using the least-singular-value estimate of
\cref{sec:means}. \Cref{sp:analysis} studies the Dyson equation of a
linear pencil and proves that harmonicity of its potential forces
invertibility; \cref{sp:support} deduces \cref{sp:main,thm:hausdorff},
and \cref{alg:algorithm} gives the algorithm of \cref{thm:algorithm}
and works it out for $c^2+c^*$. \Cref{sp:circular} proves absolute
continuity of the Brown measure for circular polynomials without
adjoints (\cref{sp:circular-main}) and gives an example in which the support of
the Brown measure is a disk.

\section{Technical preliminaries}\label{sec:prelim}
This section collects what the rest of the paper uses.
\Cref{subsec:free-inputs} recalls the Gaussian matrices, their free
limits, and the norm and concentration estimates for the random-matrix
arguments of \cref{sec:weak,sec:means}. \Cref{subsec:shared-rotation}
recalls Fuglede--Kadison determinants and Brown measures, and proves the
subharmonicity and convexity statements on which the rotation argument
rests (\cref{prop:subharm,lem:rotation-strip}); they are used for
matrices in \cref{sec:weak} and for operators in \cref{sp:support}.
\Cref{subsec:linearization} gives the linearization of $P-z$ by an
affine pencil.
For an operator between Hilbert spaces, $\|T\|$ denotes its operator
norm; for a vector, $\|v\|$ denotes its Hilbert-space norm. For a
finite-dimensional operator $T$, write
$\|T\|_{\HS}^2=\Tr(T^*T)$, with unnormalized trace. In a tracial
von Neumann algebra, write $\|T\|_{2,\tau}^2=\tau(T^*T)$.
For a selfadjoint matrix or operator $H$, we write $H>0$ if it is
positive definite and $H\ge0$ if it is positive semidefinite.
A \emph{linear pencil}, or \emph{affine pencil}, is an expression
$K\otimes1+\sum_jA_j\otimes x_j$ with coefficients
$K,A_j\in M_m(\C)$, evaluated at matrices or bounded operators $x_j$.
When matrices themselves are regarded as vectors in a Hilbert--Schmidt
space, we specify that space or use the subscript $\HS$.

For convenience, we list some notation and
where it is defined.

\begin{center}
\small
\begin{tabular}{lll}
\toprule
Notation & Meaning & Defined in\\
\midrule
$\tr_n$, $\tau_m$
  & Normalized traces
  & \eqref{eq:gue-cov}, \cref{subsec:shared-rotation}\\
$\mu_B$, $U_B$
  & Brown measure and logarithmic potential
  & \eqref{eq:matrix-potential}\\
$b(B)$
  & Second moment of $\mu_B$
  & \eqref{eq:b-definition}\\
$\sg(T)$
  & Least singular value
  & \cref{subsec:free-inputs}\\
$K_z$, $B_j$
  & Coefficients of a linear pencil
  & \cref{lem:linearization}\\
$R_\ell$
  & Coefficient matrices
  & \cref{struct:thm:main}\\
$\mathbf A$
  & Tuple $(A_1,\ldots,A_k)$
  & \cref{sp:no-intertwiner}\\
$V_\eta$
  & $\frac12\tau_m\log(T^*T+\eta^2I)$
  & \cref{sp:dyson}\\
$U_\eta$
  & $mV_\eta(L_z(s))$
  & \cref{thm:algorithm,alg:potential}\\
\bottomrule
\end{tabular}
\end{center}

\subsection{Gaussian matrices and free semicircular families}\label{subsec:free-inputs}\label{sec:gaussian}
A normalized GUE matrix $G^{(n)}$ is an $n$-by-$n$ Hermitian matrix whose
diagonal entries are independent centered real Gaussians of variance
$1/n$, and whose strictly upper triangular entries are independent
centered circular complex Gaussians of variance $1/n$, independent of
the diagonal entries. Here a circular complex Gaussian of variance $1/n$
has the form $(g+ih)/\sqrt{2n}$, with $g,h$ independent standard real
Gaussians; the lower triangular entries are determined by Hermitian
symmetry. For an independent family $G_1^{(n)},\ldots,G_k^{(n)}$, the
covariance identity and normalized trace convention are
\begin{equation}\label{eq:gue-cov}
 \E G_i^{(n)}(a,b)G_j^{(n)}(c,d)
   =\frac1n\delta_{ij}\delta_{ad}\delta_{bc},
 \qquad \tr_n=\frac1n\Tr.
\end{equation}
A normalized complex Ginibre matrix of size $n$ has independent entries
$(g_{ab}+ih_{ab})/\sqrt{2n}$, where all $g_{ab},h_{ab}$ are independent
standard real Gaussians. Its entries have mean zero and variance $1/n$,
and satisfy $\E[(C_{ab}^{(n)})^2]=0$.

We work in a finite von Neumann algebra $(\mathcal M,\tau)$ with
a faithful normal tracial state. A standard semicircular variable
has variance one, spectral distribution
$(2\pi)^{-1}\sqrt{4-x^2}\,1_{[-2,2]}(x)\dd x$, and norm two.
The joint $*$-law of a standard free semicircular family is invariant
under real orthogonal changes of variables; see \cite{Voi,AGZ}.
A standard circular variable is $c=(s+it)/\sqrt2$ for a free standard
semicircular pair. More generally, the changes of generators
\[
 C_j^{(n)}=(G_j^{(n)}+iH_j^{(n)})/\sqrt2,\qquad
 c_j=(s_j+it_j)/\sqrt2
\]
identify independent Ginibre matrices with pairs of independent GUE
matrices, and free circular variables with pairs from a jointly free
semicircular family. Taking real and imaginary parts recovers the
original generators, so these substitutions identify $*$-polynomials
in the Ginibre or circular variables with polynomials in the
corresponding selfadjoint variables.

\paragraph{Norm bounds and moment convergence.}
\begin{lemma}[{\cite{AGZ}}]\label[lemma]{lem:norm}
For a normalized GUE matrix $G^{(n)}$ there are absolute $c,C>0$ such that
$\Pp(\|G^{(n)}\|>8)\le Ce^{-cn}$. Moreover
$\sup_n\E\|G^{(n)}\|^p<\infty$ for every fixed $p>0$.
The same assertions hold for a normalized complex Ginibre matrix.
\end{lemma}

For GUE, these bounds follow from the extreme-eigenvalue and Gaussian
concentration estimates in \cite[Theorem 2.6.6 and Corollary 4.4.19]{AGZ}.
These GUE estimates give an exponential tail at every fixed threshold
greater than $2$. For Ginibre use
$C^{(n)}\stackrel{\mathrm{law}}=(G^{(n)}+iH_n)/\sqrt2$, where $G^{(n)},H_n$ are
independent normalized GUE matrices: if both norms are at most
$4\sqrt2$, then $\|C^{(n)}\|\le8$. The triangle inequality also transfers
the uniform norm-moment bounds.
\begin{proposition}[Voiculescu {\cite{Voi}}]\label[proposition]{prop:moments}
For each fixed noncommutative polynomial $Q$ in a GUE family,
\[
 \tr_n Q(G^{(n)})\longrightarrow\tau(Q(s))\quad\text{in probability}.
\]
In particular this holds for all $*$-moments of any fixed polynomial in
that family. Such a polynomial has a fixed operator-norm bound with
probability tending to one and uniformly bounded norm moments of every
fixed order.
\end{proposition}

The moment convergence follows from Voiculescu's
asymptotic freeness theorem \cite{Voi}; see also
\cite[Theorem 5.4.2]{AGZ}.
The norm bounds follow directly from \cref{lem:norm}.

\begin{lemma}
\label[lemma]{lem:gaussian-lipschitz}
On tuples with $\max_j\|X_j\|,\max_j\|Y_j\|\le M$, a fixed
polynomial $P$ satisfies
\begin{equation}\label{eq:polynomial-lipschitz}
 \begin{gathered}
 \|P(X)-P(Y)\|\le C_{P,M}\max_j\|X_j-Y_j\|,\\
 \|P(X)-P(Y)\|_{\HS}\le C_{P,M}
       \Bigl(\sum_j\|X_j-Y_j\|_{\HS}^2\Bigr)^{1/2}.
 \end{gathered}
\end{equation}
In standard real Gaussian coordinates for normalized GUE, a
$1$-Lipschitz operator-norm functional of $P(G)$ has Lipschitz
constant $C_{P,M}/\sqrt n$ on this norm-bounded set. For a selfadjoint
polynomial $Q$ and a real Lipschitz function $\psi$, the statistic
$\tr_n\psi(Q(G))$ has constant $C_{Q,M}\operatorname{Lip}(\psi)/n$.
If a statistic takes values in $[a,b]$, its restriction has an
extension $\widetilde F$ to all coordinates, with the same Lipschitz
constant $L_F$ and range, satisfying
\begin{equation}\label{eq:gaussian-extension}
 \Pp\{|\widetilde F-\E\widetilde F|>t\}
 \le2e^{-t^2/(2L_F^2)},\qquad
 |\E\widetilde F-\E F|\le(b-a)\Pp(\mathcal N^c),
\end{equation}
where $\mathcal N$ is the norm-bounded event and the mean comparison
assumes $F\in[a,b]$ on the whole probability space.
\end{lemma}
\begin{proof}
Telescoping each word gives \eqref{eq:polynomial-lipschitz}.
For Hermitian matrices, Hoffman--Wielandt \cite{HoffmanWielandt} gives
\begin{equation}\label{eq:hoffman-wielandt-trace}
 |\tr_n\psi(X)-\tr_n\psi(Y)|
 \le n^{-1/2}\operatorname{Lip}(\psi)\|X-Y\|_{\HS}.
\end{equation}
In a Hilbert--Schmidt orthonormal Hermitian basis,
\begin{equation}\label{eq:gue-coordinates}
 G^{(n)}=n^{-1/2}\sum_a\xi_aE_a,
 \qquad \xi_a\text{ independent standard real Gaussians}.
\end{equation}
Combining these estimates gives the stated
Lipschitz constants. We extend the restriction
of $F$ to $\mathcal N$ using the McShane
extension theorem and then clip the extension
to $[a,b]$. The extension agrees with $F$ on $\mathcal N$;
Gaussian concentration \cite[Corollary 4.4.19]{AGZ} gives
\eqref{eq:gaussian-extension}.
\end{proof}

\paragraph{Least singular values.}
For a matrix or bounded operator $T$, write
$\sg(T)=\inf_{\|v\|=1}\|Tv\|=\min\spec(T^*T)^{1/2}$.
In a finite von Neumann algebra, $\sg(T)>0$ is equivalent to
invertibility; then $\sg(T)=\|T^{-1}\|^{-1}$.
We also use the elementary inequalities
\begin{equation}\label{eq:singular-calculus}
 |\sg(T)-\sg(S)|\le\|T-S\|,\qquad
 \sg(TS)\ge\sg(T)\sg(S).
\end{equation}
In a finite von Neumann algebra, a one-sided inverse is a two-sided
inverse: the polar decomposition and faithfulness of the finite trace
exclude a proper isometry.

\subsection{Determinants and Brown measures}\label{subsec:shared-rotation}
We write $\dd^2z$ for planar Lebesgue measure and
$\Lap=\partial_x^2+\partial_y^2$ for the planar Laplacian, where
$z=x+iy$. The complex partial derivatives, also called Wirtinger
derivatives, are
\[
 \partial_z=\frac12(\partial_x-i\partial_y),\qquad
 \partial_{\bar z}=\frac12(\partial_x+i\partial_y),\qquad
 \Lap=4\partial_{\bar z}\partial_z.
\]
They act entrywise on matrix-valued functions. For a holomorphic
function, $\partial_{\bar z}f=0$ and $\partial_zf$ is its usual
complex derivative.
For $t>0$, put $\log^+t=\max\{0,\log t\}$ and
$\log^-t=\max\{0,-\log t\}$, with the extended values at $t=0$.
\begin{definition}[Fuglede--Kadison determinant]\label[definition]{def:fk}
For $T\in\mathcal M$ let $|T|=(T^*T)^{1/2}$, and let $\rho_{|T|}$ be
its spectral distribution under $\tau$. Define
\[
 \log\Delta_\tau(T)=\int_{[0,\infty)}\log t\dd\rho_{|T|}(t)
 =\lim_{\epsilon\downarrow0}\frac12\tau\log(T^*T+\epsilon1),
 \qquad \Delta_\tau(T)=\exp(\log\Delta_\tau(T)).
\]
The logarithm may be $-\infty$, and $\exp(-\infty)=0$.
\end{definition}
\begin{theorem}[Brown {\cite{Brown}}; Fuglede--Kadison {\cite{FK}}]\label[theorem]{thm:brown}
For bounded $T\in(\mathcal M,\tau)$ the function
$U_T:\C\to\R\cup\{-\infty\}$, given by
$U_T(z)=\log\Delta_\tau(T-z1)$, is subharmonic and locally integrable.
There is a unique probability measure $\mu_T$, supported in the spectrum
of $T$, such that
\begin{equation}\label{eq:brown-potential}
 U_T(z)=\int_\C\log|z-\lambda|\dd\mu_T(\lambda),\qquad
 \Lap U_T=2\pi\mu_T
\end{equation}
in the distributional sense, with the potential identity also understood
at points where it takes the value $-\infty$. For a normal operator this
is its spectral distribution. In $M_n(\C)$ with $\tr_n$,
\[
 \Delta_{\tr_n}(T)=|\det T|^{1/n},\qquad
 \mu_T=\frac1n\sum_{j=1}^n\delta_{\lambda_j(T)}.
\]
\end{theorem}
These are the foundational results of \cite{Brown,FK}; see also
\cite[Introduction, pp.~19--20]{HS}. We record their consequences
used below:
\begin{equation}\label{eq:brown-support-bound}
 \supp\mu_T\subset\spec T\subset\{|z|\le\|T\|\},\qquad
 b(T):=\int |z|^2\dd\mu_T(z)\le\|T\|^2.
\end{equation}
The measure is determined by the tracial $*$-law. Its potential is
harmonic off $\supp\mu_T$, and determinant differentiation gives
\begin{equation}\label{eq:brown-resolvent}
 2\partial_zU_T(z)=\tau((z-T)^{-1})\qquad(z\notin\spec T).
\end{equation}
For random probability measures, ``weakly in probability'' means
convergence in probability of their integrals against every bounded
continuous test function.

For a matrix $B\in M_d(\C)$, define its normalized eigenvalue measure
and logarithmic potential by
\begin{equation}\label{eq:matrix-potential}
 \mu_B=\frac1d\sum_{j=1}^d\delta_{\lambda_j(B)},\qquad
 U_B(z)=\frac1d\log|\det(B-zI_d)|,
\end{equation}
counting eigenvalues with algebraic multiplicity. On a fixed finite
tracial von Neumann algebra $(\mathcal M,\tau)$, define the second
moment of the Brown measure, $b:\mathcal M\to[0,\infty)$, by
\begin{equation}\label{eq:b-definition}
 b(B)=\int_\C |z|^2\dd\mu_B(z).
\end{equation}
For a finite matrix this is $d^{-1}\sum_j|\lambda_j(B)|^2$, since
$\mu_B$ is the eigenvalue measure for the normalized trace $\tr_d$. In general, $b(B)$
is different from the $*$-moment $\tau(B^*B)$.

For
$h\in C_c^\infty(\C)$,
\begin{equation}\label{eq:green}
 \int h\dd\mu_T=\frac1{2\pi}\int_\C (\Lap h)(z)U_T(z)\dd^2z.
\end{equation}
The same formula holds with $\mu_B$ and $U_B$ for every finite matrix $B$.
For background on logarithmic potentials and Girko's Hermitization,
see \cite[Section 4.1]{BC12}.

\begin{lemma}\label[lemma]{lem:kernel}
If $\nu$ is a probability measure supported in $\{|\lambda|\le M\}$,
and $u_\nu(z)=\int\log|z-\lambda|\dd\nu(\lambda)$, then for each
compact $K\subset\C$,
\begin{equation}\label{eq:kernel-bound}
 \|u_\nu\|_{L^2(K)}\le C_K\bigl(1+\log(1+M)\bigr).
\end{equation}
\end{lemma}
\begin{proof}
The negative part of the logarithmic kernel satisfies
\[
 \sup_{\lambda\in\C}\int_K(\log^-|z-\lambda|)^2\dd^2z
 \le\int_{|w|<1}|\log|w||^2\dd^2w<\infty.
\]
If $K\subset\{|z|\le R\}$ and $|\lambda|\le M$, the positive part
is bounded by $\log(1+R+M)$. Combining the two bounds and applying
Minkowski's inequality to the integral against $\nu$ proves the assertion.
\end{proof}

\begin{lemma}\label[lemma]{lem:disk}
For $R>0$ and $\lambda\in\C$,
\begin{equation}\label{eq:disk-kernel}
 \int_{|z|<R}\log|z-\lambda|\dd^2z=
 \begin{cases}
  \pi R^2\log R-\pi R^2/2+\pi|\lambda|^2/2,&|\lambda|\le R,\\
  \pi R^2\log|\lambda|,&|\lambda|>R.
 \end{cases}
\end{equation}
For a probability measure $\nu$ with finite second moment, put
\begin{equation}\label{eq:disk-TR}
 T_R(\nu)=\int_{|\lambda|>R}
  \left[|\lambda|^2-R^2-2R^2\log(|\lambda|/R)\right]\dd\nu(\lambda).
\end{equation}
Then
\begin{equation}
 \int_{|z|<R}u_\nu(z)\dd^2z
 =\pi R^2\log R-\frac{\pi R^2}{2}
   +\frac\pi2\int|\lambda|^2\dd\nu(\lambda)-\frac\pi2 T_R(\nu).
 \label{eq:disk-tail}
\end{equation}
The bracket in $T_R$ lies between zero and $|\lambda|^2$.
\end{lemma}
\begin{proof}
For $r\ne|\lambda|$, factor out the larger of $r$ and $|\lambda|$
and apply the mean-value property to $\log|1-\zeta|$, harmonic on
$|\zeta|<1$. This gives
\[
 \frac1{2\pi}\int_0^{2\pi}\log|re^{i\theta}-\lambda|\dd\theta
 =\log\max\{r,|\lambda|\}.
\]
Polar integration yields \eqref{eq:disk-kernel}; the single exceptional
radius does not affect the integral. Integrating against $\nu$ gives
\eqref{eq:disk-tail}. Fubini's theorem is justified by the local integrability of the
negative part of the logarithmic kernel and the finite second moment. Finally,
$0\le x^2-1-2\log x\le x^2$ for $x\ge1$, since the middle expression
vanishes at one and has derivative $2x-2/x\ge0$.
\end{proof}

\paragraph{Determinant identities.}
For operator matrices, $\id\otimes\tau$ is the conditional expectation
that applies $\tau$ entrywise, and
$\tau_m=m^{-1}(\Tr\otimes\tau)$ is the normalized trace.
We write $\Delta_m$ for $\Delta_{\tau_m}$, and reserve $\det$ for the usual determinant
of a finite-dimensional complex matrix.

\begin{proposition}\label[proposition]{sp:det-facts}
The Fuglede--Kadison determinant is multiplicative and invariant under
similarity. For an upper block triangular $T\in M_m(\mathcal M)$
with diagonal blocks $T_\ell$ of sizes $m_\ell$,
\[
 \log\Delta_m(T)=\sum_\ell\frac{m_\ell}{m}
                                  \log\Delta_{m_\ell}(T_\ell).
\]
For a bounded operator $B$ in $(\mathcal M,\tau)$ and a scalar polynomial
$q$, $\mu_{q(B)}=q_*\mu_B$, where $q_*\mu_B$ is the push-forward measure
defined by $(q_*\mu_B)(E)=\mu_B(q^{-1}(E))$ for every Borel set
$E\subset\C$. The same identity holds for a rational function whose
poles are outside $\spec B$.
\end{proposition}

\begin{proof}
Multiplicativity for bounded operators, including singular ones, is established in
\cite{FK,HSunbounded}. Similarity invariance follows from multiplicativity. The triangular formula can be seen directly
in the two-block case by the factorization
\[
 \begin{pmatrix}D_1&C\\0&D_2\end{pmatrix}
 =\begin{pmatrix}I&0\\0&D_2\end{pmatrix}
  \begin{pmatrix}I&C\\0&I\end{pmatrix}
  \begin{pmatrix}D_1&0\\0&I\end{pmatrix}.
\]
The middle factor has determinant one: its inverse is unitarily
conjugate to it by $\operatorname{diag}(I,-I)$, so multiplicativity
makes its positive determinant equal to its reciprocal. Induction gives
the general formula, including the extended value $-\infty$.
For a rational function $f=p/q$ with no poles on $\spec B$,
choose $p,q$ coprime, so $q(B)$ is invertible. If $f$ is constant,
the assertion is immediate. Otherwise, for every $z$ the nonzero
polynomial $p-zq$ factors over $\C$, and multiplicativity gives
\[
 \log\Delta_\tau(f(B)-z)
 =\log\Delta_\tau(p(B)-zq(B))-\log\Delta_\tau(q(B))
 =\int\log|f(\lambda)-z|\dd\mu_B(\lambda).
\]
Thus the potential of $f(B)$ is that of $f_*\mu_B$; the
uniqueness in \cref{thm:brown} proves the assertion. Since $\supp\mu_B$ is compact and $f$ is continuous
on $\supp\mu_B$, we also have
\[
\supp\mu_{f(B)}=f(\supp\mu_B).
\]
\end{proof}

\begin{lemma}\label[lemma]{sp:finite-potential}
If $Q$ is a nonconstant scalar polynomial in a free standard
semicircular family $s$, then
\[
 \log\Delta_\tau(Q(s)-z)>-\infty\qquad(z\in\C).
\]
Consequently $\mu_{Q(s)}$ has no atoms. The same conclusion holds for nonconstant
$*$-polynomials in a free standard circular family by \cref{subsec:free-inputs}.
\end{lemma}
\begin{proof}
Write $Q_d$ for the highest homogeneous part of $Q$. A free
semicircular family is algebraically free, so the nonzero polynomial
$H_z=(Q-z)^*(Q-z)$ does not vanish at $s$. Applying
\cite[Theorems 1.1(2) and 5.17]{SS} to $H_z=(Q-z)^*(Q-z)$,
whose highest homogeneous part is $Q_d^*Q_d\ne0$, shows that
its spectral law $\sigma_z$ has no atom at zero
and an integrable logarithm, so
$\log\Delta_\tau(Q(s)-z)=\tfrac12\int\log t\dd\sigma_z(t)>-\infty$.
An atom of $\mu_{Q(s)}$ would contradict \eqref{eq:brown-potential}.
\end{proof}

\paragraph{Subharmonicity.}\label{sec:preliminaries-subharm}
The following results apply to both finite
matrices and bounded operators in a tracial
von Neumann algebra.
A family $D:\Omega\to\mathcal M$, where
$\Omega\subset\C$ is open, is \emph{norm-holomorphic} if it is holomorphic
as a map into the Banach space $(\mathcal M,\|\cdot\|)$.

\begin{lemma}\label[lemma]{lem:levi}
Let $D(w)$ be a norm-holomorphic family in $(\mathcal M,\tau)$, and
$\epsilon>0$. Then
\begin{equation}\label{eq:levi}
 \partial_{\bar w}\partial_w\left[\tfrac12\tau\log(D^*D+\epsilon1)\right]
 =\frac\epsilon2\tau\left[
 (D^*D+\epsilon1)^{-1}D'^*(DD^*+\epsilon1)^{-1}D'\right]\ge0.
\end{equation}
Consequently $w\mapsto\log\Delta_\tau(D(w))$ is subharmonic or
identically $-\infty$ on each connected component.
\end{lemma}
\begin{proof}
For a differentiable positive invertible family $A$, the resolvent
integral for $\log A$ and cyclicity of the trace give
\begin{equation}\label{eq:trace-log-derivative}
 \partial\tau\log A=\tau(A^{-1}\partial A).
\end{equation}
Put $A=D^*D+\epsilon1$ and $J_\epsilon=(DD^*+\epsilon1)^{-1}$.
Differentiating once more gives
\[
 \partial_{\bar w}\partial_w\left[\tfrac12\tau\log A\right]
 =\tfrac12\tau\left[A^{-1}D'^*(1-DA^{-1}D^*)D'\right]
 =\tfrac\epsilon2\|J_\epsilon^{1/2}D'A^{-1/2}\|_{2,\tau}^2,
\]
where $1-DA^{-1}D^*=\epsilon J_\epsilon$. This proves
\eqref{eq:levi} and its positivity.
The regularized logarithms decrease as $\epsilon\downarrow0$; the
decreasing-limit theorem for subharmonic functions proves the final claim
\cite[Chapter 2]{Ransford}. For $D(w)=T-w$, see \cite{Brown,HSunbounded}.
We use $\epsilon=\eta^2$ when the regularization is introduced through
a Hermitized resolvent with regularization parameter $\eta$.
\end{proof}

We use this second-moment subharmonicity for Gaussian rotation and
convexity (\cref{sec:rotation}) and for the subharmonicity of the moments in
\cref{sp:moment-roots}.
\begin{proposition}\label[proposition]{prop:subharm}
If $B(w)$ is norm-holomorphic, then $w\mapsto b(B(w))$ is subharmonic.
\end{proposition}
\begin{proof}
Locally choose $R>\sup_w\|B(w)\|$. By \cref{lem:disk,thm:brown},
\begin{equation}\label{eq:b-disk}
 b(B(w))=\frac2\pi\int_{|z|<R}\log\Delta_\tau(B(w)-z1)\dd^2z
             -2R^2\log R+R^2.
\end{equation}
For $\epsilon>0$ set
\[
 b_\epsilon(w)=\frac1\pi\int_{|z|<R}
   \tau\log\bigl((B(w)-z1)^*(B(w)-z1)+\epsilon1\bigr)\dd^2z
   -2R^2\log R+R^2.
\]
By \cref{lem:levi}, $b_\epsilon$ is subharmonic. Moreover
$b_\epsilon(w)\downarrow b(B(w))$ as $\epsilon\downarrow0$:
for $\|B(w)\|<R$, $|z|<R$, and $0<\epsilon\le R^2$, the
integrand $\tau\log((B(w)-z)^*(B(w)-z)+\epsilon)$ is at most
$\log(5R^2)$, since $\|(B(w)-z)^*(B(w)-z)+\epsilon\|\le(2R)^2+R^2$.
After subtracting this constant, the integrands are nonpositive and
decrease as $\epsilon\downarrow0$, so monotone convergence applies. The limit is finite by
\cref{lem:kernel}, so the decreasing-limit theorem \cite[Chapter 2]{Ransford} proves the assertion.
\end{proof}

\begin{lemma}\label[lemma]{lem:rotation-strip}
Let $G,H$ be independent standard GUE families, or let $s,t$ be a
jointly free standard semicircular family. For real $x$,
\begin{equation}\label{eq:rotation-law}
 (G^{(x)},H^{(x)})=(\cos x\,G+\sin x\,H,-\sin x\,G+\cos x\,H)
 \stackrel{\mathrm d}=(G,H),
\end{equation}
and the same equality holds in joint $*$-law for $(s,t)$. For
$Z(w)=\cos w\,G+\sin w\,H$,
\begin{equation}\label{eq:rotation-height}
 Z(x+w)=\cos w\,G^{(x)}+\sin w\,H^{(x)},\qquad
 Z(x+iy)=\cosh y\,G^{(x)}+i\sinh y\,H^{(x)}.
\end{equation}
A finite subharmonic function on a horizontal strip which depends
only on $\Im w$ is convex as a function of $\Im w$.
\end{lemma}
\begin{proof}
The Gaussian covariance and the semicircular moment formula are
unchanged by real orthogonal transformations, giving
\eqref{eq:rotation-law}, and a direct calculation gives \eqref{eq:rotation-height}.
For the last assertion, let $u(x+iy)=\varphi(y)$ be such a function.
It is locally integrable. For nonnegative test functions $\chi$ on
$\R$ and $\psi$ on the interval of values of $y$, subharmonicity and
$\int\chi''\dd x=0$ give
\[
 0\le\int u(x+iy)\,\Lap\bigl(\chi(x)\psi(y)\bigr)\dd x\dd y
 =\Bigl(\int\chi\dd x\Bigr)\int\varphi(y)\psi''(y)\dd y.
\]
Thus the distributional second derivative of $\varphi$ is nonnegative,
so $\varphi$ agrees almost
everywhere with a convex function.
The area-mean characterization of
subharmonic functions makes this agreement pointwise; see
\cite[Chapter 2]{Ransford}.
\end{proof}

\subsection{Affine linearization}\label{subsec:linearization}
Fix $P\in\C\langle x_1,\ldots,x_k\rangle$.
We use an explicit version of the linearization trick; see
\cite[Section 2]{HT} and \cite[Section 4]{HMS}.
For random nonnormal polynomials, see also \cite{CGH}.
The affine pencil and normalization at the origin used in
\cite[Sections 1.1 and 2]{Han} extend here to any point $c\in\C^k$ with $P(c)\ne z$; we call such a
point a scalar basepoint for $z$. For the combination of linearization and
Hermitization in the computation of Brown measures, see \cite[Section 3]{BSS}.
\begin{lemma}\label[lemma]{lem:linearization}
Let $P\in\C\langle x_1,\ldots,x_k\rangle$. There are an integer
$m\ge1$ and matrices $B_1,\ldots,B_k\in M_m(\C)$, depending only
on $P$. Put $p_0=P(0)$ and $K_z=\operatorname{diag}(p_0-z,I_{m-1})$.
The affine pencil
\begin{equation}\label{eq:pencil}
 L_z(X)=K_z\otimes I_n+\sum_i B_i\otimes X_i
\end{equation}
satisfies, for every $n\ge1$ and $X_1,\ldots,X_k\in M_n(\C)$,
\begin{equation}\label{eq:linear-sv}
 \sg(P(X)-zI_n)\ge\sg(L_z(X)).
\end{equation}
For every scalar tuple $c$ with $P(c)\ne z$, the matrix
$K_{z,c}:=L_z(c)$ is invertible and
\begin{equation}\label{eq:recentered-pencil}
 L_z(X)=K_{z,c}\otimes I_n+\sum_i B_i\otimes(X_i-c_iI_n).
\end{equation}

For bounded operator substitutions the same pencil has the Schur
factorization \eqref{eq:schur-factorization}, where $Q$ is upper
triangular with identity diagonal and the other two factors are
unipotent. In particular its Schur complement is $P(X)-z$.
\end{lemma}
\begin{proof}
If $P$ is affine, take $m=1$ and $L_z(X)=P(X)-zI_n$.
Otherwise write
\[
 P=p_0+P_1+\sum_{w\in\mathcal W}a_w x_w,
 \qquad x_w=x_{i_1}\cdots x_{i_d},\quad d=|w|\ge2,
\]
where $P_1$ is homogeneous of degree one and $\mathcal W$ is the
finite set of words of length at least two with nonzero coefficient.
For each $w$ introduce $d-1$ auxiliary coordinates and set
\[
 Q_w=I_{d-1}-\sum_{a=1}^{d-2}E_{a,a+1}x_{i_{a+1}},\qquad
 u_w=(a_wx_{i_1},0,\ldots,0),\qquad
 v_w=(0,\ldots,0,x_{i_d})^{\mathsf T},
\]
where $E_{a,b}$ are matrix units; for $d=2$ the sum is empty.
Then $u_wQ_w^{-1}v_w=a_wx_w$. In a fixed ordering of $\mathcal W$, put
\[
 Q=\operatorname{diag}_{w\in\mathcal W}Q_w,\qquad
 u=(u_w)_{w\in\mathcal W},\qquad
 v=\operatorname{col}_{w\in\mathcal W}(v_w),\qquad
 m=1+\sum_{w\in\mathcal W}(|w|-1),
\]
where $u,v$ are the concatenated block row and column, and define
\begin{equation}\label{eq:chain-pencil}
 L_z(X)=\begin{pmatrix}
 (p_0-z)I_n+P_1(X)&-u(X)\\v(X)&Q(X)
 \end{pmatrix}.
\end{equation}
Since $Q(X)$ is always invertible, block elimination gives
\begin{equation}\label{eq:schur-factorization}
 L_z=
 \begin{pmatrix}I&-uQ^{-1}\\0&I\end{pmatrix}
 \begin{pmatrix}P(X)-zI_n&0\\0&Q\end{pmatrix}
 \begin{pmatrix}I&0\\Q^{-1}v&I\end{pmatrix}.
\end{equation}
Here $u,v,Q$ are evaluated at $X$, and
$K_z=\operatorname{diag}(p_0-z,I_{m-1})$.
The factorization identifies the Schur complement as $P(X)-zI_n$
and proves the invertibility assertions. In the invertible case,
$(P(X)-zI_n)^{-1}$ is a compression of $L_z(X)^{-1}$, giving
\eqref{eq:linear-sv}; in the singular case both sides are zero.
The affine form gives \eqref{eq:recentered-pencil}.
\end{proof}
\begin{corollary}\label[corollary]{sp:linear-det}
Let $X_1,\ldots,X_k$ be elements of a finite
von Neumann algebra $(\mathcal M,\tau)$, and let
$L_z(X)$ be the pencil constructed in
\cref{lem:linearization}. With the normalized trace $\tau_m$ defined in
\cref{subsec:shared-rotation},
\begin{equation}\label{eq:linear-determinant}
 \log\Delta_\tau(P(X)-z)=m\log\Delta_{\tau_m}(L_z(X)).
\end{equation}
For $z\ne p_0=P(0)$ define $A_j(z)=K_z^{-1}B_j$. Then
\begin{equation}\label{sp:normalized-pencil}
 U_{P(X)}(z)=\log|p_0-z|+
 m\log\Delta_{\tau_m}\left(I+\sum_j A_j(z)\otimes X_j\right).
\end{equation}
The coefficients $A_j$ are rational in $z$, and invertibility of the
operator in parentheses is equivalent to that of $P(X)-z$.
\end{corollary}
\begin{proof}
For $m=1$ the assertions are immediate. Otherwise the unipotent
factors and $Q$ in \eqref{eq:schur-factorization} have determinant one,
so \cref{sp:det-facts} gives the first identity and the factorization
gives the invertibility equivalence. Left multiplication by
$K_z^{-1}\otimes1$ preserves invertibility and gives
\eqref{sp:normalized-pencil}, since $\det K_z=p_0-z$.
\end{proof}

\section{Weak convergence for GUE polynomials}\label{sec:weak}
We give two proofs of weak convergence. The first reduces convergence of
the eigenvalue distribution to convergence of its second moment
(\cref{prop:criterion}). For GUE polynomials, convexity in the rotation
parameter reduces the latter to an auxiliary ensemble with a Ginibre
component (\cref{sec:rotation,sec:interior}). The second applies the
rotation directly to the logarithmic potential and gives all of
\cref{thm:main} (\cref{sec:almost-sure}).

\begin{remark}[Two proofs of weak convergence]\label[remark]{rem:two-proofs}
We include the first proof because it is much shorter and more
conceptual: convergence of the eigenvalue distribution is reduced to
convergence of a single number. Moreover, the second-moment criterion
applies to general sequences of random operators. For GUE polynomials, \cref{sec:almost-sure} proves all of
\cref{thm:main} directly, including almost-sure convergence and the
deviation bounds. That proof uses the rotation, the singular-value
counting estimates of \cref{lem:counts} at $y=0$, and the
least-singular-value bounds of \cref{sec:means}. It bypasses both the
second-moment criterion and the convexity argument.
\end{remark}

\subsection{A criterion using the second moment of the Brown measure}\label{sec:criterion}
The criterion below applies to random bounded operators; in the matrix
case no Gaussian assumption is imposed. Recall from \eqref{eq:b-definition} that this second moment is
\[
 b(B)=\int_\C|z|^2\dd\mu_B(z),\qquad
 b(B_n)=\frac1{d_n}\sum_{j=1}^{d_n}|\lambda_j(B_n)|^2
 \quad\text{for }B_n\in M_{d_n}(\C).
\]
The quantity to compare with $b(B)$ is the expected value
$\E b(B_n)$, since the operators may be random.
The regularized-logarithm argument underlying upper semicontinuity of
the determinant is standard; see \cite[Lemma 2]{Sniady}.
The passage from potentials to measures uses the distributional
Laplacian, as in \cite[Section 4.1]{BC12}. We combine these facts
with the disk identity to obtain the following criterion.

\begin{theorem}[Second-moment criterion]\label[theorem]{prop:criterion}
Let $B_n$ be strongly measurable random variables with values in
the operator-norm Banach spaces of finite von Neumann algebras
$(\mathcal M_n,\tau^{(n)})$, with faithful normal normalized traces. Let
$B\in(\mathcal M,\tau)$ be a fixed bounded operator. Assume:
\begin{enumerate}[label=(\roman*),leftmargin=*]
\item all tracial $*$-moments $\tau^{(n)}(w(B_n,B_n^*))$ converge
in probability to $\tau(w(B,B^*))$;
\item for some fixed $R_0>0$, $\Pp(\|B_n\|>R_0)\to0$;
\item $\sup_n\E\|B_n\|^p<\infty$ for some $p>2$.
\end{enumerate}
The following are equivalent:
\begin{enumerate}[label=(\alph*),leftmargin=*]
\item $\E b(B_n)\to b(B)$;
\item $\E\int_K|U_{B_n}-U_B|\dd^2z\to0$ for every compact $K$;
\item $\mu_{B_n}\Rightarrow\mu_B$ weakly in probability.
\end{enumerate}
\end{theorem}
The operator norm is finite for each realization; no deterministic
uniform bound is assumed. Finite random matrices are the case
$(\mathcal M_n,\tau^{(n)})=(M_{d_n}(\C),\tr_{d_n})$.
\begin{proof}
\textit{(a)$\Longrightarrow$(b).}
Fix a compact $K$ and choose $R>\max\{R_0,\|B\|\}$ with
$K\subset D_R:=\{|z|<R\}$. By (ii)--(iii) and H\"older's inequality,
\begin{equation}\label{eq:criterion-tail}
 \E[\|B_n\|^2 1_{\{\|B_n\|>R\}}]
 \le (\E\|B_n\|^p)^{2/p}
       \Pp(\|B_n\|>R)^{1-2/p}\longrightarrow0.
\end{equation}
For $\epsilon>0$ set
\[
 U_{B_n,\epsilon}(z)=\tfrac12\tau^{(n)}\log((B_n-z)^*(B_n-z)+\epsilon1),
 \quad U_{B,\epsilon}(z)=\tfrac12\tau\log((B-z)^*(B-z)+\epsilon1).
\]
On $\{\|B_n\|\le R\}$, polynomial approximation to $\log(x+\epsilon)$
on $[0,(R+|z|)^2]$ and (i) give convergence in probability
of $U_{B_n,\epsilon}(z)$ to $U_{B,\epsilon}(z)$. The bounds
\[
 \tfrac12\log\epsilon\le U_{B_n,\epsilon}(z)
 \le\tfrac12\log((\|B_n\|+|z|)^2+\epsilon)
\]
and (iii) give uniform integrability for each $z$ and a uniform
bound on the expectations of their absolute values for $z\in D_R$.
Uniform integrability upgrades convergence in probability of
$U_{B_n,\epsilon}(z)$ to convergence in mean for each $z$.
Dominated convergence in $z$ therefore gives
\begin{equation}\label{eq:reg-L1}
 \E\int_{D_R}|U_{B_n,\epsilon}-U_{B,\epsilon}|\dd^2z\longrightarrow0.
\end{equation}
For the unregularized potentials, \cref{lem:kernel} and (iii)
justify the use of Fubini's theorem. Since $U_{B_n}\le U_{B_n,\epsilon}$,
\[
 (U_{B_n}-U_B)_+\le|U_{B_n,\epsilon}-U_{B,\epsilon}|
                              +(U_{B,\epsilon}-U_B).
\]
Thus \eqref{eq:reg-L1} gives
\[
 \limsup_n\E\int_{D_R}(U_{B_n}-U_B)_+\dd^2z
 \le\int_{D_R}(U_{B,\epsilon}-U_B)\dd^2z.
\]
The right side decreases to zero as $\epsilon\downarrow0$, dominated
for $\epsilon\le1$ by the integrable function $U_{B,1}-U_B$.
Consequently,
\begin{equation}\label{eq:positive-defect}
 \E\int_{D_R}(U_{B_n}-U_B)_+\dd^2z\longrightarrow0.
\end{equation}
The disk identity \cref{lem:disk} gives
\begin{equation}\label{eq:criterion-disk}
 \E\int_{D_R}(U_{B_n}-U_B)\dd^2z
 =\frac\pi2\bigl(\E b(B_n)-b(B)-t_{n,R}\bigr),\qquad
 t_{n,R}=\E T_R(\mu_{B_n}).
\end{equation}
By \cref{lem:disk},
$0\le t_{n,R}\le\E[\|B_n\|^2 1_{\{\|B_n\|>R\}}]$.
The right side tends to zero by \eqref{eq:criterion-tail}. Hence
\begin{equation}\label{eq:negative-defect}
 \E\int_{D_R}|U_{B_n}-U_B|\dd^2z
 =2\E\int_{D_R}(U_{B_n}-U_B)_+\dd^2z
   +\frac\pi2\bigl(b(B)-\E b(B_n)+t_{n,R}\bigr).
\end{equation}
By (a) and \eqref{eq:positive-defect}, this tends to zero, proving (b).

\textit{(b)$\Longrightarrow$(c).} Choose $R>\max\{R_0,\|B\|\}$.
For $\varphi\in C_c^\infty(\C)$, Green's identity \eqref{eq:green} gives
\[
 \E\left|\int\varphi\dd(\mu_{B_n}-\mu_B)\right|
 \le\frac{\|\Lap\varphi\|_\infty}{2\pi}
       \E\int_{\supp\varphi}|U_{B_n}-U_B|\dd^2z\longrightarrow0.
\]
On $\{\|B_n\|\le R\}$ both measures are supported in
$\{|z|\le R\}$. Uniform approximation there by smooth compactly
supported functions, Markov's inequality, and (ii) extend the
convergence to every bounded continuous test function.

\textit{(c)$\Longrightarrow$(a).} Choose $R>\max\{R_0,\|B\|\}$.
For the bounded continuous function $h_R(z)=\min\{|z|^2,R^2\}$,
(c) gives convergence of its integrals in probability and in expectation.
Since $\int h_R\dd\mu_B=b(B)$ and
\[
 0\le\E b(B_n)-\E\int h_R\dd\mu_{B_n}
 \le\E[\|B_n\|^2 1_{\{\|B_n\|>R\}}]\longrightarrow0
\]
by \eqref{eq:criterion-tail}, (a) follows.
\end{proof}

\begin{remark}[Sharp potential-deficit bound]\label[remark]{ex:jordan}
Under (i)--(iii), even without (a), \eqref{eq:positive-defect} and
\eqref{eq:criterion-disk} give $\limsup_n\E b(B_n)\le b(B)$.
Moreover, \eqref{eq:negative-defect} implies, for every compact $K$,
\begin{equation}\label{eq:deficit}
 \limsup_n\E\int_K|U_{B_n}-U_B|\dd^2z
 \le\frac\pi2\bigl(b(B)-\liminf_n\E b(B_n)\bigr).
\end{equation}
The constant is sharp. For the Jordan shifts $J_n$ with ones on the
first upper diagonal, $\|J_n\|\le1$ and the normalized $*$-moments
converge to those of a Haar unitary $u$: a word of length $\ell$
has the bilateral-shift action except on at most $2\ell$ boundary
basis vectors. This is the standard discontinuity example of
\cite[Section 2.6]{Sniady}. Here $b(J_n)=0$, $b(u)=1$,
$U_{J_n}(z)=\log|z|$, and $U_u(z)=\log^+|z|$, so
\[
 \int_{|z|\le1}|U_{J_n}-U_u|\dd^2z
 =2\pi\int_0^1(-\log r)r\dd r=\pi/2.
\]
\end{remark}

\subsection{Complex Gaussian rotation}\label{sec:rotation}
Let $G,H$ be independent normalized GUE families, and let $s,t$ be a joint
free standard semicircular family of size $2k$. Recall the rotation
\eqref{eq:rotation-intro} and define its free counterpart, for $w\in\C$, by
\begin{equation}\label{eq:rotation-family}
 Z_i^{(n)}(w)=\cos(w)G_i+\sin(w)H_i,
 \qquad Z_i^{(\infty)}(w)=\cos(w)s_i+\sin(w)t_i.
\end{equation}
Both are polynomials in $\cos w$ and $\sin w$ with fixed coefficients,
so they are holomorphic in $w$ on all of $\C$, with values in $M_n(\C)$
and in $\mathcal M$ respectively; in particular they are
norm-holomorphic in the sense of \cref{sec:preliminaries-subharm}.
Define $f_n,f:\R\to[0,\infty)$ by
\begin{equation}\label{eq:fn-def}
 f_n(y)=\E b(P(Z^{(n)}(iy))),\qquad
 f(y)=b(P(Z^{(\infty)}(iy))).
\end{equation}

\Needspace{6\baselineskip}
\begin{theorem}[Convexity along the rotation]\label[theorem]{thm:convex}
The functions $f_n,f$ are finite, even, and convex on $\R$.
For a constant $C_P$ independent of $n$,
\begin{equation}\label{eq:boundary}
 0\le f_n(y)-f_n(0)\le C_P|y|,\qquad
 0\le f(y)-f(0)\le C_P|y|\qquad(|y|\le1).
\end{equation}
For $0<y\le1$,
\begin{equation}\label{eq:boundary-error}
 |f_n(0)-f(0)|\le |f_n(y)-f(y)|+2C_Py.
\end{equation}
In particular, convergence at every fixed positive $y$ implies
$f_n(0)\to f(0)$.
\end{theorem}

\begin{proof}
\textit{Step 1: subharmonicity in the complex rotation parameter.}
The constant case is immediate, so assume $\deg P\ge1$.
Define $g_n,g:\C\to[0,\infty)$ by
\begin{equation}\label{eq:gn-def}
 g_n(w)=\E b(P(Z^{(n)}(w))),\qquad
 g(w)=b(P(Z^{(\infty)}(w))).
\end{equation}
By \cref{prop:subharm}, $g$ is subharmonic. For the matrix model,
\cref{prop:subharm} gives pathwise subharmonicity of
$v_n(w)=b(P(Z^{(n)}(w)))$.
On each compact set $K\subset\C$, the bound $b(B)\le\|B\|^2$ gives
\[
 \sup_{w\in K}v_n(w)
 \le C_{P,K}\left(1+\sum_i(\|G_i^{(n)}\|+\|H_i^{(n)}\|)\right)^{2\deg P}.
\]
For fixed $n$, $v_n$ is continuous because the eigenvalue multiset
depends continuously on the matrix entries. The right side is
integrable by \cref{lem:norm}, so dominated convergence makes $g_n$
continuous. Since $v_n\ge0$, we may take the expectation of the
submean inequality for $v_n$ and exchange the expectation with the
angular integral by Tonelli's theorem. For every $w_0\in\C$ and $r>0$,
this gives
\[
 g_n(w_0)\le\frac1{2\pi}\int_0^{2\pi}
                  g_n(w_0+re^{i\theta})\dd\theta.
\]
Thus $g_n$ is subharmonic.

\textit{Step 2: horizontal invariance gives convexity.}
By \cref{lem:rotation-strip}, the law of $Z^{(n)}(x+iy)$ and the
$*$-law of $Z^{(\infty)}(x+iy)$ depend only on $y$. Brown measure
is determined by the $*$-law, so
\[
 g_n(x+iy)=f_n(y),\qquad g(x+iy)=f(y).
\]
The convexity assertion of \cref{lem:rotation-strip} applies to both functions.
Changing the sign of the second input family proves that $f_n,f$
are even.

\textit{Step 3: a uniform estimate at the semicircular endpoint.}
By \cref{lem:norm} for the matrices, and boundedness of the free inputs,
$\max\{\sup_n f_n(1),f(1)\}\le C_P$.
Since these functions are even and convex, for $|y|\le1$,
\[
 f_n(0)\le\tfrac12\bigl(f_n(y)+f_n(-y)\bigr)=f_n(y),
 \qquad
 f_n(y)\le(1-|y|)f_n(0)+|y|f_n(1).
\]
Consequently, using $f_n(0)\ge0$,
\[
 0\le f_n(y)-f_n(0)
 \le |y|\bigl(f_n(1)-f_n(0)\bigr)\le C_P|y|.
\]
The same argument applies to $f$. Fix $0<y\le1$. Then
\[
 \begin{aligned}
 |f_n(0)-f(0)|
 &\le |f_n(0)-f_n(y)|+|f_n(y)-f(y)|+|f(y)-f(0)|\\
 &\le |f_n(y)-f(y)|+2C_Py.
 \end{aligned}
\]
If $f_n(y)\to f(y)$ for every fixed $y>0$, this yields
$\limsup_n|f_n(0)-f(0)|\le2C_Py$. Now let $y\downarrow0$.
\end{proof}

\subsection{Convergence in the auxiliary ensemble and removal of the deformation}\label{sec:interior}
Fix $y>0$ and write
\begin{equation}\label{eq:Any-def}
 A_{n,y}=P(\cosh(y)G^{(n)}+i\sinh(y)H^{(n)}),\qquad
 A_y=P(\cosh(y)s+i\sinh(y)t).
\end{equation}
We first prove convergence at each positive $y$, and then use
\cref{thm:convex} to reach $y=0$.
We may suppose that $P$ is nonconstant,
since a constant polynomial has the same deterministic value in every model.

For fixed $y>0$, the input tuples satisfy
\begin{equation}\label{eq:elliptic}
 \cosh(y)G+i\sinh(y)H\stackrel{\mathrm d}=G'+\lambda_yC,
 \qquad \lambda_y=\sqrt2\sinh(y),
\end{equation}
where $G'_i$ are independent normalized GUE, $C_i$ are independent
normalized complex Ginibre, and the two families are independent.
Indeed, realize $C_i=(G''_i+iH''_i)/\sqrt2$ using two further independent
GUE families, independent of $G'$. Then $G'_i+\sinh(y)G''_i$ has the
law of $\cosh(y)G_i$, since $1+\sinh^2(y)=\cosh^2(y)$, and is
independent of $H''_i$ and of the other input indices. The same argument
with free semicircular families gives the corresponding equality of
joint $*$-laws for semicircular and circular inputs. The independent
Ginibre decomposition and conditioning device also occur in
\cite[proofs of Propositions 5.1--5.2]{BCKP}.

After linearization and conditioning on $G'$, \eqref{eq:elliptic}
gives a tensorized Ginibre matrix with a deterministic mean in the
coefficient span. \Cref{struct:thm:main}, proved in \cref{sec:means},
controls these means uniformly. Its case $\lambda=1$ suffices here;
its explicit dependence on $\lambda$ will be used in \cref{sec:almost-sure}.

The estimate below needs a scalar basepoint $c\in\C^k$ with
$P(c)\ne z$, as in \cref{lem:linearization}. Define the scalar
exceptional set
\begin{equation}\label{eq:exceptional-set}
 \mathcal E_P=\{z\in\C:P(c)=z\text{ for every }c\in\C^k\}.
\end{equation}
This set is empty when the scalar evaluation $c\mapsto P(c)$
is nonconstant, and equals $\{a\}$ when $P(c)=a$ for every scalar
tuple $c$. The latter can occur even for a nonconstant polynomial,
as for a commutator. For $z\notin\mathcal E_P$ we can choose
a scalar basepoint $c$ with $P(c)\ne z$.
Let $m$ be the size of the scalar coefficient matrices in the pencil
$L_z(X)\in M_m(\C)\otimes M_n(\C)$ of \cref{lem:linearization},
and let $a(m)$ be the exponent in \cref{struct:thm:main}.

\begin{proposition}\label[proposition]{prop:interior-lsv}
For fixed $y>0$ and $z\notin\mathcal E_P$, there are $C,c>0$
and $n_0$ such that, for every $t>0$ and $n\ge n_0$,
\begin{equation}\label{eq:interior-tail}
 \Pp\{\sg(A_{n,y}-zI_n)\le t\}
 \le C\bigl(n^{a(m)}t^{2/m}+e^{-cn}\bigr).
\end{equation}
The constants $C,c$ and the threshold $n_0$ may depend on $P,y,z$
and a chosen scalar basepoint $c_0\in\C^k$ with $P(c_0)\ne z$.
For any $q>0$, one can choose $M$, depending only on $m,q$, such that
\begin{equation}\label{eq:interior-lsv}
 \Pp\{\sg(A_{n,y}-zI_n)<n^{-M}\}=O(n^{-q}).
\end{equation}
\end{proposition}
The estimate at fixed $y$ suffices for convergence in mean. The
dependence on $\lambda$ in \cref{struct:thm:main} supplies the
uniformity as $y\downarrow0$ used in \cref{sec:almost-sure}.
\begin{proof}
Choose $c_0\in\C^k$ with $P(c_0)\ne z$ and realize the input
as $G'+\lambda_yC$ using \eqref{eq:elliptic}. Put
\begin{equation}\label{eq:aux-sigma}
 \mathcal F'_n=\sigma(G'_1,\ldots,G'_k),\qquad
 \mathcal N_{n,\mathrm{aux}}=\{\max_i\|G'_i\|\le8\}.
\end{equation}
By \cref{lem:norm}, $\Pp(\mathcal N_{n,\mathrm{aux}}^c)\le Ce^{-cn}$.
Multiplying \eqref{eq:recentered-pencil} on the left gives
\begin{gather*}
 T_{n,y}:=(K_{z,c_0}^{-1}\otimes I_n)L_z(G'+\lambda_yC)
     =I+\sum_i\lambda_yR_i\otimes(C_i+D_i),\\
 R_i=K_{z,c_0}^{-1}B_i,\qquad
 D_i=(G'_i-c_{0,i}I_n)/\lambda_y.
\end{gather*}
Conditionally on $\mathcal F'_n$, the $C_i$ remain independent
normalized Ginibre matrices and the $D_i$ are deterministic.
On $\mathcal N_{n,\mathrm{aux}}$ they satisfy
$\max_i\|D_i\|\le \omega:=\max_i(8+|c_{0,i}|)/\lambda_y$.
Apply the case $\lambda=1$ to the fixed tuple $(\lambda_yR_i)$.
The uniformity over the means in \cref{struct:thm:main} gives,
almost surely on $\mathcal N_{n,\mathrm{aux}}$,
\[
 \Pp\{\sg(T_{n,y})\le u\mid\mathcal F'_n\}
 \le C\bigl(n^{a(m)}u^{2/m}+e^{-cn}\bigr)\qquad(u>0).
\]
By \cref{lem:linearization}, in this realization,
\begin{equation}\label{eq:linear-lsv-transfer}
 \sg(A_{n,y}-zI_n)
 \ge\sg(L_z(G'+\lambda_yC))
 \ge\sg(K_{z,c_0})\sg(T_{n,y}).
\end{equation}
Substitute $u=t/\sg(K_{z,c_0})$ in the conditional estimate,
average over $G'$, and add $\Pp(\mathcal N_{n,\mathrm{aux}}^c)$. This proves
\eqref{eq:interior-tail}. Any $M>m(a(m)+q)/2$ gives
\eqref{eq:interior-lsv}.
\end{proof}

We next control the number of small singular values.
For fixed $y,z$, define the selfadjoint polynomial
\begin{equation}\label{eq:Qyz}
 Q_{y,z}(x,x')=
 \big(P(\cosh(y)x+i\sinh(y)x')-z\big)^*
 \big(P(\cosh(y)x+i\sinh(y)x')-z\big)
\end{equation}
in $2k$ selfadjoint indeterminates. The squared singular values of
$A_{n,y}-zI_n$ are the eigenvalues of $Q_{y,z}(G,H)$.
Denote their empirical law by
$\rho_{n,y,z}$, and denote the spectral law of $(A_y-z)^*(A_y-z)$ by
$\rho_{y,z}$. These are measures on $[0,\infty)$, and \cref{def:fk} gives
\begin{equation}\label{eq:singular-log-laws}
 U_{A_y}(z)=\frac12\int\log x\dd\rho_{y,z}(x),\qquad
 U_{A_{n,y}}(z)=\frac12\int\log x\dd\rho_{n,y,z}(x).
\end{equation}

\begin{lemma}\label[lemma]{lem:counts}
For fixed $y\ge0,z\in\C$ there are $C,\alpha,\beta>0$ such that
\begin{equation}\label{eq:counts}
 \rho_{y,z}([0,x])\le Cx^\alpha,\qquad
 \E\rho_{n,y,z}([0,x])\le Cx^\alpha+Cn^{-\beta}
 \quad(0<x<1).
\end{equation}
The averaged distribution functions also satisfy
\begin{equation}\label{eq:counts-kolmogorov}
 \sup_{t\in\R}|\E\rho_{n,y,z}((-\infty,t])-
                       \rho_{y,z}((-\infty,t])|\le Cn^{-\beta}.
\end{equation}
In particular, $\int\log x\dd\rho_{y,z}(x)$ is finite, including
at $y=0$, and $\rho_{n,y,z}$ converges weakly in probability to $\rho_{y,z}$.
\end{lemma}
\begin{proof}
Recall that $P$ is nonconstant, and let $P_d$ be its top homogeneous
part, with $d=\deg P\ge1$. After setting $x'=0$, the degree-$2d$
part of $Q_{y,z}(x,0)$ is $(\cosh y)^{2d}P_d^*P_d\ne0$, since the
free polynomial algebra has no zero divisors: in the
degree-lexicographic order on words, the leading word of a product is
the product of the leading words. Thus $Q_{y,z}$ is
a nonconstant selfadjoint polynomial.

The conjugate variables of the semicircular family, in the sense of
Voiculescu, are $\xi_j=s_j$, and they are Lipschitz in the sense
recalled in \cite{BM}: their free difference quotients
$\partial_is_j=\delta_{ij}1\otimes1$ are bounded.
Hence \cite[Theorem 1.1(ii)]{BM} applies to this
nonconstant selfadjoint polynomial. Its matrix evaluation is a
selfadjoint polynomial in independent GUE matrices, so the averaged
law estimate \cite[Corollary 6.7]{BM} applies as well. Thus its spectral laws satisfy,
for some $C,\alpha,\beta>0$,
\begin{gather*}
 \rho_{y,z}([0,x])\le Cx^\alpha\quad(0<x<1),\\
 \sup_{t\in\R}\left|
 \E\rho_{n,y,z}(( -\infty,t])-\rho_{y,z}(( -\infty,t])
 \right|\le Cn^{-\beta}.
\end{gather*}
Since both laws are supported on $[0,\infty)$, these estimates give
\eqref{eq:counts}.

By \cref{prop:moments}, both laws are supported in a common bounded
interval with probability tending to one, and their moments converge
in probability. Polynomial approximation on that interval therefore
gives weak convergence in probability.

Finally, for $0<x<\delta<1$,
\begin{equation}\label{eq:log-tail-identity}
 -\log x=-\log\delta+\int_x^\delta t^{-1}\dd t.
\end{equation}
Writing $F(t)=\rho_{y,z}([0,t])$, Tonelli's theorem and \eqref{eq:counts} give
\begin{equation}\label{eq:free-log-tail}
 \begin{aligned}
 \int_{(0,\delta)}|\log x|\dd\rho_{y,z}(x)
 &\le |\log\delta|F(\delta)+\int_0^\delta\frac{F(t)}t\dd t\\
 &\le C\delta^\alpha\bigl(|\log\delta|+\alpha^{-1}\bigr).
 \end{aligned}
\end{equation}
There is no atom at zero by \eqref{eq:counts}, and bounded support
controls the positive logarithm. This proves finiteness.
\end{proof}

\begin{proposition}\label[proposition]{prop:interior}
For every fixed $y>0$ and $z\notin\mathcal E_P$,
\begin{equation}\label{eq:interior-pointwise}
 U_{A_{n,y}}(z)\longrightarrow U_{A_y}(z)
 \quad\text{in probability}.
\end{equation}
For every compact $K\subset\C$,
\begin{equation}\label{eq:interior-local}
 \E\int_K|U_{A_{n,y}}-U_{A_y}|\dd^2z\longrightarrow0.
\end{equation}
Consequently $\mu_{A_{n,y}}\Rightarrow\mu_{A_y}$ weakly in probability and
\begin{equation}\label{eq:interior-second}
 \E b(A_{n,y})\longrightarrow b(A_y).
\end{equation}
\end{proposition}
\begin{proof}
\textit{Step 1: pointwise convergence of the potentials.}
Fix $y>0$ and $z\notin\mathcal E_P$, and suppress $y,z$ from
$\rho_n,\rho$. All constants in this step may depend on $y,z$.
Set $F_n(x)=\rho_n([0,x])$. By \cref{prop:interior-lsv}, for a fixed
large $M$ the event
$\mathcal S_n=\{\sg(A_{n,y}-zI_n)\ge n^{-M}\}$ has probability
tending to one. On this event, all squared singular values are at
least $n^{-2M}$. For fixed $0<\delta<1$ and sufficiently large $n$,
\cref{eq:log-tail-identity,eq:counts-kolmogorov} and Tonelli's theorem give
\begin{equation}\label{eq:log-tail-good}
\begin{aligned}
 &\E\left[1_{\mathcal S_n}
       \int_{(0,\delta)}|\log x|\dd\rho_n(x)\right]\\
 &\qquad\le |\log\delta|\E F_n(\delta)
      +\int_{n^{-2M}}^\delta\frac{\E F_n(x)}x\dd x\\
 &\qquad\le C\delta^\alpha\bigl(|\log\delta|+\alpha^{-1}\bigr)
       +Cn^{-\beta}\bigl(|\log\delta|+2M\log n\bigr).
\end{aligned}
\end{equation}

By \cref{prop:moments}, choose $R>1$ such that
$\supp\rho\subset[0,R]$ and the event
$\mathcal B_n=\{\|Q_{y,z}(G,H)\|\le R\}$ has probability tending to one.
For $0<\delta<1$,
the bounded continuous function
$\ell_{\delta,R}(x)=\log(\min\{R,\max\{\delta,x\}\})$ satisfies
\[
 \int\ell_{\delta,R}\dd(\rho_n-\rho)\longrightarrow0
 \quad\text{in probability}
\]
by \cref{lem:counts}. On $\mathcal S_n\cap\mathcal B_n$,
\eqref{eq:singular-log-laws} and $0\le\ell_{\delta,R}(x)-\log x
\le |\log x|$ for $0<x<\delta$ give
\begin{equation}\label{eq:interior-truncation}
 2|U_{A_{n,y}}(z)-U_{A_y}(z)|
 \le\left|\int\ell_{\delta,R}\dd(\rho_n-\rho)\right|
   +\int_{(0,\delta)}|\log x|\dd\rho_n(x)
   +\int_{(0,\delta)}|\log x|\dd\rho(x).
\end{equation}
For any $\epsilon>0$, choose $\delta$ so that the last term is
at most $\epsilon/2$, using \eqref{eq:free-log-tail}. Markov's
inequality gives
\begin{align*}
 \Pp\{|U_{A_{n,y}}(z)-U_{A_y}(z)|>\epsilon\}
 &\le\Pp(\mathcal S_n^c)+\Pp(\mathcal B_n^c)
   +\Pp\left\{\left|\int\ell_{\delta,R}\dd(\rho_n-\rho)\right|
           >\epsilon\right\}\\
 &\quad+\frac2\epsilon\E\left[1_{\mathcal S_n}
                    \int_{(0,\delta)}|\log x|\dd\rho_n(x)\right].
\end{align*}
At fixed $\delta$, \cref{eq:interior-lsv,eq:log-tail-good} therefore yield
\[
 \limsup_{n\to\infty}
 \Pp\{|U_{A_{n,y}}(z)-U_{A_y}(z)|>\epsilon\}
 \le\frac{2C}{\epsilon}\delta^\alpha
                         (|\log\delta|+\alpha^{-1}).
\]
Letting $\delta\downarrow0$ proves \eqref{eq:interior-pointwise}.

\textit{Step 2: integrate in $z$.}
The exceptional set $\mathcal E_P$ is empty or a singleton by
\eqref{eq:exceptional-set}. Thus on each compact set $K$, Step 1 applies for almost every $z$. Dominated convergence of the deviation
probabilities gives convergence in measure on the product of $K$
with the probability space. By \cref{lem:kernel,prop:moments},
\[
 \sup_n\E\|U_{A_{n,y}}\|_{L^2(K)}^2
 \le C_K\sup_n\E(1+\log(1+\|A_{n,y}\|))^2<\infty .
\]
This gives uniform integrability on the product space and proves
\eqref{eq:interior-local}. The moment and norm hypotheses of
\cref{prop:criterion} follow from \cref{prop:moments} applied to
$(G,H)$. The criterion now gives weak convergence in probability
and \eqref{eq:interior-second}.
\end{proof}

\begin{proof}[Proof of the convergence-in-mean assertion in \cref{thm:main}]
By \cref{prop:interior}, $\E b(A_{n,y})\to b(A_y)$ at every fixed
$y>0$. Convexity and the uniform bound in \cref{thm:convex} give the
same convergence at $y=0$, namely \eqref{eq:main-second}. Applying
\cref{prop:criterion} at $y=0$ proves \eqref{eq:main-potential}
and weak convergence in probability. The constant case is immediate.
\end{proof}

\begin{corollary}\label[corollary]{cor:ginibre-tail}
For independent normalized complex Ginibre matrices $C_i^{(n)}$
and a fixed polynomial $P$ without adjoints, for fixed $z\notin\mathcal E_P$,
\[
 \Pp\{\sg(P(C^{(n)})-zI_n)\le t\}
 \le C\bigl(n^{a(m)}t^{2/m}+e^{-cn}\bigr),
\]
for all $t>0$ and $n\ge n_0$, where $m$ is the linearization size;
$C,c>0$ and $n_0$ may depend on $P,z$, and the exponents depend
only on $m$.
\end{corollary}
\begin{proof}
Choose $c_0$ with $P(c_0)\ne z$, normalize the
pencil at $c_0$, and apply \cref{struct:thm:main,lem:linearization}
with coefficients $K_{z,c_0}^{-1}B_i$ and means $D_i=-(c_0)_iI_n$.
\end{proof}

\begin{remark}
When the scalar evaluation of $P$ is nonconstant, \cref{cor:ginibre-tail}
treats every $z$, including $P(0)$. This answers the Gaussian case of
\cite[Question 1.6]{Han} under this assumption.
\end{remark}

\begin{remark}[The exceptional value]\label[remark]{rem:exception}
Suppose that $P$ is nonconstant but $P(c)=a$ for every scalar tuple
$c$. Equivalently, $P-a$ lies in the two-sided ideal generated by the
commutators $x_ix_j-x_jx_i$. Then
the singular-value estimates of \cref{prop:interior-lsv,cor:ginibre-tail}
leave $z=a$ untreated. This single point does not affect the planar
integration in \cref{prop:interior}, and \cref{thm:main} holds
without exception. Only the pointwise statements at $z=a$ are
affected, and one inequality survives there. At every $z$,
the regularized potential
$U_{n,\epsilon}(z)=\tfrac12\tr_n\log((A_n-z)^*(A_n-z)+\epsilon)$
satisfies $U_{A_n}(z)\le U_{n,\epsilon}(z)$. At fixed $\epsilon>0$,
moment convergence and the uniform norm moments give
$\E U_{n,\epsilon}(z)\to U_{A,\epsilon}(z)$: use polynomial
approximation on a norm-bounded event, and uniform integrability
of $1+\log(1+\|A_n\|)$ on its complement. Letting
$\epsilon\downarrow0$ gives
$\limsup_n\E U_{A_n}(z)\le U_A(z)$, allowing $-\infty$ expectations.
The argument supplies no matching lower bound at the exceptional value.
\end{remark}

\subsection{Quantitative estimates for logarithmic potentials and almost sure convergence}
\label{sec:almost-sure}
In \cref{sec:rotation,sec:interior} we used the rotation
$Z^{(n)}(w)$ of \eqref{eq:rotation-family} to control the expected second
moment at a fixed height $w=iy$. Here we apply the rotation to the
logarithmic potential at a height $y_n\to0$.

Throughout this subsection $P$ is nonconstant, $A_n=P(G^{(n)})$,
$A=P(s)$, and we write $U_n=U_{A_n}$ and $U=U_A$. Let $D(w,r)$
denote the open disk with center $w$ and radius $r$.
\Cref{as:pointwise} bounds the deviations of $U_n(z)$ at fixed $z$.
Its proof studies the function
$h_n(w)=\frac1n\log|\det(P(Z^{(n)}(w))-z)|$ of the rotation parameter,
defined in \eqref{as:rotation-potential}; it is subharmonic, and
$h_n(0)=U_n(z)$. There are three steps. First, for fixed $z\notin\mathcal E_P$ and
$y_n=(\log n)^{-L}$, an upper bound for $h_n$ on a disk of radius
$4y_n$ around $iy_n$ and a lower bound at its center hold with
probability $1-O(n^{-q})$ (\cref{as:anchor}). The upper bound comes
from concentration for a fixed selfadjoint polynomial at $w=0$ and an
operator-norm comparison; the lower bound uses the least-singular-value
estimate and the small-singular-value count at $iy_n$ (\cref{as:log-lsv}), where a Ginibre component is
present. Second, a Cartan-type estimate (\cref{as:cartan}) transfers
this information to most points $x$ of the real interval
$[-y_n/2,y_n/2]$. Third, $h_n(x)$ has the same law as $h_n(0)$ for
every real $x$, so an average over that interval gives a probability
bound at $x=0$. A separate integration argument (\cref{as:integrating})
then gives local $L^1$ bounds.

\begin{theorem}[Concentration of the potential at a fixed point]\label[theorem]{as:pointwise}
For every fixed $z\in\C$ there is $C_{P,z}>0$ such that, for every $q>0$,
\begin{equation}\label{as:upper-tail}
 \Pp\{U_n(z)>U(z)+C_{P,z}(\log n)^{-6}\}=O_{P,z,q}(n^{-q}).
\end{equation}
For every fixed $z\notin\mathcal E_P$ and $q>0$ there is $K_{P,z,q}>0$
such that
\begin{equation}\label{as:lower-tail}
 \Pp\{U_n(z)<U(z)-K_{P,z,q}(\log n)^{-4}\}=O_{P,z,q}(n^{-q}).
\end{equation}
In particular, for such $z,q$, enlarging $K_{P,z,q}$ if necessary gives
\begin{equation}\label{as:pointwise-rate}
 \Pp\{|U_n(z)-U(z)|>K_{P,z,q}(\log n)^{-4}\}=O_{P,z,q}(n^{-q}).
\end{equation}
Consequently $U_n(z)\to U(z)$ almost surely for each fixed
$z\notin\mathcal E_P$, and for every fixed $t,q>0$,
\begin{equation}\label{as:pointwise-bound}
 \Pp\{|U_n(z)-U(z)|>t\}=O_{P,z,t,q}(n^{-q}).
\end{equation}
\end{theorem}

The structured-means theorem \cref{struct:thm:main} allows
$\lambda$ to shrink with $n$. Its losses are polynomial in
$1/\lambda$, which is sufficient for the choice $y_n=(\log n)^{-L}$ in
\cref{as:log-lsv,as:anchor}.

\begin{corollary}
\label[corollary]{as:log-lsv}
Fix $z\notin\mathcal E_P$ and $L>0$. Put
$y_n=(\log n)^{-L}$ for $n\ge3$ and
\[
 A_{n,y_n}=P(\cosh(y_n)G^{(n)}+i\sinh(y_n)H^{(n)}).
\]
For every $q>0$ there exists a fixed $M>0$ such that
\begin{equation}\label{as:log-lsv-bound}
 \Pp\{\sg(A_{n,y_n}-z)<n^{-M}\}=O_{P,z,L,q}(n^{-q}).
\end{equation}
\end{corollary}
\begin{proof}
Choose $c\in\C^k$ with $P(c)\ne z$, and use the
recentered linearization of \cref{lem:linearization}. By
\eqref{eq:rotation-intro}, $A_{n,y}$ has the law of $P(G'+\lambda_y C)$,
where $G'$ and $C$ are independent. For this representative,
\eqref{eq:recentered-pencil} gives the normalized pencil
\[
 (K_{z,c}^{-1}\otimes I_n)L_z(G'+\lambda_y C)
 = I+\sum_j\lambda_y R_j\otimes(C_j+D_j),\qquad
 R_j=K_{z,c}^{-1}B_j,\quad D_j=(G'_j-c_jI_n)/\lambda_y .
\]
Conditionally on $G'$, the means $D_j$ are deterministic and satisfy
$\max_j\|D_j\|\le\omega_0/\lambda_y$ on $\{\max_j\|G'_j\|\le8\}$,
where $\omega_0=8+\max_j|c_j|$. Apply \cref{struct:thm:main} to the
displayed pencil on this event. Equations \eqref{eq:linear-sv} and
\eqref{eq:singular-calculus} bound $\sg(P(G'+\lambda_yC)-z)$ from
below by $\sg(K_{z,c})$ times the least singular value of the normalized
pencil. Integrating the conditional bound and applying \cref{lem:norm}
on the complement gives, for $0<\lambda_y\le1$ and
$n\ge C_*\lambda_y^{-D_m}$,
\[
 \Pp\{\sg(A_{n,y}-z)\le t\}
 \le C\lambda_y^{-2}n^{a(m)}t^{2/m}+Ce^{-cn}.
\]
All constants now depend on the fixed $P,z,c$, not on $y$.
Here $a(m)$ is as in \cref{struct:thm:main}.
At $y=y_n$, $\lambda_{y_n}\asymp(\log n)^{-L}$, so the threshold
holds eventually and $\lambda_{y_n}^{-2}=O((\log n)^{2L})$.
Choose $M$ so that $2M/m>a(m)+q+2$ and set $t=n^{-M}$.
\end{proof}

Fix any $z\in\C$; the restriction $z\notin\mathcal E_P$ is needed only
for the lower bound \eqref{as:anchor-bound}. Let
\begin{equation}\label{as:fixed-polynomial}
 Q_n=(P(G^{(n)})-z)^*(P(G^{(n)})-z),\qquad
 Q=(P(s)-z)^*(P(s)-z),
\end{equation}
and let $\rho_n,\rho$ be their spectral laws, with $\rho_n$ random.
By \cref{lem:counts} at $y=0$, these laws satisfy
\eqref{eq:counts} and \eqref{eq:counts-kolmogorov} with constants
$C,\alpha,\beta>0$ depending only on $P,z$, and $U(z)=\tfrac12\int\log x\dd\rho(x)$ is finite.
For the rotation $Z^{(n)}(w)$ in \eqref{eq:rotation-family}, put
\begin{equation}\label{as:rotation-potential}
 T_n(w)=P(Z^{(n)}(w))-zI_n,\qquad Q_n(w)=T_n(w)^*T_n(w),\qquad
 h_n(w)=\frac1n\log|\det T_n(w)|.
\end{equation}
The determinant is jointly continuous in its matrix entries and $w$;
its logarithm, with the convention $\log0=-\infty$, and hence $h_n$,
is jointly measurable with values in $[-\infty,\infty)$.
Thus $Q_n(0)=Q_n$, $h_n(0)=U_n(z)$ and
$T_n(iy)=A_{n,y}-zI_n$.
By \cref{lem:norm}, we have
\begin{equation}\label{as:norm-event}
 \Omega_n=\{\max_j(\|G_j^{(n)}\|,\|H_j^{(n)}\|)\le8\},
 \qquad \Pp(\Omega_n^c)\le Ce^{-cn}.
\end{equation}
The telescoping bounds \eqref{eq:polynomial-lipschitz} show that there are
fixed $R>1$ and $C_0$ such that
\begin{equation}\label{as:perturb}
 \|Q\|<R,\quad
 \sup_{|w|\le1}\|Q_n(w)\|<R,\quad
 \|Q_n(w)-Q_n(0)\|\le C_0|w|
 \quad\text{on }\Omega_n,\quad |w|\le1 .
\end{equation}
The constants $R,C_0$ depend only on $P,z$.

For $0<\delta<1$, let $\ell_{\delta,R}$ be the logarithmic cutoff used
in \eqref{eq:interior-truncation} and set
\begin{equation}\label{as:floor-statistic}
 F_{n,\delta}=\tfrac12\tr_n\ell_{\delta,R}(Q_n).
\end{equation}
Also let $\chi_\delta$ equal one on $(-\infty,2\delta]$, decrease
linearly to zero on $[2\delta,4\delta]$, and vanish thereafter, and put
\begin{equation}\label{as:count-statistic}
 N_{n,\delta}=\tr_n\chi_\delta(Q_n).
\end{equation}

\begin{lemma}
\label[lemma]{as:fixed-concentration}
For all $0<\delta<1/4$,
\begin{align}
 0\le \tfrac12\int\ell_{\delta,R}\dd\rho-U(z)
                  &\le C\delta^\alpha,\label{as:free-floor}\\
 \left|\E F_{n,\delta}
       -\tfrac12\int\ell_{\delta,R}\dd\rho\right|
        &\le C n^{-\beta}(1+|\log\delta|),\label{as:mean-floor}\\
 \E N_{n,\delta}&\le C\delta^\alpha+Cn^{-\beta}.
                                                    \label{as:mean-count}
\end{align}
For $S=F_{n,\delta}$ or $N_{n,\delta}$, the restriction of $S$ to
$\Omega_n$ is $C/(n\delta)$-Lipschitz in the standard Gaussian
coordinates. It agrees there with a bounded function $\widetilde S$
on the whole coordinate space having the same Lipschitz bound. Moreover,
\[
 |\E\widetilde S-\E S|\le C(1+|\log\delta|)e^{-cn},\qquad
 \Pp\{|\widetilde S-\E\widetilde S|>u\}
 \le 2\exp(-cn^2\delta^2u^2)\quad(u>0).
\]
\end{lemma}
\begin{proof}
To prove \eqref{as:free-floor}, use $\|Q\|<R$, Tonelli's theorem,
and \eqref{eq:counts} at $y=0$ to obtain
\[
 \frac12\int_{(0,\delta)}\log\frac{\delta}{x}\dd\rho(x)
 =\frac12\int_0^\delta\frac{\rho([0,t])}{t}\dd t
 \le\frac{C}{2\alpha}\delta^\alpha.
\]
There is no atom at zero. Since $\ell_{\delta,R}$ is bounded
and monotone with total variation $\log(R/\delta)$, integration against
the difference of distribution functions in \eqref{eq:counts-kolmogorov} gives
\eqref{as:mean-floor}. The same argument applies to the monotone
cutoff $\chi_\delta$, whose total variation is one and whose free
integral is at most $\rho([0,4\delta])$. This gives \eqref{as:mean-count}.

Both scalar cutoffs are $\delta^{-1}$-Lipschitz. Hence
\cref{lem:gaussian-lipschitz} gives the stated Lipschitz bounds on
$\Omega_n$. The two statistics are defined globally and their ranges
have length at most $C(1+|\log\delta|)$. Applying
\eqref{eq:gaussian-extension} to their restrictions to $\Omega_n$ and
using \eqref{as:norm-event} proves the claims about $\widetilde S$.
\end{proof}

\begin{proposition}
\label[proposition]{as:anchor}
For each fixed $z\in\C$ there are constants $L,C_{P,z}>0$ such that,
with $y_n=(\log n)^{-L}$, the estimate
\begin{equation}\label{as:anchor-upper}
 \sup_{w\in D(iy_n,4y_n)} h_n(w)\le U(z)+C_{P,z}(\log n)^{-6}
\end{equation}
holds outside an event of probability $O_{P,z,q}(n^{-q})$ for every $q>0$.
If $z\notin\mathcal E_P$, then for every $q>0$ there is a constant
$C_{P,z,q}>0$ such that, with probability $1-O_{P,z,q}(n^{-q})$,
\eqref{as:anchor-upper} holds and
\begin{equation}\label{as:anchor-bound}
 h_n(iy_n)\ge U(z)-C_{P,z,q}(\log n)^{-5}.
\end{equation}
\end{proposition}
\begin{proof}
Choose an integer $J$ with $J\alpha\ge6$ and set
\begin{equation}\label{as:scales}
 \delta_n=(\log n)^{-J},\qquad
 y_n=(\log n)^{-(J+6)},\qquad
 v_n=(\log n)^{-6}.
\end{equation}
We take $L=J+6$. Since $\delta_n^\alpha\le v_n$ and
$|\log\delta_n|=J\log\log n$, \cref{as:free-floor,as:mean-floor} give
\[
 |\E F_{n,\delta_n}-U(z)|\le C\delta_n^\alpha+
 Cn^{-\beta}(1+J\log\log n)=O(v_n).
\]
Likewise \eqref{as:mean-count} gives
$\E N_{n,\delta_n}\le C(\delta_n^\alpha+n^{-\beta})=O(v_n)$.
In \cref{as:fixed-concentration} the change of mean from extension is
$o(v_n)$. Taking $u=v_n$ there, and using \eqref{as:norm-event},
shows that for a sufficiently large $C$ the event
\begin{equation}\label{as:two-good}
 \mathcal H_n=\Omega_n\cap
 \{|F_{n,\delta_n}-U(z)|\le C v_n\}\cap\{N_{n,\delta_n}\le C v_n\}
\end{equation}
has complement of probability at most
$C\exp(-cn^2/(\log n)^{2J+12})+Ce^{-cn}$.

For $|w|\le5y_n$, Weyl's inequality, the bound
$\operatorname{Lip}(\ell_{\delta_n,R})\le1/\delta_n$, and
\eqref{as:perturb} give
\begin{equation}\label{as:floor-comparison}
 \left|\tfrac12\tr_n\ell_{\delta_n,R}(Q_n(w))
                  -F_{n,\delta_n}\right|
 \le C y_n/\delta_n=Cv_n
 \quad\text{on }\Omega_n.
\end{equation}
On $\mathcal H_n$ all eigenvalues of $Q_n(w)$ are below $R$.
Truncating the logarithm from below increases it, so
\cref{as:two-good,as:floor-comparison} give
\begin{equation}\label{as:disk-upper}
 h_n(w)=\tfrac12\tr_n\log Q_n(w)
       \le \tfrac12\tr_n\ell_{\delta_n,R}(Q_n(w))
       \le U(z)+Cv_n.
\end{equation}
Every $w\in D(iy_n,4y_n)$ satisfies $|w|\le5y_n$, proving
\eqref{as:anchor-upper} for every fixed $z$.

Now fix $z\notin\mathcal E_P$ and $q>0$. Choose $M$ as in
\cref{as:log-lsv} and set
\[
 \mathcal G_n(q)=\mathcal H_n\cap\{\sg(T_n(iy_n))\ge n^{-M}\}.
\]
Since $T_n(iy_n)=A_{n,y_n}-zI_n$, \eqref{as:log-lsv-bound} and
the union bound give $\Pp(\mathcal G_n(q)^c)=O_{P,z,q}(n^{-q})$.
Write $\rho_{n,y_n,z}$ for the spectral law of $Q_n(iy_n)$.
Since $C_0y_n\le\delta_n$ eventually, \eqref{as:perturb}
and Weyl's inequality give, on $\mathcal H_n$,
\[
 \rho_{n,y_n,z}([0,\delta_n])
 \le\rho_n([0,2\delta_n])
 \le N_{n,\delta_n}\le Cv_n .
\]
On $\mathcal G_n(q)$ every eigenvalue of $Q_n(iy_n)$ is also at least
$n^{-2M}$. Hence
\begin{align}
 0&\le \tfrac12\tr_n\ell_{\delta_n,R}(Q_n(iy_n))-h_n(iy_n)
       \nonumber\\
  &=\tfrac12\int_{(0,\delta_n)}
                     \log\frac{\delta_n}{x}\dd\rho_{n,y_n,z}(x)
       \nonumber\\
  &\le M\log n\,\rho_{n,y_n,z}([0,\delta_n])
       \le C M(\log n)^{-5}.\label{as:log-tail}
\end{align}
Together with \eqref{as:two-good} and
\eqref{as:floor-comparison}, this yields
\[
 h_n(iy_n)\ge U(z)-C(1+M)(\log n)^{-5}.
\]
Since $M$ is fixed for each $q$, this proves \eqref{as:anchor-bound}.
\end{proof}

The next estimate is of Cartan type. We include
its potential-theoretic proof, based on the Riesz decomposition,
Harnack's inequality, and exponential integrability of logarithmic
kernels; see \cite[Chapters 1--3]{Ransford} for these tools and
\cite[Lecture 11]{LevinEntire} for the classical Cartan estimate.
\begin{lemma}\label[lemma]{as:cartan}
There are absolute constants $c,C>0$ such that the following holds.
If $u\le0$ is subharmonic on $D(i,4)$ and $u(i)\ge-\epsilon$ for some
$\epsilon>0$, then
\begin{equation}\label{as:cartan-bound}
 \bigl|\{x\in[-1/2,1/2]:u(x)<-t\}\bigr|
       \le C\exp(-ct/\epsilon)\qquad(t>0).
\end{equation}
Values $-\infty$ at other points are allowed. 
\end{lemma}
\begin{proof}
First suppose $u$ is smooth and set
$\dd\nu=(2\pi)^{-1}\Lap u\dd^2z$. Jensen's formula at $i$
with radius $3$ gives
\[
 \int_{D(i,3)}\log\frac3{|\zeta-i|}\dd\nu(\zeta)
 =\frac1{2\pi}\int_0^{2\pi}u(i+3e^{i\theta})\dd\theta-u(i)
 \le\epsilon.
\]
Consequently $m=\nu(D(i,2))\le\epsilon/\log(3/2)$.
On $D(i,2)$ write $u=H-V$, where $H\le0$ is the harmonic
extension of the boundary values and
\[
 V(z)=\int_{D(i,2)}G(z,\zeta)\dd\nu(\zeta)\ge0,\qquad
 G(z,\zeta)=\log\left|
    \frac{4-(z-i)\overline{(\zeta-i)}}{2(z-\zeta)}\right|.
\]
Since $H(i)=u(i)+V(i)\ge-\epsilon$, Harnack's inequality
\cite[Section 1.3]{Ransford} for $-H$ gives $-H(x)\le C_0\epsilon$ on $[-1/2,1/2]$.
There $G(x,\zeta)\le\log(4/|x-\zeta|)$.
When $m>0$, Jensen's inequality yields
\[
 e^{V(x)/(2m)}
 \le\frac1m\int_{D(i,2)}
                \left(\frac4{|x-\zeta|}\right)^{1/2}\dd\nu(\zeta).
\]
As $|x-\zeta|\ge|x-\Re\zeta|$, the integral of $|x-\zeta|^{-1/2}$
over $[-1/2,1/2]$ is uniformly bounded in $\zeta$. Tonelli's theorem
therefore gives $\int_{-1/2}^{1/2}e^{V(x)/(2m)}\dd x\le C$.
For $t\ge2C_0\epsilon$, the inequality $u(x)<-t$ implies $V(x)>t/2$.
Markov's inequality gives a bound $Ce^{-t/(4m)}$, which implies
\eqref{as:cartan-bound} since $m\le\epsilon/\log(3/2)$. When $m=0$,
$u=H$ and Harnack's bound suffices. For $t<2C_0\epsilon$ the interval
has length one, so increasing $C$ gives the same estimate.

For a general subharmonic function $u$, the finite value at $i$
excludes the identically $-\infty$ case. Fix a nonnegative smooth radial
function $\varphi$, supported in $D(0,1)$ and of integral one, and set
$\varphi_r(\xi)=r^{-2}\varphi(\xi/r)$ for $0<r<1$. Put
\[
 u_r(\zeta)=\int u(\zeta-\xi)\varphi_r(\xi)\dd^2\xi,
 \qquad \zeta\in D(i,4-r).
\]
Each $u_r$ is smooth, subharmonic and nonpositive. Circular means of
$u$ are nondecreasing in the radius, so $u_r$ decreases as
$r\downarrow0$. The submean inequality and upper semicontinuity give
the pointwise limit $u_r\downarrow u$, including where $u=-\infty$.
In particular
$u_r(i)\ge u(i)\ge-\epsilon$. The proof in the smooth case used the
function only on a neighborhood of $\overline{D(i,3)}$, so it applies
to $u_r$:
\[
 \bigl|\{x\in[-1/2,1/2]:u_r(x)<-t\}\bigr|\le C\exp(-ct/\epsilon).
\]
Since $u_r$ decreases to $u$, these sets increase to
$\{x\in[-1/2,1/2]:u(x)<-t\}$ as $r\downarrow0$, and the bound passes
to the limit.
\end{proof}

\begin{proof}[Proof of \cref{as:pointwise}]
For each fixed $z\in\C$, $0\in D(iy_n,4y_n)$ and $h_n(0)=U_n(z)$,
so \eqref{as:anchor-upper} gives \eqref{as:upper-tail}.

For the lower tail, fix $z\notin\mathcal E_P$ and $q>0$. By
\cref{as:anchor}, with probability $1-O_{P,z,q}(n^{-q-1})$ both
\eqref{as:anchor-upper} and \eqref{as:anchor-bound} hold with $q+1$
in place of $q$. On this event $w\mapsto\det T_n(w)$ is holomorphic
and nonzero at $iy_n$. Thus the function
\[
 u(\zeta)=h_n(y_n\zeta)-U(z)-C_{P,z}(\log n)^{-6}
\]
is subharmonic and nonpositive on $D(i,4)$, and
$u(i)\ge-D_{P,z,q}(\log n)^{-5}$ for a suitable constant $D_{P,z,q}$.
By \cref{as:cartan}, for every $t>0$,
\begin{equation}\label{as:cartan-applied}
 \frac1{y_n}\left|\{x\in[-y_n/2,y_n/2]:h_n(x)<U(z)-t\}\right|
 \le C\exp\bigl(-c_{P,z,q}\,t(\log n)^5\bigr)
\end{equation}
on the same event.

For every real $x$, $\cos(x)G+\sin(x)H$ has the same law as $G$, so
$h_n(x)$ has the same law as $U_n(z)$. By Tonelli's theorem and
\eqref{as:cartan-applied}, bounding the relative length by one on the
exceptional event, we obtain
\begin{equation}\label{as:tonelli-rotation}
 \begin{aligned}
 \Pp\{U_n(z)<U(z)-t\}
 &=\frac1{y_n}\E\int_{-y_n/2}^{y_n/2}
       \mathbf1_{\{h_n(x)<U(z)-t\}}\dd x\\
 &\le C\exp\bigl(-c_{P,z,q}\,t(\log n)^5\bigr)+O_{P,z,q}(n^{-q-1}).
 \end{aligned}
\end{equation}
Take $t=K_{P,z,q}(\log n)^{-4}$ with $c_{P,z,q}K_{P,z,q}\ge q+1$
to obtain \eqref{as:lower-tail}. Combining the two tails proves
\eqref{as:pointwise-rate}. The threshold there tends to zero, giving
\eqref{as:pointwise-bound} for every fixed $t>0$. Taking $q=2$ in
\eqref{as:pointwise-rate} and applying Borel--Cantelli gives almost-sure
convergence at the fixed point $z$.
\end{proof}

We next integrate the probability bounds, allowing
their constants to depend on $z$.

\begin{lemma}\label[lemma]{as:integrating}
Let $D\subset\C$ be measurable with finite area. Let $V_n$ be
jointly measurable random functions and $V$ a measurable deterministic
function, all finite almost everywhere on $D$.
Suppose for almost every $z\in D$, every $\delta,q>0$,
\[
 \Pp\{|V_n(z)-V(z)|>\delta\}=O_{z,\delta,q}(n^{-q}).
\]
Suppose also that on events $\mathcal N_n$ with
$\Pp(\mathcal N_n^c)\le Ce^{-cn}$,
$\|V_n-V\|_{L^2(D)}\le M_2$ for a fixed $M_2$.
Then for every $t,q>0$,
\[
 \Pp\{\|V_n-V\|_{L^1(D)}>t\}=O_{t,q}(n^{-q}).
\]
\end{lemma}
\begin{proof}
The case $\operatorname{Area}(D)=0$ is immediate. Otherwise fix
$\delta=t/(4\operatorname{Area}(D))$, and write
$p_n(z)=\Pp\{|V_n(z)-V(z)|>\delta\}$.
For integers $M\ge1$, let
\[
 D_M=\{z\in D:p_n(z)\le M n^{-q}\text{ for every }n\ge M\}.
\]
These measurable sets increase to a full-measure subset of $D$.
Choose a fixed $M$ with $M_2\operatorname{Area}(D\setminus D_M)^{1/2}\le t/4$.
Let
$B_n=\{z\in D_M:|V_n(z)-V(z)|>\delta\}$.
On $\mathcal N_n$,
\[
 \|V_n-V\|_{L^1(D)}
 \le t/4+\delta\operatorname{Area}(D)+M_2\operatorname{Area}(B_n)^{1/2}
 \le t/2+M_2\operatorname{Area}(B_n)^{1/2}.
\]
For $n\ge M$, Tonelli's theorem gives
$\E\operatorname{Area}(B_n)\le M\operatorname{Area}(D)n^{-q}$.
Markov's inequality now gives
\[
 \Pp\{\|V_n-V\|_{L^1(D)}>t\}
 \le Ce^{-cn}+\frac{4M_2^2M\operatorname{Area}(D)}{t^2}\,n^{-q}.
\]
The assertion is also immediate if $M_2=0$.
\end{proof}

\begin{corollary}
\label[corollary]{as:integrated}
For every fixed polynomial $P$, compact $K\subset\C$, and $t,q>0$,
\begin{align}
 \Pp\{\|U_{A_n}-U_A\|_{L^1(K)}>t\}
   &=O_{P,K,t,q}(n^{-q}),\label{as:integrated-bound}\\
 \Pp\{|b(A_n)-b(A)|>t\}
   &=O_{P,t,q}(n^{-q}).\label{as:second-bound}
\end{align}
\end{corollary}
\begin{proof}
For constant $P$ the measures are identical and both errors vanish.
Suppose $P$ is nonconstant. On
$\mathcal N_n=\{\max_j\|G_j^{(n)}\|\le8\}$, all eigenvalues of
$A_n$ and the support of $\mu_A$ lie in a fixed disk.
By \cref{lem:norm}, $\Pp(\mathcal N_n^c)\le Ce^{-cn}$. On
$\mathcal N_n$, \cref{lem:kernel} bounds $\|U_{A_n}-U_A\|_{L^2(K)}$
by a deterministic constant. The fixed-threshold estimate
\eqref{as:pointwise-bound} holds for almost every $z$, since
$\mathcal E_P$ has area zero. Thus \cref{as:integrating} gives
\eqref{as:integrated-bound}.

Choose $R$ larger than the common support bound on $\mathcal N_n$.
The disk identity \eqref{eq:b-disk} gives, on $\mathcal N_n$,
\[
 |b(A_n)-b(A)|\le\frac2\pi
       \|U_{A_n}-U_A\|_{L^1(\{|z|\le R\})}.
\]
Equation \eqref{as:second-bound} follows from
\eqref{as:integrated-bound} and the exponential norm tail.
\end{proof}

\begin{proof}[Completion of the proof of \cref{thm:main}]
Apply \cref{as:integrated} with $q=2$ to disks of positive integer
radius and positive rational tolerances. The Borel--Cantelli lemma gives one
event of probability one on which $U_{A_n}\to U_A$ locally in $L^1$.
This proves \eqref{eq:main-potential-as}.
By \cref{lem:norm} and the Borel--Cantelli lemma, the GUE norms are at most
eight eventually on another probability-one event. On the intersection,
the supports are eventually contained in a fixed disk. The identity \eqref{eq:green} gives
convergence against every compactly supported smooth test function;
uniform approximation on the fixed disk then gives convergence
against every bounded continuous function, proving \eqref{eq:main-weak}.
To obtain convergence in
mean without the second-moment criterion, put
$X_n=\|U_{A_n}-U_A\|_{L^1(K)}$. By \cref{lem:kernel},
\[
 X_n\le C_K(1+\log(1+\|A_n\|)),\qquad \sup_n\E X_n^2<\infty
\]
by \cref{prop:moments}. Thus $X_n$ is uniformly integrable and
$\E X_n\to0$, proving \eqref{eq:main-potential}.
The deviation assertions in \cref{thm:main} are precisely
\cref{as:integrated}. The Borel--Cantelli lemma
applied to \eqref{as:second-bound} with $q=2$ also gives $b(A_n)\to b(A)$
almost surely. All assertions of \cref{thm:main} are now proved.
\end{proof}

\section{Least singular value estimates for tensorized Ginibre matrices with structured means}\label{sec:means}
This section proves the least-singular-value estimate used in
\cref{sec:interior,sec:almost-sure}. After conditioning on its GUE
part, the linearized matrix there is a tensorized Ginibre matrix with
a deterministic mean; \cref{struct:thm:main} bounds its least singular
value from below, uniformly over bounded means and with polynomial
dependence on the size $\lambda$ of the Ginibre component.
The coefficient size $m$ is fixed and only $n$ grows.

A coefficient family on $\C^m$ is \emph{irreducible} if it has
no common invariant subspace other than zero and the whole space.
A \emph{composition series} is a flag of common invariant subspaces
whose successive quotients are irreducible for the unital coefficient
algebra. Such a flag exists by finite-dimensional induction and has
length at most $m$.

\begin{theorem}[Least singular value with structured means]
\label[theorem]{struct:thm:main}
Fix $R_1,\ldots,R_k\in M_m(\C)$ and $\omega_0\ge0$. Let $s_1$ be the
number of nonzero quotient actions in a composition series of the
coefficient family. There are an exponent $a=a(m)>0$ and constants
$C_*,C,c>0$, depending only on the base tuple and $\omega_0$, such that
for $0<\lambda\le1$, $n\ge C_*\lambda^{-D_m}$,
$D_m=(2m^2)^{m^2}$, and every deterministic tuple with
$\max_j\|D_j\|\le \omega_0/\lambda$, the matrix
\[
 T_{\lambda,D}=I_m\otimes I_n+
          \sum_j\lambda R_j\otimes(C_j+D_j)
\]
with independent normalized complex Ginibre $C_j$ satisfies, if $s_1>0$,
\begin{equation}\label{eq:struct-1.1}
 \Pp\{\sg(T_{\lambda,D})\le t\}
 \le C\lambda^{-2}n^a t^{2/s_1}+Ce^{-cn}\qquad(t>0).
\end{equation}
In particular $2/m$ is always an admissible exponent, and a nonzero
irreducible tuple has exponent $2$. If $s_1=0$, there is $t_0>0$
such that the bound $Ce^{-cn}$ holds for $0<t<t_0$, uniformly in
$\lambda$. When $\lambda=1$ we write
$T_D=T_{1,D}$. The exponent $a(m)$ is specified in
\eqref{eq:struct-dimension-exponent}.
\end{theorem}

Suppose that the coefficient tuple is nonzero and irreducible and that $\lambda=1$.
Condition on all but one block column. The projection of that
column onto the orthogonal complement of the range of the other columns
is then a Gaussian matrix of fixed size. Conditionally on the other columns,
the mean of the exposed column enters only the conditional mean
\eqref{eq:struct-column-mean}. The conditional variances
\eqref{eq:struct-column-variance} depend on the coefficients and the
complement. A Gaussian matrix whose
covariance is bounded below has a small least singular value only with
small probability, whatever its mean (\cref{struct:lem:gaussian}).
The covariance can degenerate only if some Hilbert--Schmidt unit map
from the coefficient space $V$ into that orthogonal complement has all its $n$
blocks close to the orthogonal complement of the span of the
coefficients \eqref{eq:struct-coefficient-spaces}. The
restricted adjoint estimate (\cref{struct:thm:dual}), proved by a
covering argument, excludes this up to an exponentially small
probability. This gives \cref{prop:irreducible}. General tuples are
reduced to irreducible ones through a composition series
(\cref{subsec:triangular-reduction}), and a rescaling gives the
dependence on $\lambda$ (\cref{subsec:scalar-rescaling}).

For centered matrices the estimate is due to Han
\cite[Theorem 1.5]{Han}, and the coefficient geometry and
composition-series reduction follow \cite[Sections 3--4]{Han}.
\Cref{struct:thm:main} adds deterministic means, with a bound uniform
over bounded means, and a polynomial dependence on $\lambda$, which
allows $\lambda$ to tend to zero with $n$: in the almost-sure proof
$\lambda\asymp(\log n)^{-L}$. The exponent $2$ comes from the Gaussian
column-distance shortcut suggested in \cite[Remark 3.9]{Han}.
Independent random means are treated by conditioning and adding
the probability that their norm bound fails.

We prove the estimate for a generic tuple $R$ and mean bound $\omega$,
displaying the dependence on $L,M_1,\sigma_0,\vartheta$ below.
These statements can then be applied to $\lambda R$, $\omega=\omega_0/\lambda$
without repeating their proofs. Hilbert--Schmidt norms are unnormalized.

\subsection{The coefficient geometry}\label{subsec:coefficient}
Let $R_1,\ldots,R_k$ act irreducibly on a finite-dimensional complex
Hilbert space $V$, and suppose they are not all zero. On $\End(V)$
use the Hilbert--Schmidt inner product
$\ip AB=\Tr_V(A^*B)$, linear in the second variable, and set
\begin{equation}\label{eq:struct-coefficient-spaces}
 \Sigma=\Span\{R_1,\ldots,R_k\},\qquad H=\Sigma^\perp.
\end{equation}
We distinguish coefficient matrices from operators on $\End(V)$ by
boldface: left multiplication and its Hilbert--Schmidt adjoint are
\begin{equation}\label{eq:struct-left-multiplication}
 \mathbf R_\ell(A)=R_\ell A,\qquad
 \mathbf R_\ell^*(A)=R_\ell^*A,\qquad
 \|\mathbf R_\ell\|=\|R_\ell\|.
\end{equation}
Thus products such as $R_\ell^*A$ are ordinary matrix products, whereas
$\mathbf R_\ell^*\Gamma\mathbf R_\ell$ is a composition of operators
on $\End(V)$.
Write $P_E$ for orthogonal projection onto a subspace
$E\subset\End(V)$, and $\Gr(d,H)$ for the Grassmannian of
$d$-dimensional complex subspaces of $H$. All distances between
matrices use the Hilbert--Schmidt norm. An $\epsilon$-net covers a set
by balls of radius $\epsilon$; it is \emph{internal} if its centers
belong to that set. The case $H=0$ is treated directly in
\cref{subsec:exposed}.

Suppose henceforth that $H\ne0$, and put
\begin{equation}\label{eq:struct-pq}
 p=\dim\End(V)=(\dim V)^2,\qquad q=\dim H<p.
\end{equation}
For a nonzero subspace $E\subset H$, write $d=\dim E$ and define
\begin{equation}\label{eq:struct-coefficient-map}
 \Psi_E:E^k\longrightarrow\End(V),\qquad
 \Psi_E(x_1,\ldots,x_k)=\sum_{\ell=1}^kR_\ell^*x_\ell.
\end{equation}
For $\sigma>0$, retain the left singular vectors of $\Psi_E$ whose
singular values are at least $\sigma$:
\begin{equation}\label{eq:struct-active-image}
 \begin{gathered}
 Y_\sigma(E)=\Span\{\text{left singular vectors of $\Psi_E$
                            with singular value at least $\sigma$}\},\\
 \Psi_E^{<\sigma}=(I-P_{Y_\sigma(E)})\Psi_E.
 \end{gathered}
\end{equation}
We call $Y_\sigma(E)$ the active image; $\Psi_E^{<\sigma}$ contains
the remaining singular directions.

For subspaces $E,Y\subset\End(V)$ and $0<\vartheta<1$, let
$C_\vartheta(E,Y)$ be the span of the right singular vectors of
$P_{Y^\perp}|_E:E\to Y^\perp$ with singular value at most
$\vartheta$, including zero singular values. For $x\in E$, write
\begin{equation}\label{eq:struct-angle-decomposition}
 \begin{gathered}
 x=x_0+x_1,\qquad x_0\in C_\vartheta(E,Y),\quad
 x_1\in E\ominus C_\vartheta(E,Y),\\
 \|P_{Y^\perp}x_0\|_{\HS}\le\vartheta\|x_0\|_{\HS},\qquad
 \|P_{Y^\perp}x_1\|_{\HS}\ge\vartheta\|x_1\|_{\HS}.
 \end{gathered}
\end{equation}
Thus $C_\vartheta(E,Y)$ describes the directions in $E$ close to $Y$.
The following dimension inequality is the geometric input to the
covering argument.

\begin{lemma}[{\cite[Proposition 4.1]{Han}}]\label[lemma]{struct:lem:surplus}
There are constants $\sigma_0\in(0,1/2]$ and
$\vartheta\in(0,1/2)$ such that, for every nonzero $E\subset H$
and every $0<\sigma\le\sigma_0$,
\begin{equation}
 \dim Y_\sigma(E)\ge\dim C_\vartheta(E,Y_\sigma(E))+1.
 \label{eq:struct-2.1}
\end{equation}
In particular, the largest singular value of $\Psi_E$ is at least
$\sigma_0$.
\end{lemma}
\begin{proof}
We give the dimension-surplus argument of \cite[Proposition 4.1]{Han}
in the present notation. No nonzero subspace $E\subset H$ is invariant
under every $\mathbf R_j^*$. Otherwise Burnside's theorem makes $E$
a left ideal of $\End(V)$. For $B\in E$ and arbitrary $A\in\End(V)$,
orthogonality to $R_j$ gives
$0=\Tr(R_j^*AB)=\Tr(BR_j^*A)$, hence $BR_j^*=0$ for all $j$.
Thus $\ran B^*$ lies in the common kernel of the $R_j$, which is zero
by nonzero irreducibility. This contradicts $E\ne0$.
Compactness of the finitely many Grassmannians therefore gives
\begin{equation}\label{eq:struct-infimum}
 \inf_{0\ne E\subset H}\max_j
       \|P_{E^\perp}\mathbf R_j^*P_E\|>0.
\end{equation}
Put $Y=Y_\sigma(E)$ and $C=C_\vartheta(E,Y)$. Since
$\|P_{Y^\perp}|_C\|\le\vartheta<1$, projection into $Y$ is
injective, so $\dim C\le\dim Y$. If equality held, the elementary
principal-angle estimate would give
$\|P_{C^\perp}|_Y\|\le\vartheta/\sqrt{1-\vartheta^2}$.
Decomposing $\mathbf R_j^*|_E$ into its active and inactive parts
would then imply
\[
 \|P_{E^\perp}\mathbf R_j^*P_E\|
 \le \sigma+\|\Psi_E\|\vartheta/\sqrt{1-\vartheta^2}.
\]
The norms $\|\Psi_E\|$ are uniformly bounded. Taking $\sigma_0$
and $\vartheta$ sufficiently small contradicts
\eqref{eq:struct-infimum} and proves the strict dimension inequality. The proof remains
valid if either threshold is decreased. The largest singular value
must be at least $\sigma_0$, since otherwise $Y_{\sigma_0}(E)=0$.
\end{proof}
We fix these thresholds for the rest of the argument on $V$.

\Needspace{12\baselineskip}
\begin{lemma}\label[lemma]{struct:lem:nets}
Let $E,Y\subset\End(V)$, put $d=\dim E$, and let
$0<\vartheta<1$, $0<\epsilon<1$, and $h>0$. The set
\begin{equation}\label{eq:struct-tube}
 \mathcal T(E,Y,h)=
 \{x\in E:\|x\|_{\HS}\le2,\ \operatorname{dist}(x,Y)\le h\}
\end{equation}
has an internal $\epsilon$-net of size at most
\begin{equation}
 (C/\epsilon)^{2\dim C_\vartheta(E,Y)}
 \left(C\max\{h/\vartheta,\epsilon\}/\epsilon\right)^{2(d-\dim C_\vartheta(E,Y))}.
 \label{eq:struct-2.2}
\end{equation}
For $1\le d\le q$ and $0<\eta<1$, every subset of $\Gr(d,H)$
has an internal $\eta$-net in the metric $\|P_E-P_{E'}\|$, of size
at most $(C/\eta)^{2q^2}$. The constant $C$ depends only on $p,q$.
\end{lemma}
\begin{proof}
To prove \eqref{eq:struct-2.2}, decompose $x=x_0+x_1$ as in
\eqref{eq:struct-angle-decomposition}. The singular-value
decomposition makes $P_{Y^\perp}x_0$ and $P_{Y^\perp}x_1$ orthogonal.
Consequently, for $x\in\mathcal T(E,Y,h)$,
\[
 \|x_0\|_{\HS}\le2,\qquad
 \vartheta\|x_1\|_{\HS}\le\|P_{Y^\perp}x_1\|_{\HS}
 \le\|P_{Y^\perp}x\|_{\HS}\le h.
\]
Take the product of $\epsilon/4$-nets in the two indicated balls.
Since a complex space of dimension $b$ has real dimension $2b$,
the usual volume bound gives \eqref{eq:struct-2.2} for its cardinality.
This product is an $\epsilon/2$-net of the tube. Replacing each center
whose ball meets the tube by a point in that intersection gives an
internal $\epsilon$-net of no larger size.

For the Grassmannian estimate, $E\mapsto P_E$ embeds the given
subset in the operator-norm unit ball of the real $q^2$-dimensional
space of Hermitian operators on $H$. A maximal $\eta$-separated
subset is an internal $\eta$-net, and volume packing bounds its size
by $(1+2/\eta)^{q^2}\le(C/\eta)^{2q^2}$.
\end{proof}

\subsection{A restricted adjoint estimate and the covering argument}\label{subsec:dual}
Let $C_\ell,D_\ell\in M_n(\C)$ be as in \cref{struct:thm:main}.
Reordering tensor coordinates, write the pencil on $V^n$ as
\begin{equation}\label{eq:struct-block-pencil}
 (T_D)_{ij}=\delta_{ij}I_V+
       \sum_{\ell=1}^k(C_\ell(i,j)+D_\ell(i,j))R_\ell.
\end{equation}
For an index $a$, let $T_a:V\to V^n$ be block column $a$, and
let $T_{-a}:V^{n-1}\to V^n$ contain the remaining block columns.
Set
\begin{equation}\label{eq:struct-column-geometry}
 \begin{gathered}
 W_a=\ker T_{-a}^*=(\ran T_{-a})^\perp,\qquad
 F_a=P_{W_a}T_a:V\to W_a,\\
 \mathcal F_a=\sigma(C_\ell(i,j):1\le\ell\le k,\ 1\le i,j\le n,\ j\ne a).
 \end{gathered}
\end{equation}
A map $A:V\to V^n$, with components $A_i\in\End(V)$, takes values
in $W_a$ exactly when it satisfies the adjoint equations
\begin{equation}\label{eq:struct-normal-equations}
 T_{-a}^*A=0
 \quad\Longleftrightarrow\quad
 A_j+\sum_{i=1}^n\sum_{\ell=1}^k
 \overline{C_\ell(i,j)+D_\ell(i,j)}R_\ell^*A_i=0
 \quad(j\ne a).
\end{equation}
These are the adjoint kernel equations referred to in \cref{sec:interior}.
Conditionally on $\mathcal F_a$, $W_a$ is fixed and the entries
$C_\ell(i,a)$ remain independent circular Gaussians of variance $1/n$.
Thus the conditional mean of $F_a$ is
\begin{equation}\label{eq:struct-column-mean}
 \mu_a=\E(F_a\mid\mathcal F_a)
 =P_{W_a}\left(\delta_{ia}I_V+\sum_\ell D_\ell(i,a)R_\ell\right)_{i=1}^n.
\end{equation}
For an $\mathcal F_a$-measurable map $A:V\to W_a$, the identity
$P_{W_a}A=A$ gives
\[
 \langle A,F_a-\mu_a\rangle
 =\sum_{i,\ell}C_\ell(i,a)\langle A_i,R_\ell\rangle.
\]
Here maps from $V$ to $W_a$ carry their Hilbert--Schmidt inner product.
Independence yields
\begin{equation}\label{eq:struct-column-variance}
 \operatorname{Var}(\langle A,F_a\rangle\mid\mathcal F_a)
 =\E\bigl(|\langle A,F_a-\mu_a\rangle|^2\mid\mathcal F_a\bigr)
 =\frac1n\sum_{i,\ell}|\langle A_i,R_\ell\rangle|^2.
\end{equation}
By \eqref{eq:struct-coefficient-spaces}, this variance vanishes exactly
when all $A_i$ lie in $H$. We will show that a unit solution of
\eqref{eq:struct-normal-equations} cannot lie too close to $H^n$.
This is the purpose of the restricted adjoint estimate below.

It is convenient to prove that estimate with unconjugated entries.
Define
\begin{equation}\label{eq:struct-3.1}
 \begin{gathered}
 \mathcal L_{a,D}:H^n\longrightarrow\End(V)^{n-1},\qquad
 B=(B_1,\ldots,B_n),\\
 (\mathcal L_{a,D}B)_j=B_j+
 \sum_{i=1}^n\sum_{\ell=1}^k(C_\ell(i,j)+D_\ell(i,j))R_\ell^*B_i,
 \qquad j\ne a.
 \end{gathered}
\end{equation}
The term $B_j$ uses the inclusion $H\subset\End(V)$. Replacing
$C_\ell,D_\ell$ by their entrywise conjugates makes
\eqref{eq:struct-3.1} the restriction to $H^n$ of the map in
\eqref{eq:struct-normal-equations}. Entrywise conjugation preserves
the Ginibre law and the operator-norm bound on the means.
For tuples, use $\|B\|_{\HS}^2=\sum_i\|B_i\|_{\HS}^2$.

Introduce the norm events once for use throughout the proof:
\begin{equation}\label{eq:struct-norm-events}
 \Omega=\{\max_\ell\|C_\ell\|\le8\},\qquad
 \Omega_a=\{\max_\ell\|C_\ell|_{\{j\ne a\}}\|\le8\}.
\end{equation}
Here $C_\ell|_{\{j\ne a\}}:\C^{n-1}\to\C^n$ denotes the
matrix obtained by deleting column $a$. Restriction cannot increase the operator norm,
so $\Omega\subset\Omega_a$. By \cref{lem:norm},
\begin{equation}\label{eq:struct-norm-tails}
 \Pp(\Omega_a^c)\le\Pp(\Omega^c)
 \le\sum_{\ell=1}^k\Pp\{\|C_\ell\|>8\}\le Ce^{-cn}.
\end{equation}
The event $\Omega_a$ is $\mathcal F_a$-measurable. On $\Omega_a$,
\eqref{eq:struct-3.1} gives
\begin{equation}\label{eq:struct-3.3}
 \|\mathcal L_{a,D}\|\le L:=1+(8+\omega)\sum_\ell\|R_\ell\|.
\end{equation}
The same bound holds on all of $\End(V)^n$, and with conjugated
entries. Also set
\begin{equation}\label{eq:struct-mean-size}
 M_1=\sqrt{k}(8+\omega).
\end{equation}

For $B\in H^n$, let $E_B=\Span\{B_1,\ldots,B_n\}$ and define
\begin{equation}\label{eq:struct-tuple-covariance}
 \Gamma_Bx=\sum_i B_i\ip{B_i}{x}\qquad(x\in\End(V)).
\end{equation}
This positive operator measures how $B$ spans $E_B$:
\begin{equation}\label{eq:struct-covariance-properties}
 \langle x,\Gamma_Bx\rangle=\sum_i|\langle B_i,x\rangle|^2,
 \qquad\ran\Gamma_B=E_B,\qquad
 \Tr_H(\Gamma_B|_H)=\|B\|_{\HS}^2.
\end{equation}
It vanishes on $E_B^\perp$.

\begin{theorem}\label[theorem]{struct:thm:dual}
There are an integer $d_*=d_*(p,q)>0$, constants $c,C>0$, and
an integer $n_0$ such that
\begin{equation}
 \Pp\{\sg(\mathcal L_{a,D})\le n^{-d_*}\}\le Ce^{-cn}
 \label{eq:struct-3.2}
\end{equation}
for every $n\ge n_0$, every deleted index $a$, and every deterministic
tuple with $\max_\ell\|D_\ell\|\le \omega$.
More precisely, with $\ell$ from \eqref{eq:struct-gap-exponents}, put
\begin{equation}\label{eq:struct-dual-exponents}
 u_1=1,\qquad u_{j+1}=2\ell(u_j)+4,\qquad
 d_*=\ell(u_q)+1.
\end{equation}
One can take
$n_0=\lceil\max\{\max_j n_0(u_j),q,3,\sqrt{2L\sqrt q}\}\rceil$.
The exponential constants depend only on $k,p,q$; all coefficient and
mean dependence is in this threshold.
\end{theorem}
For $D=0$, this is the dual estimate of \cite[Theorem 3.3]{Han}.
We first control tuples satisfying $\Gamma_B\ge\mu P_{E_B}$ at a
fixed scale (\cref{struct:lem:regular}), then choose a scale for each
subspace (\cref{struct:lem:gap}). Finally, projecting onto the
eigendirections of $\Gamma_B$ with sufficient mass reduces every
unit tuple to one of these cases.

For a Borel class $\mathcal G$ of nonzero subspaces of $H$, write
\begin{equation}\label{eq:struct-regular-event}
 \mathcal A(\mathcal G,\mu,\delta)=
 \left\{\exists B\in H^n:
 \begin{array}{l}
 E_B\in\mathcal G,\quad\|B\|_{\HS}\le2,\\
 \Gamma_B\ge\mu P_{E_B},\quad\|\mathcal L_{a,D}B\|_{\HS}\le\delta
 \end{array}\right\}.
\end{equation}
The event in \eqref{eq:struct-regular-event} is measurable in the
completed Gaussian probability space: it is the projection of a Borel
set in the finite-dimensional entry and tuple variables, hence analytic
and universally measurable.
The next lemma combines the covering argument of
\cite[Proposition 4.5]{Han} with the treatment of small singular
directions in \cite[proof of Proposition 4.8, (4.28)--(4.31)]{Han}.
We impose a single bound on $\Psi_E^{<\sigma}$. The structured means
satisfy the same bound on these directions as the Gaussian terms;
on the active image, they change only the mean of the Gaussian law.

\begin{lemma}\label[lemma]{struct:lem:regular}
Fix $\sigma_0,\vartheta$ as in \cref{struct:lem:surplus}.
There is $C_{\mathrm{reg}}\ge1$, depending only on
$p,q,L,\vartheta$, with the following property. Let $n\ge2$,
$0<\mu\le1$, $0<\sigma\le\sigma_0$, and $0<\delta\le\mu$.
Let $1\le d\le q$ and let $\mathcal G$ be a Borel class of
$d$-dimensional subspaces of $H$ satisfying
\begin{equation}
 \|\Psi_E^{<\sigma}\|\le\delta/(100M_1),\qquad E\in\mathcal G.
 \label{eq:struct-3.5}
\end{equation}
Then, with $\alpha_{\mathrm{net}}=2q+2q^2$,
\begin{equation}
 \Pp\bigl(\Omega_a\cap\mathcal A(\mathcal G,\mu,\delta)\bigr)
 \le\delta^{2(n-1)-\alpha_{\mathrm{net}}}
       \left(C_{\mathrm{reg}}n^{2p}\mu^{-p}\sigma^{-2p}\right)^n.
 \label{eq:struct-3.6}
\end{equation}
The estimate is uniform in $a$ and in the deterministic means with
$\max_\ell\|D_\ell\|\le \omega$.
\end{lemma}
\begin{proof}
If $\mathcal G=\varnothing$, the claim is immediate. Otherwise, for a
tuple in $\mathcal A(\mathcal G,\mu,\delta)$,
\begin{equation}\label{eq:struct-frame-bound}
 \langle x,\Gamma_Bx\rangle
 =\sum_i|\langle B_i,x\rangle|^2
 \ge\mu\|x\|_{\HS}^2\qquad(x\in E_B).
\end{equation}
We work on $\Omega_a$ when estimating approximation errors.

\textit{Step 1: approximate the subspace.}
Recall the operator bound $\|\mathcal L_{a,D}\|\le L$ from
\eqref{eq:struct-3.3}. By \cref{struct:lem:nets}, choose an internal $\eta$-net
$\mathcal N_{\mathcal G}$ of $\mathcal G$, where $\eta=\delta/(100L)$.
Choose $E'\in\mathcal N_{\mathcal G}$ with $\|P_{E_B}-P_{E'}\|\le\eta$
and put $B'_i=P_{E'}B_i$. Then
\[
 \|B-B'\|_{\HS}\le2\eta,\quad \|B'\|_{\HS}\le2,\quad
 \|\mathcal L_{a,D}B'\|_{\HS}\le\delta+2L\eta<2\delta.
\]
Since $\Gamma_B$ vanishes on $E_B^\perp$, \eqref{eq:struct-frame-bound}
gives, as an operator on $E'$,
\begin{equation}
 \Gamma_{B'}=P_{E'}\Gamma_BP_{E'}
 \ge\mu P_{E'}P_{E_B}P_{E'}
 \ge\mu(1-\eta^2)P_{E'}\ge(\mu/2)P_{E'}.
 \label{eq:struct-3.7}
\end{equation}
Indeed, $\|P_{E_B^\perp}x\|_{\HS}\le\eta\|x\|_{\HS}$ for
$x\in E'$, and $\eta\le1/100$ since $\delta\le1$ and $L\ge1$.
The positive lower bound in \eqref{eq:struct-3.7} also shows that
$E_{B'}=E'$.

\textit{Step 2: cover tuples with small $\|\mathcal L_{a,D}B'\|_{\HS}$.}
Recall the active image and the operator $\Psi_E^{<\sigma}$ from
\eqref{eq:struct-active-image}. On $\Omega_a$, for $E\subset H$
and $B\in E^n$, we have
\begin{equation}\label{eq:struct-3.4}
 \begin{aligned}
 &\left\|\left(\Psi_E^{<\sigma}
       \left(\sum_i(C_\ell(i,j)+D_\ell(i,j))B_i\right)_{\ell=1}^k
                    \right)_{j\ne a}\right\|_{\HS}\\
 &\quad\le\|\Psi_E^{<\sigma}\|
   \left(\sum_{\ell=1}^k\sum_{j\ne a}
     \left\|\sum_i(C_\ell(i,j)+D_\ell(i,j))B_i\right\|_{\HS}^2\right)^{1/2}\\
 &\quad\le\|\Psi_E^{<\sigma}\|\sqrt{k}(8+\omega)\|B\|_{\HS}
   =M_1\|\Psi_E^{<\sigma}\|\,\|B\|_{\HS}.
 \end{aligned}
\end{equation}
The second inequality uses
$\|((C_\ell+D_\ell)|_{\{j\ne a\}})^{\mathsf T}\|\le8+\omega$
on $\Omega_a$, also after tensoring by $E$.
Projecting \eqref{eq:struct-3.1} onto $Y_\sigma(E')^\perp$
eliminates the active parts of both the random term and the mean.
For $j\ne a$, \eqref{eq:struct-3.4} therefore gives
\[\begin{aligned}
 \operatorname{dist}(B'_j,Y_\sigma(E'))
 &=\|P_{Y_\sigma(E')^\perp}B'_j\|_{\HS}\\
 &=\left\|P_{Y_\sigma(E')^\perp}(\mathcal L_{a,D}B')_j
   -\Psi_{E'}^{<\sigma}
    \left(\sum_i(C_\ell(i,j)+D_\ell(i,j))B'_i\right)_{\ell=1}^k
   \right\|_{\HS}\\
 &\le\|(\mathcal L_{a,D}B')_j\|_{\HS}
       +M_1\|\Psi_{E'}^{<\sigma}\|\,\|B'\|_{\HS}
 \le2\delta+\delta/50<3\delta,
\end{aligned}\]
by \eqref{eq:struct-3.4} and \eqref{eq:struct-3.5}.
Use \cref{struct:lem:nets} to choose an internal $\epsilon$-net
$\mathcal N_{\mathrm{tube}}(E')$ of
$\mathcal T(E',Y_\sigma(E'),3\delta)$, where
$\epsilon=\delta/(100L\sqrt n)$.
Choose an internal $\epsilon$-net $\mathcal N_a(E')$ of
$\{x\in E':\|x\|_{\HS}\le2\}$ for the deleted coordinate.
Approximating each $B'_j$ by a center in the corresponding net gives
\begin{equation}
 \|W-B'\|_{\HS}\le\sqrt n\epsilon,\quad
 \|W\|_{\HS}\le3,\quad \|\mathcal L_{a,D}W\|_{\HS}\le3\delta.
 \label{eq:struct-3.8}
\end{equation}
Moreover,
\[
 \|\Gamma_W-\Gamma_{B'}\|
 \le(\|W\|_{\HS}+\|B'\|_{\HS})\|W-B'\|_{\HS}
 \le5\sqrt n\epsilon=\frac{\delta}{20L}\le\frac\mu4.
\]
Together with \eqref{eq:struct-3.7}, this gives
$\Gamma_W\ge(\mu/4)P_{E'}$.
Let $\mathcal N(E')$ consist of the tuples in the product net
satisfying this covariance bound. Every $B'$ under consideration
has an approximant in $\mathcal N(E')$.

\textit{Step 3: compare the Gaussian probability with the net size.}
For a fixed $E'\in\mathcal N_{\mathcal G}$, put
\begin{equation}\label{eq:struct-b-cstar}
 b_*=\dim Y_\sigma(E'),\qquad
 c_*=\dim C_\vartheta(E',Y_\sigma(E')).
\end{equation}
For a fixed $W\in\mathcal N(E')$ and $j\ne a$, write
\begin{equation}\label{eq:struct-mean-fluctuation}
 (\mathcal L_{a,D}W)_j=m_j+\xi_j,\qquad
 m_j=W_j+\sum_{i,\ell}D_\ell(i,j)R_\ell^*W_i,\qquad
 \xi_j=\sum_{i,\ell}C_\ell(i,j)R_\ell^*W_i.
\end{equation}
Independence and $\E|C_\ell(i,j)|^2=1/n$ give
\begin{equation}\label{eq:struct-3.9}
 \begin{aligned}
 \operatorname{Cov}(\xi_j)
 &=\frac1n\sum_\ell \mathbf R_\ell^*\Gamma_W\mathbf R_\ell
 \ge\frac\mu{4n}\sum_\ell \mathbf R_\ell^* P_{E'}\mathbf R_\ell\\
 &=\frac\mu{4n}\Psi_{E'}\Psi_{E'}^*
 \ge\frac{\mu\sigma^2}{4n}P_{Y_\sigma(E')}.
 \end{aligned}
\end{equation}
Here the boldface operators are defined in
\eqref{eq:struct-left-multiplication}; the last inequality uses
\eqref{eq:struct-active-image}.
Since $\xi_j$ is a complex-linear combination of independent circular
Gaussians, the projection of $m_j+\xi_j$ onto $Y_\sigma(E')$ has a
Gaussian density bounded by $(4n/(\pi\mu\sigma^2))^{b_*}$, uniformly
in its mean. Multiplying by the volume of a radius-$3\delta$ ball
in $\C^{b_*}$ gives
\begin{equation}\label{eq:struct-3.10}
 \begin{aligned}
 \Pp\{\|(\mathcal L_{a,D}W)_j\|_{\HS}\le3\delta\}
 &\le\Pp\{\|P_{Y_\sigma(E')}(m_j+\xi_j)\|_{\HS}\le3\delta\}\\
 &\le\left(\frac{4n}{\pi\mu\sigma^2}\right)^{b_*}
       \times\frac{\pi^{b_*}(3\delta)^{2b_*}}{b_*!}
 \le\left(\frac{Cn\delta^2}{\mu\sigma^2}\right)^{b_*}.
 \end{aligned}
\end{equation}
The vectors $\xi_j$, $j\ne a$, are independent because they use
different Gaussian columns. Hence, without conditioning on $\Omega_a$,
\[
 \Pp\{\|\mathcal L_{a,D}W\|_{\HS}\le3\delta\}
 \le\prod_{j\ne a}\Pp\{\|(\mathcal L_{a,D}W)_j\|_{\HS}\le3\delta\}
 \le\left(\frac{Cn\delta^2}{\mu\sigma^2}\right)^{b_*(n-1)}.
\]
By \cref{struct:lem:nets},
\[
 \#\mathcal N_{\mathrm{tube}}(E')
 \le(C\sqrt n/\delta)^{2c_*}(C\sqrt n)^{2(d-c_*)},
 \qquad \#\mathcal N_a(E')\le(C\sqrt n/\delta)^{2d}.
\]
In particular,
\[
 \#\mathcal N(E')
 \le\#\mathcal N_a(E')\,\#\mathcal N_{\mathrm{tube}}(E')^{n-1}
 \le C^n n^{dn}\delta^{-2c_*(n-1)-2d}.
\]
All these nets are deterministic, and Steps 1--2 give the inclusion
\[
 \Omega_a\cap\mathcal A(\mathcal G,\mu,\delta)
 \subset\bigcup_{E'\in\mathcal N_{\mathcal G}}\ \bigcup_{W\in\mathcal N(E')}
       \{\|\mathcal L_{a,D}W\|_{\HS}\le3\delta\}.
\]
Since $\#\mathcal N_{\mathcal G}\le(CL/\delta)^{2q^2}$, the union bound yields
\begin{align*}
 \Pp\bigl(\Omega_a\cap\mathcal A(\mathcal G,\mu,\delta)\bigr)
 &\le\sum_{E'\in\mathcal N_{\mathcal G}}\#\mathcal N(E')
             \left(\frac{Cn\delta^2}{\mu\sigma^2}\right)^{b_*(n-1)}\\*
 &\le C^n\max_{E'\in\mathcal N_{\mathcal G}}
 \frac{n^{dn+b_*(n-1)}\delta^{2(b_*-c_*)(n-1)-2d-2q^2}}
      {(\mu\sigma^2)^{b_*(n-1)}}\\*
 &\le\delta^{2(n-1)-\alpha_{\mathrm{net}}}(Cn^{2p}\mu^{-p}\sigma^{-2p})^n.
\end{align*}
For the last inequality, $E'\subset H\subset\End(V)$ and
$Y_\sigma(E')\subset\End(V)$ imply $d\le q<p$ and $b_*\le p$;
\cref{struct:lem:surplus} gives $b_*-c_*\ge1$. Thus, since
$\delta,\mu,\sigma\le1$,
\[
 \delta^{2(b_*-c_*)(n-1)-2d-2q^2}\le\delta^{2(n-1)-\alpha_{\mathrm{net}}},\qquad
 n^{dn+b_*(n-1)}\le n^{2pn},\qquad
 (\mu\sigma^2)^{-b_*(n-1)}\le(\mu^{-p}\sigma^{-2p})^n.
\]
The nets use only $p,q,L,\vartheta$, and the active Gaussian density
bound keeps $\sigma$ explicit. Thus their combined constant is
$C_{\mathrm{reg}}=C_{\mathrm{reg}}(p,q,L,\vartheta)$; neither $M_1$
nor a separate coefficient bound enters it. This proves \eqref{eq:struct-3.6}.
\end{proof}

The nonzero singular values of $\Psi_E$ need not be bounded away
from zero as $E$ varies. The next lemma chooses a scale at which
\cref{struct:lem:regular} applies. Its conclusion corresponds to
\cite[Proposition 4.9]{Han}; in place of the gap classification and
dyadic discretization in \cite[Facts 4.6--4.7 and Proposition 4.8]{Han},
we use a sequence of $p+1$ thresholds. The thresholds depend on $n$,
but their number depends only on $p=\dim\End(V)$.

\begin{lemma}\label[lemma]{struct:lem:gap}
Fix $\sigma_0\in(0,1/2]$ and $\vartheta\in(0,1/2)$ as in
\cref{struct:lem:surplus}.
For each real $u\ge1$, there is an integer $\ell(u)$, depending only
on $p,u$, and a threshold $n_0(u)$, depending only on $u$, the coefficient
tuple $(R_\ell)$ and $\omega$, such that, for $n\ge n_0(u)$,
\begin{equation}\label{eq:struct-gap-probability}
 \Pp\left(\Omega_a\cap\bigcup_{d=1}^q
 \mathcal A(\Gr(d,H),n^{-u},n^{-\ell(u)})\right)\le e^{-n\log n}.
\end{equation}
An explicit choice, nondecreasing in $u$, is
\begin{equation}\label{eq:struct-gap-exponents}
 \begin{gathered}
 \kappa(u)=2p(u+2)+10,\quad
 \beta(p,u)=(\kappa(u)+2)\sum_{j=0}^{p-1}(2p)^j-2,\quad
 \ell(u)=\lceil\beta(p,u)\rceil+1,\\
 n_0(u)=\left\lceil\max\{3,\,2(\alpha_{\mathrm{net}}+2),\,
 10\sqrt{M_1},\,\sigma_0^{-(2p)^p},\,C_{\mathrm{reg}},\,pq\}\right\rceil.
 \end{gathered}
\end{equation}
Thus the coefficient dependence is explicit through
$M_1,\sigma_0,C_{\mathrm{reg}}$.
\end{lemma}
\begin{proof}
\textit{Step 1: partition the subspaces according to a singular-value gap.}
Put $\kappa(u)=2p(u+2)+10$ and define
\begin{equation}\label{eq:struct-3.11}
 t_0=\sigma_0,\qquad t_{j+1}=t_j^{2p}n^{-\kappa(u)-2},\qquad
 \delta_j=t_j^{2p}n^{-\kappa(u)}\quad(0\le j<p).
\end{equation}
Since $0<t_{j+1}<t_j\le\sigma_0$, \cref{struct:lem:surplus} gives
\[
 1\le\dim Y_{t_0}(E)\le\dim Y_{t_1}(E)\le\cdots
   \le\dim Y_{t_p}(E)\le p\qquad(0\ne E\subset H).
\]
At least one of the $p$ consecutive pairs must therefore be equal.
For $1\le d\le q$ and $0\le j<p$, define
\begin{equation}\label{eq:struct-gap-classes}
 \mathcal G_{d,j}=\{E\in\Gr(d,H):
       \dim Y_{t_j}(E)=\dim Y_{t_{j+1}}(E)\}.
\end{equation}
These classes satisfy
\[
 \Gr(d,H)=\bigcup_{j=0}^{p-1}\mathcal G_{d,j},\qquad
 \|\Psi_E^{<t_j}\|<t_{j+1}=\delta_j/n^2
 \le\delta_j/(100M_1)\quad(E\in\mathcal G_{d,j})
\]
for $n^2\ge100M_1$: equality of dimensions excludes singular values
in $[t_{j+1},t_j)$. Each class is Borel, because the singular values
of $(x_\ell)\mapsto\sum_\ell R_\ell^*P_Ex_\ell$ depend continuously
on $P_E$. Thus these classes satisfy the hypotheses of
\cref{struct:lem:regular}.

\textit{Step 2: bound the probabilities and choose a common threshold.}
Apply \eqref{eq:struct-3.6} with $\mu=n^{-u}$, $\sigma=t_j$,
and $\delta=\delta_j$. Since $t_j\le1$ and $\kappa(u)>u$,
$\delta_j\le\mu$ for $n\ge2$. For $n\ge2(\alpha_{\mathrm{net}}+2)$,
\begin{align*}
 \Pp\bigl(\Omega_a\cap\mathcal A(\mathcal G_{d,j},n^{-u},\delta_j)\bigr)
 &\le\delta_j^{3n/2}(Cn^{p(u+2)}t_j^{-2p})^n\\
 &=\bigl(Ct_j^p n^{-2p(u+2)-15}\bigr)^n\le n^{-10n}
\end{align*}
for $n\ge C_{\mathrm{reg}}$, after using the constant in
\eqref{eq:struct-3.6}. All size conditions are included in
\eqref{eq:struct-gap-exponents}.
The thresholds $\delta_j$ decrease with $j$, and the recursion
\eqref{eq:struct-3.11} gives
\[
 \min_{0\le j<p}\delta_j=\delta_{p-1}
 =n^2t_p=\sigma_0^{(2p)^p}n^{-\beta(p,u)},\qquad
 \beta(p,u)=(\kappa(u)+2)\sum_{j=0}^{p-1}(2p)^j-2.
\]
Set $\ell(u)=\lceil\beta(p,u)\rceil+1$. For
$n\ge\sigma_0^{-(2p)^p}$ we have $n^{-\ell(u)}\le\delta_j$ for
all $0\le j<p$. Hence
\begin{align*}
 \Pp\left(\Omega_a\cap\bigcup_{d=1}^q
       \mathcal A(\Gr(d,H),n^{-u},n^{-\ell(u)})\right)
 &\le\sum_{d=1}^q\sum_{j=0}^{p-1}
       \Pp\bigl(\Omega_a\cap\mathcal A(\mathcal G_{d,j},n^{-u},\delta_j)\bigr)\\
 &\le pq\,n^{-10n}\le e^{-n\log n}
\end{align*}
for $n\ge n_0(u)$. The integer $\ell(u)$ is nondecreasing in $u$ and depends
only on $p,u$; $\sigma_0$ and the other coefficient-dependent constants
enter only the lower threshold $n_0(u)$.
\end{proof}

\begin{proof}[Proof of \cref{struct:thm:dual}]
\textit{Step 1: choose thresholds for every possible dimension.}
Use the exponents in \eqref{eq:struct-dual-exponents}; then
$u_1<\cdots<u_q$.
By \cref{struct:lem:gap}, there is an event $\mathcal E\subset\Omega_a$
with $\Pp(\Omega_a\setminus\mathcal E)\le q e^{-n\log n}$ such that,
on $\mathcal E$, for every $1\le j\le q$ and $W\in H^n$,
\begin{equation}\label{eq:struct-regular-implication}
 \left.
 \begin{gathered}
 E_W\ne0,\quad \|W\|_{\HS}\le2,\\
 \Gamma_W\ge n^{-u_j}P_{E_W}
 \end{gathered}\right\}
 \quad\Longrightarrow\quad \|\mathcal L_{a,D}W\|_{\HS}>n^{-\ell(u_j)}.
\end{equation}
Work on $\mathcal E$ for the rest of the deterministic argument.

\textit{Step 2: retain directions with sufficient mass.}
For $B\in H^n$ with $\|B\|_{\HS}=1$, choose an orthonormal eigenbasis
$e_1,\ldots,e_q$ of $\Gamma_B|_H$, with eigenvalues
$\lambda_1\ge\cdots\ge\lambda_q\ge0$. Since
\[
 \sum_{j=1}^q\lambda_j=\Tr_H\Gamma_B=1,\qquad
 \lambda_1\ge1/q\ge n^{-u_1}\quad(n\ge q),
\]
the integer
\[
 d=\max\{j\in\{1,\ldots,q\}:\lambda_j\ge n^{-u_j}\}
\]
is well defined and at least one. Put $E'=\Span\{e_1,\ldots,e_d\}$
and $B'_i=P_{E'}B_i$. Then $\|B'\|_{\HS}\le1$ and
\[
 \Gamma_{B'}=P_{E'}\Gamma_BP_{E'}
 =\sum_{j=1}^d\lambda_j e_je_j^*
 \ge n^{-u_d}P_{E'},\qquad E_{B'}=E'.
\]
For $d<q$, the definition of $d$ and monotonicity of $u_j$ give
\[
 \|B-B'\|_{\HS}^2=\sum_{j>d}\lambda_j
 <\sum_{j>d}n^{-u_j}\le qn^{-u_{d+1}}.
\]

\textit{Step 3: transfer the bound to the original tuple.}
For $d<q$, \eqref{eq:struct-3.3} and
$u_{d+1}=2\ell(u_d)+4$ imply
\[
 \|\mathcal L_{a,D}(B-B')\|_{\HS}
 \le L\sqrt q\,n^{-u_{d+1}/2}
 =L\sqrt q\,n^{-\ell(u_d)-2}\le\tfrac12n^{-\ell(u_d)}
\]
for large $n$. Apply \eqref{eq:struct-regular-implication} to $W=B'$:
\[
 \|\mathcal L_{a,D}B\|_{\HS}
 \ge\|\mathcal L_{a,D}B'\|_{\HS}
      -\|\mathcal L_{a,D}(B-B')\|_{\HS}
 >\tfrac12n^{-\ell(u_d)}\ge\tfrac12n^{-\ell(u_q)}.
\]
If $d=q$, $B'=B$ and the same lower bound follows directly.
Taking the infimum over unit tuples, for $n>2$ and
$d_*=\ell(u_q)+1$ we obtain
\[
 \sg(\mathcal L_{a,D})\ge\tfrac12n^{-\ell(u_q)}>n^{-d_*}
 \qquad\text{on }\mathcal E.
\]
Adding the exceptional probability in \eqref{eq:struct-norm-tails},
we obtain
\begin{align*}
 \Pp\{\sg(\mathcal L_{a,D})\le n^{-d_*}\}
 &\le\Pp(\Omega_a^c)+\Pp(\Omega_a\setminus\mathcal E)\\*
 &\le Ce^{-cn}+q e^{-n\log n}\le C'e^{-c'n}.
\end{align*}
This proves \eqref{eq:struct-3.2}.
\end{proof}

\subsection{Exposing one block column}\label{subsec:exposed}
We now pass from the restricted adjoint estimate to the least singular
value of $T_D$, following \cite[Sections 3.4--3.6]{Han}. For the
earlier use of block-column projections in polynomial linearizations,
see \cite{CGH}. We first isolate the Gaussian estimate used in this
step. It implements the shortcut suggested in \cite[Remark 3.9]{Han}
by the Gaussian smoothed-analysis column-distance method
\cite[Section 3]{SST}; compare
\cite[Lemma 3.5 and Section 3.4]{RV}.
The elementary proof below uses a complex Gaussian small ball, which
gives a quadratic rather than a linear power of the radius.

A complex Gaussian vector $Z$ is \emph{proper} if
$\E[(Z-\E Z)(Z-\E Z)^{\mathsf T}]=0$. Its covariance is
$\E[(Z-\E Z)(Z-\E Z)^*]$. We apply this terminology to a Gaussian
matrix through its vectorization, which lists its entries in a fixed order.

\begin{lemma}\label[lemma]{struct:lem:gaussian}
Let $X\in M_r(\C)$ be a proper Gaussian matrix with arbitrary mean.
If $\operatorname{Cov}(\operatorname{vec}X)\ge\lambda I_{r^2}$
for some $\lambda>0$, then
\begin{equation}\label{eq:struct-4.1}
 \Pp\{\sg(X)\le u\}\le r^2u^2/\lambda\qquad(u>0).
\end{equation}
\end{lemma}
\begin{proof}
Split the covariance to write, in distribution,
$X=\E X+\sqrt\lambda G+G_0$, where $G$ has independent standard
circular entries and $G_0$ is an independent centered proper Gaussian
matrix with covariance
$\operatorname{Cov}(\operatorname{vec}X)-\lambda I_{r^2}$.
Conditionally on $G_0$, the columns $x_1,\ldots,x_r$ are independent,
each with covariance $\lambda I_r$. Put
$E_j=\Span\{x_i:i\ne j\}$.
If $\sg(X)\le u$, a unit vector $v$ with $\|Xv\|\le u$ has a
coordinate $|v_j|\ge r^{-1/2}$. Projection onto $E_j^\perp$ gives
$|v_j|\operatorname{dist}(x_j,E_j)\le u$, whence
\[
 \{\sg(X)\le u\}\subseteq\bigcup_{j=1}^r
       \{\operatorname{dist}(x_j,E_j)\le\sqrt r\,u\}.
\]
For each $j$, condition further on the other columns and choose a
measurable unit vector $q_j\in E_j^\perp$. The scalar $q_j^*x_j$
is circular Gaussian with variance $\lambda$ and arbitrary mean.
Its density is bounded by $1/(\pi\lambda)$, so
\begin{align*}
 \Pp\{\operatorname{dist}(x_j,E_j)\le\sqrt r\,u
                          \mid G_0,(x_i)_{i\ne j}\}
 &\le\Pp\{|q_j^*x_j|\le\sqrt r\,u\mid G_0,(x_i)_{i\ne j}\}\\*
 &\le ru^2/\lambda.
\end{align*}
Taking expectations and summing over $j$ proves
\eqref{eq:struct-4.1}.
\end{proof}

\Needspace{12\baselineskip}
\begin{proposition}\label[proposition]{prop:irreducible}
Let $R_1,\ldots,R_k\in\End(V)$ act irreducibly and be not all zero,
and fix $\omega<\infty$. Define $T_D$ as in \cref{struct:thm:main}.
There are an integer $d_*\ge0$, constants $C_R,C,c>0$, and an
integer $n_0$ such that, for every $n\ge n_0$, every deterministic
tuple with $\max_\ell\|D_\ell\|\le \omega$, and every $t>0$,
\begin{equation}
 \Pp\{\sg(T_D)\le t\}\le C_R n^{2d_*+3}t^2+Ce^{-cn}.
 \label{eq:struct-5.5}
\end{equation}
The exponent $d_*$ depends only on $p=(\dim V)^2$ and
$q=\dim\Span\{R_1,\ldots,R_k\}^\perp$; the constants $C_R,C,c,n_0$
depend only on $(R_\ell)$ and $\omega$. More precisely, if $r=\dim V$
and $s_R$ is the least singular value of the coefficient frame
$B\mapsto(\langle B,R_\ell\rangle)_\ell$ on $\Sigma$, one may take
\begin{equation}\label{eq:struct-column-constant}
 C_R=r^2(L^2+1)/s_R^2.
\end{equation}
The threshold is that of \cref{struct:thm:dual}, and the exponential
constants depend only on $k,p,q$. For $H=0$ no dual threshold is needed.
\end{proposition}
When $H\ne0$, we use $d_*$ from \cref{struct:thm:dual}.
When $H=0$, the proof gives $d_*=0$ directly.

\begin{proof}
\textit{Step 1: the deleted-column space.}
Let $r=\dim V$, and recall $T_a,T_{-a},W_a,F_a$ and $\mathcal F_a$
from \eqref{eq:struct-column-geometry}.
The determinant of $T_D$ is a polynomial in the Ginibre entries
which equals one at $C_\ell=-D_\ell$. Its zero set has Lebesgue
measure zero, so $T_D$ is almost surely invertible. It follows that
\[
 \operatorname{rank}T_{-a}=(n-1)r,\qquad\dim W_a=r.
\]
Thus $F_a:V\to W_a$ is a square map of fixed dimension $r$.

\textit{Step 2: the conditional covariance.}
Let $s_R>0$ be the least singular value of the conjugate-linear map
\begin{equation}\label{eq:struct-frame}
 \Sigma\longrightarrow\C^k,\qquad
 B\longmapsto(\langle B,R_\ell\rangle)_\ell.
\end{equation}
This map is injective because the $R_\ell$ span $\Sigma$, so $s_R$
is positive. First assume $H\ne0$, and write
\begin{equation}\label{eq:struct-5.1}
 \begin{gathered}
 \widetilde{\mathcal L}_{a,D}:\End(V)^n\longrightarrow\End(V)^{n-1},
 \\
 (\widetilde{\mathcal L}_{a,D}A)_j
 =A_j+\sum_{i,\ell}\overline{C_\ell(i,j)+D_\ell(i,j)}R_\ell^*A_i.
 \end{gathered}
\end{equation}
This is the map $A\mapsto T_{-a}^*A$ in
\eqref{eq:struct-normal-equations}. Its restriction to $H^n$ is
\eqref{eq:struct-3.1} with entrywise conjugated Gaussian entries
and means. Consequently, \cref{struct:thm:dual} applies to this
restriction, and \eqref{eq:struct-3.3} gives
$\|\widetilde{\mathcal L}_{a,D}\|\le L$ on $\Omega_a$.

Condition on $\mathcal F_a$ and work on the event
\begin{equation}\label{eq:struct-column-good-event}
 \mathcal G_a=\Omega_a\cap
 \{\sg(\widetilde{\mathcal L}_{a,D}|_{H^n})>n^{-d_*}\}
 \cap\{\operatorname{rank}T_{-a}=(n-1)r\}.
\end{equation}
It is $\mathcal F_a$-measurable, and
$\Pp(\mathcal G_a^c)\le Ce^{-cn}$ by \cref{struct:thm:dual}
and \eqref{eq:struct-norm-tails}.
For a Hilbert--Schmidt unit map $A:V\to W_a$, decompose each block as
$A_i=B_i+E_i$, where $B_i\in H$ and $E_i\in\Sigma$, and put
$\zeta=\|E\|_{\HS}$. Orthogonality gives
$\|B\|_{\HS}^2=1-\zeta^2$. Since
$\widetilde{\mathcal L}_{a,D}A=0$, the lower bound in
\eqref{eq:struct-column-good-event} and the upper bound
\eqref{eq:struct-3.3} imply
\begin{equation}\label{eq:struct-5.2}
 \begin{gathered}
 n^{-d_*}\sqrt{1-\zeta^2}
 \le\|\widetilde{\mathcal L}_{a,D}B\|_{\HS}
 =\|\widetilde{\mathcal L}_{a,D}E\|_{\HS}\le L\zeta,\\
 \zeta^2\ge\frac{n^{-2d_*}}{L^2+1}.
 \end{gathered}
\end{equation}
Using \eqref{eq:struct-column-variance}, $B_i\perp\Sigma$, and
the definition of $s_R$, we obtain
\begin{equation}\label{eq:struct-5.3}
 \begin{aligned}
 \E\bigl(|\langle A,F_a-\mu_a\rangle|^2\mid\mathcal F_a\bigr)
 &=\frac1n\sum_{i,\ell}|\langle E_i,R_\ell\rangle|^2\\*
 &\ge\frac{s_R^2}{n}\|E\|_{\HS}^2
 \ge\frac{s_R^2}{L^2+1}n^{-(2d_*+1)}.
 \end{aligned}
\end{equation}

Choose a fixed orthonormal basis of $V$ and an
$\mathcal F_a$-measurable orthonormal basis of $W_a$. In these
bases put
\begin{equation}\label{eq:struct-Za}
 Z_a=\operatorname{vec}(F_a-\mu_a),\qquad
 Q_a=\E(Z_aZ_a^*\mid\mathcal F_a),\qquad
 \lambda_n=\frac{s_R^2}{L^2+1}n^{-(2d_*+1)}.
\end{equation}
The conditional Gaussian law of $Z_a$ is proper, since $Z_a$ is a
complex-linear function of independent centered circular Gaussians.
For $u\in\C^{r^2}$, let $A:V\to W_a$ be the map with
$\operatorname{vec}A=u$. By homogeneity, \eqref{eq:struct-5.3} gives
\[
 u^*Q_au
 =\E(|\langle A,F_a-\mu_a\rangle|^2\mid\mathcal F_a)
 \ge\lambda_n\|A\|_{\HS}^2=\lambda_n\|u\|^2.
\]
Hence $Q_a\ge\lambda_n I_{r^2}$ on $\mathcal G_a$.

If $H=0$, all blocks $A_i$ already lie in $\Sigma$, and
\eqref{eq:struct-column-variance} directly gives
\[
 \E(|\langle A,F_a-\mu_a\rangle|^2\mid\mathcal F_a)
 \ge\frac{s_R^2}{n}\|A\|_{\HS}^2.
\]
The same conclusion therefore holds with $d_*=0$,
$\lambda_n=s_R^2/n$, and
$\mathcal G_a=\{\operatorname{rank}T_{-a}=(n-1)r\}$, an
$\mathcal F_a$-measurable event of probability one by Step 1.

\textit{Step 3: from the projected column to $T_D$.}
If $\sg(T_D)\le t$, choose a unit vector $v=(v_1,\ldots,v_n)$
with $\|T_Dv\|\le t$. Some block satisfies $\|v_a\|\ge n^{-1/2}$.
Since projection onto $W_a$ annihilates the other block columns,
\begin{equation}\label{eq:struct-5.4}
 \|F_av_a\|=\|P_{W_a}T_av_a\|=\|P_{W_a}T_Dv\|\le t,
 \qquad\sg(F_a)\le\sqrt n\,t.
\end{equation}
Conditionally on $\mathcal F_a$, Step 2 verifies the hypotheses of
\cref{struct:lem:gaussian} for $F_a$ on $\mathcal G_a$. Taking
$u=\sqrt n\,t$ in \eqref{eq:struct-4.1} gives
\[
 \mathbf1_{\mathcal G_a}
 \Pp\{\sg(F_a)\le\sqrt n\,t\mid\mathcal F_a\}
 \le\frac{r^2nt^2}{\lambda_n}\le C_Rn^{2d_*+2}t^2.
\]
By \eqref{eq:struct-5.4} and a union bound over the $n$ block columns,
\begin{align*}
 \Pp\{\sg(T_D)\le t\}
 &\le\sum_{a=1}^n\Pp\{\sg(F_a)\le\sqrt n\,t\}\\*
 &\le\sum_{a=1}^n\left(
   \E\!\left[\mathbf1_{\mathcal G_a}
      \Pp\{\sg(F_a)\le\sqrt n\,t\mid\mathcal F_a\}\right]
     +\Pp(\mathcal G_a^c)\right)\\*
 &\le C_Rn^{2d_*+3}t^2+Cn e^{-cn}
 \le C_Rn^{2d_*+3}t^2+C'e^{-c'n}.
\end{align*}
This proves \eqref{eq:struct-5.5}, including the case $H=0$.
\end{proof}

\subsection{Upper triangular reduction and completion}\label{subsec:triangular-reduction}
We first treat $\lambda=1$, following the composition-series
reduction of \cite[Section 3.1]{Han}; the explicit constants in
\eqref{eq:struct-gap-exponents} and \eqref{eq:struct-column-constant}
will then give the statement for general $\lambda$.
\begin{proof}[Proof of \cref{struct:thm:main} for $\lambda=1$]
\textit{Step 1: the irreducible diagonal blocks.}
Choose a composition series for the unital algebra generated by the
coefficient matrices,
\[
 0=V_0\subsetneq V_1\subsetneq\cdots\subsetneq V_s=\C^m,
 \qquad s\le m,
\]
and set $E_j=V_j\cap V_{j-1}^\perp$. The quotient $V_j/V_{j-1}$
identifies isometrically with $E_j$ by $[v]\mapsto P_{E_j}v$.
Its coefficient action is represented by
\[
 R_\ell^{(j)}=P_{E_j}R_\ell|_{E_j},\qquad
 \|R_\ell^{(j)}\|\le\|R_\ell\|.
\]
Although $E_j$ need not be invariant, the invariance of $V_{j-1}$
and $V_j$ makes this compression the quotient action.

In orthonormal bases adapted to
$\C^m=E_1\oplus\cdots\oplus E_s$, all $R_\ell$ are block upper
triangular. The corresponding unitary change of basis for $T_D$
preserves its singular values and gives diagonal blocks
\[
 (T_D)_{jj}=I_{E_j}\otimes I_n+
             \sum_\ell R_\ell^{(j)}\otimes(C_\ell+D_\ell).
\]
Each quotient action is irreducible. If all its coefficients vanish,
the diagonal block is the identity; otherwise
\cref{prop:irreducible} applies. These blocks share the Gaussian
matrices, but only a union bound will be needed.

\textit{Step 2: an inverse bound counting only nonidentity blocks.}
On $\Omega$ from \eqref{eq:struct-norm-events}, write $T_D=Q+U$,
where $Q$ is block diagonal and $U$ strictly block upper triangular.
Each diagonal block is a compression of $T_D$, so
$\|Q\|\le\|T_D\|\le L$ and $\|U\|\le2L$. The inverse expansion is
\begin{equation}\label{eq:struct-triangular-inverse}
 T_D^{-1}=\sum_{j=0}^{s-1}(-Q^{-1}U)^jQ^{-1}.
\end{equation}
Expand each block of this sum along strictly increasing sequences
of block indices. Each path visits any diagonal block at most once.
If every nonidentity diagonal block has least singular value at
least $h\le1$, a path contains at most $s_1$ factors of norm $h^{-1}$;
the identity blocks do not contribute to the power of $h^{-1}$. There are finitely many paths,
depending only on $s$. Hence, for example with
$C_{\mathrm{tri}}=s^2 2^s(2L)^{s-1}$,
\begin{equation}\label{eq:struct-triangular-bound}
 \|T_D^{-1}\|\le C_{\mathrm{tri}}h^{-s_1}.
\end{equation}
For $s_1>0$ and sufficiently small $t$, take
$h=(2C_{\mathrm{tri}}t)^{1/s_1}\le1$. Then
$\sg(T_D)\le t$ on $\Omega$ forces a nonidentity diagonal block
to have least singular value less than $h$. A union bound and
\cref{prop:irreducible} give $Cn^at^{2/s_1}+Ce^{-cn}$.
If $s_1=0$, \eqref{eq:struct-triangular-inverse} instead bounds
$\|T_D^{-1}\|$ by $C_{\mathrm{tri}}$ on $\Omega$, proving the
exponential estimate for $t<C_{\mathrm{tri}}^{-1}$. Enlarging $C$
handles the remaining $t$ when $s_1>0$.

\textit{Step 3: a dimension-only exponent.}
Set $d_*(p,0)=0$, and for $q\ge1$ use
\cref{eq:struct-gap-exponents,eq:struct-dual-exponents}. One may take
\begin{equation}\label{eq:struct-dimension-exponent}
 a(m)=\max_{1\le p\le m^2,\ 0\le q<p}\{2d_*(p,q)+3\}.
\end{equation}
This covers every nonzero irreducible quotient. The number of inputs
affects the constants and the threshold, but not $a(m)$. Since
$s_1\le m$, the exponent $2/m$ also follows, with the constant enlarged
for $t>1$. The all-zero tuple gives $T_D=I$ directly.
\end{proof}

\begin{remark}[Exponents]\label[remark]{struct:rem:exponents}
The dimension exponent is not optimized: already $2d_*+3=42123$
for $p=4,q=1$, and it is about $3.55\cdot10^{93}$ for $p=9,q=8$.
The power $2/m$ of the small-ball variable has a simple worst-case
illustration. At $n=1$, take $R=I_m+N$, with $N$ a nilpotent Jordan
block, and mean $D=-1$. Then $T=cI_m+(c-1)N$, with $c$ a complex
Gaussian. Near $c=0$ its determinant is $c^m$ and its other $m-1$
singular values stay bounded above and below. Thus
$\sg(T)\asymp|c|^m$ and its small-ball probability is of order
$t^{2/m}$. This illustrates the loss from triangular blocks; it is
not an optimality assertion about the dimension exponent or the
exponential remainder in \eqref{eq:struct-1.1}.
\end{remark}

\subsection{Uniformity under a common scalar rescaling}
\label{subsec:scalar-rescaling}
\begin{proof}[Completion of the proof of \cref{struct:thm:main}]
Apply \cref{struct:thm:dual,prop:irreducible} and the case $\lambda=1$ of
\cref{subsec:triangular-reduction} to $\lambda R$, $\omega=\omega_0/\lambda$.
The coefficient spaces and composition series do not change.
For fixed base thresholds $\sigma_*,\vartheta_*$, the active image
at threshold $\sigma$ is the base active image at $\sigma/\lambda$.
Thus all parameter dependence is given by
\begin{equation}\label{as:scaled-norms}
 \begin{gathered}
 \sigma_0=\lambda\sigma_*,\quad\vartheta=\vartheta_*,\quad
 M_1=\sqrt k(8+\omega_0/\lambda)\le C/\lambda,\\
 L=1+(8+\omega_0/\lambda)\sum_j\|\lambda R_j\|
   \le L_0:=1+(8+\omega_0)\sum_j\|R_j\|,\quad
 s_{\lambda R}=\lambda s_R.
 \end{gathered}
\end{equation}
The constant $C_{\mathrm{reg}}(p,q,L,\vartheta)$ is therefore bounded
uniformly in $\lambda$. The explicit threshold
\eqref{eq:struct-gap-exponents} is at most
$C\lambda^{-(2p)^p}$; the extra conditions in
\cref{struct:thm:dual} are uniform because $L\le L_0$.
The exponents $\ell,u_j,d_*$ do not change. The frame bound
\eqref{eq:struct-column-constant} gives
\begin{equation}\label{as:scaled-irreducible}
 \Pp\{\sg(T_{\lambda,D})\le t\}
 \le C\lambda^{-2}n^{2d_*+3}t^2+Ce^{-cn}
\end{equation}
for a nonzero irreducible base tuple. The case $H=0$ follows from
the same frame bound with $d_*=0$.

For a general tuple, \eqref{eq:struct-triangular-bound} has a constant
depending only on $s,L_0$. Applying \eqref{as:scaled-irreducible} to
its nonzero quotient actions gives \eqref{eq:struct-1.1}, and gives
the stated exponential alternative when $s_1=0$. Every quotient has
$p\le m^2$, so $n\ge C_*\lambda^{-(2m^2)^{m^2}}$ covers all
thresholds. This proves the theorem.
\end{proof}

\section{Harmonicity and invertibility of linear pencils}\label{sp:analysis}
Recall that $(\mathcal M,\tau)$ is a finite von Neumann algebra with a faithful
normal tracial state, and $s_1,\ldots,s_k$ is a free standard
semicircular family. Throughout this section the coefficient matrices
are arbitrary complex matrices; the pencil $T$ defined below need not
be selfadjoint.
For $A_1,\ldots,A_k\in M_m(\C)$, define
\begin{equation}\label{sp:pencil-def}
\begin{aligned}
 T&=I+\sum_{j=1}^k A_j\otimes s_j,&
 V&=\log\Delta_m(T)=\lim_{\eta\downarrow0}
       \frac12\tau_m\log(T^*T+\eta^2I),
\end{aligned}
\end{equation}
with the trace convention of \cref{subsec:shared-rotation}.
Irreducibility has the
meaning fixed in \cref{sec:means}; a tuple is nonzero if at least
one coefficient is nonzero.
Invertibility of $T$ is understood in $M_m(\mathcal M)$,
so its inverse is bounded.

\begin{definition}\label[definition]{sp:no-intertwiner}
A coefficient tuple $\mathbf A=(A_1,\ldots,A_k)$ satisfies the
\emph{no-intertwiner condition} if
\begin{equation}\label{sp:admissible}
 \{H\in M_m(\C):A_jH=HA_j^*\text{ for every }j\}=\{0\}.
\end{equation}
\end{definition}
This condition is open, since it is injectivity of the linear map
$H\mapsto(A_jH-HA_j^*)_{j=1}^k$. Irreducibility is open as well:
by Burnside's theorem, finitely many words span $M_m(\C)$ at an
irreducible tuple, and their linear independence persists nearby.

For a word $w=(j_1,\ldots,j_\ell)$, write
\begin{equation}\label{sp:words}
 A_w=A_{j_1}\cdots A_{j_\ell},\qquad A_\varnothing=I.
\end{equation}

The following result shows a strong rigidity property of the harmonicity outside the domain. 
\begin{theorem}\label[theorem]{sp:harmonic-disk}
Let $A_j\colon\Omega\to M_m(\C)$ be holomorphic on a connected
domain $\Omega\subset\C$. Assume that $(A_j(z))_{j=1}^k$ is nonzero
and irreducible for every $z\in\Omega$. Let $D\subset\Omega$ be a
connected open set on which $V$ is harmonic and the no-intertwiner
condition holds.

If $T$ is noninvertible at some point of $D$, then all traces
$\Tr A_w$ are constant on $\Omega$, and $T$ is noninvertible
throughout $\Omega$. In particular, if $T(z_*)$ is invertible
for some $z_*\in\Omega$, then $T$ is invertible throughout $D$.
\end{theorem}

We will show that noninvertibility on a harmonic disk forces
the traces $\Tr A_w$ of all words to be constant. These
traces determine an irreducible tuple up to a common change of basis.
The proof has two analytic cases, distinguished by the limiting
Hermitized resolvent: a positive definite limit of its negative
imaginary part, treated in
\cref{sp:positive-rigidity}, and a vanishing imaginary part at the
spectral-radius threshold, treated in \cref{sp:critical-rigidity}.
If $T$ is noninvertible somewhere in $D$, these cases lead to the
same algebraic conclusion:
\[
 \begin{gathered}
 \Tr A_w\text{ is constant on a subdisk for every word }w\\
 \Longrightarrow\quad \Tr A_w\text{ is constant on }\Omega
                  \text{ for every word }w\\
 \Longrightarrow\quad T\text{ is nowhere invertible on }\Omega.
 \end{gathered}
\]
The last implication uses \cref{sp:traces}; the existence of an
invertible value $T(z_*)$ therefore excludes noninvertibility on $D$. We first establish the algebraic reductions,
then the estimates and regularity needed for these two cases.

The three cases are distinguished by the limits of
$\widehat Y_\eta$ and $\alpha_\eta$ as $\eta\downarrow0$;
see \cref{sp:compact,sp:limit-cases}.
\begin{center}
\begin{tabular}{llll}
\toprule
Limit & $r(\Phi_{R\mathbf A})$ & Pencil $T$ & Treated in\\
\midrule
$\widehat Y>0$, $\alpha=0$ & & not invertible & \cref{sp:positive-rigidity}\\
$\widehat Y=0$, $\alpha>0$ & $<1$ & invertible & \\
$\widehat Y=0$, $\alpha=0$ & $=1$ & not invertible & \cref{sp:critical-rigidity}\\
\bottomrule
\end{tabular}
\end{center}

\subsection{Irreducible blocks and traces of words}\label{sp:algebra}
\begin{lemma}[Burnside {\cite{LR}}]\label[lemma]{lem:burnside}
Let $\mathbb K$ be an algebraically closed field and let $W$ be a
nonzero finite-dimensional vector space over $\mathbb K$.
Every irreducible unital $\mathbb K$-subalgebra of $\End_{\mathbb K}(W)$
is $\End_{\mathbb K}(W)$ itself.
\end{lemma}
We use this theorem over $\overline{\C(z)}$ in the block decomposition
below. Over $\C$, it identifies an irreducible coefficient algebra
with the full matrix algebra.

Recall the word notation $A_w$ from \eqref{sp:words}.
We call $\Tr A_w$ a \emph{trace of a word in $A_1,\ldots,A_k$}.

\begin{lemma}\label[lemma]{sp:blocks}
Let $A_j(z)\in M_m(\C(z))$, and fix bounded inputs $a_j\in\mathcal M$.
Set $T(z)=I+\sum_j A_j(z)\otimes a_j$ and $V(z)=\log\Delta_m(T(z))$.
There is a finite set $\Lambda\subset\C$
such that on every simply connected domain $\Omega\subset\C\setminus \Lambda$
one can choose a holomorphic $S\colon\Omega\to\mathrm{GL}_m(\C)$ with
\begin{align*}
 S^{-1}A_jS&=
 \begin{pmatrix}
 A_j^{(1)}&*&*\\0&\ddots&*\\0&0&A_j^{(q)}
 \end{pmatrix},&
 (S^{-1}\otimes1)T(S\otimes1)&=
 \begin{pmatrix}
 T_1&*&*\\0&\ddots&*\\0&0&T_q
 \end{pmatrix},
\end{align*}
where $T_\ell=I_{m_\ell}+\sum_jA_j^{(\ell)}\otimes a_j$.
Each diagonal tuple is irreducible at every point of $\Omega$.
Its entries are holomorphic and algebraic over $\C(z)$.
For each $\ell$, the tuple is either nonzero everywhere or
identically zero with $m_\ell=1$, in which case $T_\ell=1$. Moreover,
\begin{equation}\label{sp:block-potential}
 V=\sum_{\ell=1}^q\frac{m_\ell}{m}\log\Delta_{m_\ell}(T_\ell),
 \qquad
 T\text{ invertible}\ \Longleftrightarrow\
 T_\ell\text{ invertible for every }\ell.
\end{equation}
\end{lemma}
    \begin{proof}
The composition-filtration reduction follows Han
\cite[Section~3.1, equations (3.2)--(3.5)]{Han}.
We explain the additional specialization in the parameter $z$.

Let $F=\overline{\C(z)}$. Choose a maximal flag of common invariant
subspaces
\[
0=F_0\subsetneq F_1\subsetneq\cdots\subsetneq F_q=F^m.
\]
Maximality makes each quotient $F_\ell/F_{\ell-1}$ irreducible.
Choose a basis adapted to this flag and let
$S=[\,S_1\ \cdots\ S_q\,]\in\mathrm{GL}_m(F)$, where $S_\ell$ contains
the $m_\ell=\dim F_\ell-\dim F_{\ell-1}$ new basis vectors added
at stage $\ell$. Invariance gives
\[
A_jS_\ell=\sum_{i\leq\ell}S_iB^{(j)}_{i\ell},
\qquad
S^{-1}A_jS=
\begin{pmatrix}
A_j^{(1)}&*&\cdots&*\\
0&A_j^{(2)}&\cdots&*\\
\vdots&\ddots&\ddots&\vdots\\
0&\cdots&0&A_j^{(q)}
\end{pmatrix},
\]
where $A_j^{(\ell)}=B^{(j)}_{\ell\ell}\in M_{m_\ell}(F)$
represents the quotient action.

By Burnside's theorem (\cref{lem:burnside}), for each $\ell$
choose $m_\ell^2$ words whose evaluations on the $\ell$th block
form a basis of $M_{m_\ell}(F)$. The determinant of their vectorized columns is a nonzero element $d_\ell\in F$. Its nonvanishing will
certify irreducibility after specialization.

The entries of $S$, $S^{-1}$, and the diagonal blocks lie in a
finite extension $K=\C(z)(\theta)$. Let $p(z,t)\in\C[z,t]$ be the minimal polynomial of $\theta$ after clearing denominators.
Exclude the poles of the original coefficients and the zeros
of the leading coefficient and discriminant of $p$ in $t$.
Every needed entry and each $d_\ell$ can be written as
\[
f=\frac{U_f(z,\theta)}{Q_f(z,\theta)},
\qquad U_f,Q_f\in\C[z,t],\qquad Q_f(z,\theta)\ne0\text{ in }K.
\]
For every denominator, every $d_\ell$, and one nonzero entry from
each nonzero diagonal tuple, take its minimal polynomial over
$\C(z)$ and clear denominators. Its constant coefficient is a
nonzero polynomial. Add the zeros of these constant coefficients
to the exceptional set $\Lambda$, which remains finite.

On a simply connected $\Omega\subset\C\setminus \Lambda$, the simple
roots of $p$ admit holomorphic branches. Choose one branch
$\theta(z)$ and use it in every expression above. Polynomial
identities satisfied by $\theta$ remain valid after substitution.
Indeed, their polynomial numerators are divisible by $p$ over
$\C(z)$, so they vanish after substitution wherever the coefficients
are defined, and then everywhere by holomorphic continuation.
In particular, each evaluated denominator satisfies its chosen
polynomial equation; its nonvanishing constant coefficient prevents
it from being zero. Thus all entries are holomorphic and algebraic
on $\Omega$.

The identities defining $S^{-1}$ and the triangular form remain valid after substitution. Thus $S(z)$ is invertible and simultaneously puts
all $A_j(z)$ into the asserted block form throughout $\Omega$.
Moreover, $d_\ell(z)\ne0$, so the chosen word matrices span
$M_{m_\ell}(\C)$ and each diagonal tuple remains irreducible.
For $m_\ell>1$, irreducibility excludes the zero tuple.
For $m_\ell=1$, the selected nonzero coefficient stays nonzero
throughout $\Omega$, unless all coefficients vanish identically,
in which case $T_\ell=1$.

Conjugating $T$ by $S\otimes1$ gives the corresponding upper triangular
operator pencil. Similarity invariance and the triangular formula
for the Fuglede--Kadison determinant in \cref{sp:det-facts} yield
\[
V=\sum_{\ell=1}^q\frac{m_\ell}{m}
  \log\Delta_{m_\ell}(T_\ell).
\]
For invertibility, first consider
\[
U=\begin{pmatrix}D_1&C\\0&D_2\end{pmatrix}.
\]
If $U$ is invertible, the lower-right block of $UU^{-1}=I$ gives
a right inverse for $D_2$. A right-invertible element of a finite
von Neumann algebra is invertible. Block elimination therefore gives
\[
\begin{pmatrix}I&-CD_2^{-1}\\0&I\end{pmatrix}
U
=
\begin{pmatrix}D_1&0\\0&D_2\end{pmatrix}.
\]
The elimination factor is invertible, with inverse
$\begin{pmatrix}I&CD_2^{-1}\\0&I\end{pmatrix}$.
Thus the block diagonal matrix is invertible, and so $D_1$ is
invertible. Conversely, if both diagonal blocks are invertible,
the same identity expresses $U$ as a product of invertible matrices.
Induction, grouping the first $q-1$ blocks together and separating
the last block, proves the invertibility equivalence in
\eqref{sp:block-potential}.
\end{proof}

\begin{remark}\label[remark]{sp:joining-domain}
Given a disk with closure outside the finite set $\Lambda$ and a point $z_*\notin \Lambda$,
one can join them by a simple path avoiding $\Lambda$ and take a simply
connected neighborhood of their union. The same holomorphic blocks
are then available near both the disk and $z_*$.
\end{remark}

\begin{lemma}\label[lemma]{sp:traces}
Let $(A_j(z))_{j=1}^k$ be holomorphic and irreducible on a connected
domain $\Omega\subset\C$. Suppose that $\Tr A_w(z)$ is constant
on a nonempty open subset for every word $w$. For any $z_0,z\in\Omega$
there is $S_z\in\mathrm{GL}_m(\C)$ such that
\[
 A_j(z)=S_zA_j(z_0)S_z^{-1}\qquad(1\le j\le k).
\]
The same matrix $S_z$ conjugates every coefficient. Consequently, invertibility of
$I+\sum_jA_j(z)\otimes a_j$ is independent of $z$, for any fixed
bounded operators $a_j$.
\end{lemma}
\begin{proof}
For every word $w$, the function $z\mapsto\Tr A_w(z)$ is holomorphic.
The identity theorem therefore makes it constant throughout $\Omega$.
Fix $z_0,z\in\Omega$. By Procesi's theorem
\cite[Theorem~1.3]{Procesi}, traces of words generate the polynomial
invariants of simultaneous conjugation on $M_m(\C)^k$. Thus the
tuples at $z_0$ and $z$ have the same polynomial invariants.

Each tuple defines a representation of the free algebra
$\C\langle x_1,\ldots,x_k\rangle$ on $\C^m$, with $x_j$ acting
by $A_j$. Irreducibility makes this representation simple, hence
semisimple, so its simultaneous-conjugation orbit is Zariski closed.
The invariant quotient separates closed orbits
\cite[pp.~245--246]{LeBruyn}. Consequently, the two tuples belong
to the same orbit. There is therefore a single $S_z\in\mathrm{GL}_m(\C)$
such that
\[
A_j(z)=S_zA_j(z_0)S_z^{-1}
\qquad(1\le j\le k).
\]

For fixed bounded inputs $a_j$, it follows that
\[
I+\sum_jA_j(z)\otimes a_j
=(S_z\otimes1)
 \left(I+\sum_jA_j(z_0)\otimes a_j\right)
 (S_z^{-1}\otimes1).
\]
Thus invertibility is independent of $z$. Only pointwise similarity
is needed; no holomorphic choice of $S_z$ is required.
\end{proof}

\subsection{Hermitized Dyson equation and potential identities}\label{sp:dyson}
Hermitization replaces a nonselfadjoint pencil by a selfadjoint block
operator, whose operator-valued Cauchy transform satisfies a Dyson
equation. For the general operator-valued framework see
\cite{VoiculescuOV,SpeicherOV}; matrix and operator-valued stability
theory is developed in \cite{AjankiEK,AltEKDyson}, and the Hermitization
of non-Hermitian Kronecker matrices is treated in \cite{AEKN}.

In this subsection we allow an arbitrary matrix constant term:
\begin{equation}\label{sp:pencil-u}
 T(u)=K(u)\otimes1+\sum_j A_j(u)\otimes s_j,
 \qquad K,A_j:\mathcal O\to M_m(\C),\quad \mathcal O\subset\R^p\text{ open}.
\end{equation}
No invertibility of $K$ is assumed. The subsequent support analysis
will specialize to $K=I$.
In \cref{sp:rigidity} we return to holomorphic coefficients
$A_j(z)$ and take $u=(\Re z,\Im z)$. For each $u$, define the
linear maps $\Phi_u,\Theta_u:M_m(\C)\to M_m(\C)$ and their
Hilbert--Schmidt adjoints by
\begin{equation}\label{sp:covariance-maps}
\begin{aligned}
 \Phi_u(H)&=\sum_jA_j(u)HA_j(u)^*,&
 \Phi_u^*(H)&=\sum_jA_j(u)^*HA_j(u),\\
 \Theta_u(H)&=\sum_jA_j(u)HA_j(u),&
 \Theta_u^*(H)&=\sum_jA_j(u)^*HA_j(u)^*.
\end{aligned}
\end{equation}
The stars on these maps denote Hilbert--Schmidt adjoints;
in particular $\Tr(\Phi_u(H)K)=\Tr(H\Phi_u^*(K))$.
For $\eta>0$, put
$V_\eta(u)=\tfrac12\tau_m\log(T(u)^*T(u)+\eta^2I)\in\R$ and
define the matrix-valued functions $X_\eta,Y_\eta,R_\eta:\mathcal O\to M_m(\C)$ by
\begin{equation}\label{sp:XYR-def}
\begin{aligned}
 X_\eta(u)&=(\id\otimes\tau)[\eta(T(u)^*T(u)+\eta^2I)^{-1}],\\
 Y_\eta(u)&=(\id\otimes\tau)[\eta(T(u)T(u)^*+\eta^2I)^{-1}],\\
 R_\eta(u)&=(\id\otimes\tau)[(T(u)^*T(u)+\eta^2I)^{-1}T(u)^*].
\end{aligned}
\end{equation}
We have $X_\eta,Y_\eta>0$ and
\begin{equation}\label{sp:trace-balance}
 \Tr X_\eta=\Tr Y_\eta,
\end{equation}
because $T^*T$ and $TT^*$ have the same spectral distribution under the
finite trace, so their regularized inverse traces agree.
Write $\mathcal H=\left(\begin{smallmatrix}0&T\\T^*&0\end{smallmatrix}\right)$,
$J_K=\left(\begin{smallmatrix}0&K\\K^*&0\end{smallmatrix}\right)$, and define
$M_\eta:\mathcal O\to M_{2m}(\C)$ and, at each $u$, the covariance map
$\mathcal S:M_{2m}(\C)\to M_{2m}(\C)$ by
\begin{equation}\label{sp:M-def}
\begin{aligned}
 M_\eta&=(\id\otimes\tau)[(i\eta I-\mathcal H)^{-1}]
       =\begin{pmatrix}-iY_\eta&-R_\eta^*\\-R_\eta&-iX_\eta\end{pmatrix},\\
 \mathcal S(H)&=
 \begin{pmatrix}
 \Phi_u(H_{22})&\Theta_u(H_{21})\\
 \Theta_u^*(H_{12})&\Phi_u^*(H_{11})
 \end{pmatrix}.
\end{aligned}
\end{equation}
In particular, the covariance map has the representation
\begin{equation}\label{sp:calA-def}
 \mathcal A_j=\begin{pmatrix}0&A_j\\A_j^*&0\end{pmatrix},\qquad
 \mathcal H-J_K=\sum_j\mathcal A_j\otimes s_j,\qquad
 \mathcal S(H)=\sum_j\mathcal A_jH\mathcal A_j.
\end{equation}
We note that the representation $\mathcal S(H)=\sum_j\mathcal A_jH\mathcal A_j$ gives the bilinear trace symmetry
\begin{equation}\label{sp:covariance-trace}\Tr_{2m}(A\mathcal S(B))=\Tr_{2m}(B\mathcal S(A)).\end{equation}
We suppress $u$ when it is fixed, and suppress $\eta$ within calculations
at fixed $\eta>0$. The traces $\Tr,\Tr_{2m}$ are unnormalized; a subscript specifies
the matrix size. The normalized traces retain the conventions of \cref{subsec:shared-rotation}.

\begin{proposition}[Dyson equation]\label[proposition]{sp:mde}
The matrix $M_\eta$ is the unique solution with $\Im M<0$ of
\begin{equation}\label{sp:mde-eq}
 M^{-1}=i\eta I-J_K-\mathcal S(M).
\end{equation}
For $C=\eta I+\Phi(X)$ and $B=K-\Theta(R)$, its block form is
\begin{equation}\label{sp:mde-blocks}
 M^{-1}=\begin{pmatrix}
 iC&-B\\-B^*&i(\eta I+\Phi^*(Y))
 \end{pmatrix}.
\end{equation}
Block inversion therefore gives
\begin{equation}\label{sp:schur}
 R=XB^*C^{-1},\qquad Y=C^{-1}-R^*X^{-1}R,\qquad
 X^{-1}=\eta I+\Phi^*(Y)+B^*C^{-1}B.
\end{equation}
\end{proposition}
\begin{proof}
The operator $\mathcal H-J_K$ is semicircular over $M_{2m}(\C)$,
with covariance map $\mathcal S$. Its operator-valued $R$-transform
is therefore $\mathcal S$; see \cite[Section 8]{VoiculescuOV},
\cite[Chapter 4]{SpeicherOV}, and \cite[Section 1]{HRS}.
This gives \eqref{sp:mde-eq}. For uniqueness, put $W=iM$. Then
\[
 W^{-1}=(\eta I/2+iJ_K)+(\mathcal S(W)+\eta I/2),\qquad
 \Re(\mathcal S(W)+\eta I/2)=\mathcal S(\Re W)+\eta I/2>0.
\]
The second summand is an analytic map, bounded on bounded sets,
whose real part is at least $\eta I/2$ on $\Re W>0$.
These are the hypotheses of \cite[Theorem 2.1]{HRS};
the map need not be linear or the covariance nondegenerate.
\cite[Theorem 2.1]{HRS} therefore gives uniqueness on $\Re W>0$.

Substitution of the blocks of $M$ gives \eqref{sp:mde-blocks}.
Since $C>0$, the Schur complement of its upper-left block $iC$ is
\[
 i(\eta I+\Phi^*(Y))-(-B^*)(iC)^{-1}(-B)
 =i\bigl(\eta I+\Phi^*(Y)+B^*C^{-1}B\bigr).
\]
Its inverse is the lower-right block $-iX$ of $M$, proving the
last identity in \eqref{sp:schur}. The remaining inverse blocks are
\[
 [M]_{21}=-XB^*C^{-1},\qquad
 [M]_{11}=-iC^{-1}+iC^{-1}BXB^*C^{-1}.
\]
Comparing them with $-R$ and $-iY$, respectively, gives
$R=XB^*C^{-1}$ and $Y=C^{-1}-R^*X^{-1}R$.
\end{proof}

The following determinant identity was obtained independently
in the present work. It also appears in two recent works
\cite[Theorem 3.1]{Kargin} and \cite[version 1, equation (2.8)]{YinFK}, which became available
while this manuscript was being prepared. The centered determinant-capacity formula is due to Mai--Speicher \cite{MaiSpeicher}. Our proof follows
the regularization-and-differentiation strategy of
Haagerup--Schultz \cite[Lemma 4.14]{HSunbounded} in the
$R$-diagonal setting and Zhong
\cite[Lemma 3.12]{ZhongElliptic} for elliptic operators.

We give the block calculation in the form needed for our
Brown measure algorithm. Combined with linearization, the
identity expresses the regularized logarithmic potential
directly in terms of the solution $M_\eta$ of the Dyson equation in \cref{sp:mde}.
Recall that $V_\eta=\tfrac12\tau_m\log(T^*T+\eta^2I)$.

\begin{lemma}\label[lemma]{sp:action}
For every $u$ and $\eta>0$,
\begin{equation}\label{sp:kargin-identity}
 2mV_\eta=-\log|\det M_\eta|
       +\tfrac12\Re\Tr_{2m}(M_\eta\mathcal S(M_\eta)).
\end{equation}
Equivalently, in block form,
\begin{equation}\label{sp:action-eq}
\begin{aligned}
    2mV_\eta(u) &=\log|\det M^{-1}|-\Tr(Y\Phi(X))
+\operatorname{Re}\Tr(R\Theta(R))\\
 &=\log\frac{\det(\eta I+\Phi_u(X_\eta(u)))}{\det X_\eta(u)}
       +\Re\Tr(K(u)R_\eta(u))+\eta\Tr X_\eta(u)-m.
\end{aligned}
\end{equation}
If $K(u)$ and the coefficients $A_j(u)$ are $C^1$, then, for each coordinate $u_a$,
\begin{equation}\label{sp:gradient}
 2m\partial_{u_a}V_\eta
 =\Tr\bigl(Y_\eta(\partial_{u_a}\Phi_u)(X_\eta)\bigr)
       -\Re\Tr\bigl(R_\eta(\partial_{u_a}\Theta_u)(R_\eta)\bigr)
       +2\Re\Tr(R_\eta\partial_{u_a}K).
\end{equation}
The derivatives of the maps hold their argument fixed; explicitly,
\begin{align*}
 (\partial_{u_a}\Phi_u)(H)
 &=\sum_j\bigl((\partial_{u_a}A_j)HA_j^*
                       +A_jH(\partial_{u_a}A_j)^*\bigr),\\
 (\partial_{u_a}\Theta_u)(H)
 &=\sum_j\bigl((\partial_{u_a}A_j)HA_j
                       +A_jH(\partial_{u_a}A_j)\bigr).
\end{align*}
\end{lemma}
\begin{proof}
Write $M=M_\eta$ and set
\[
 \mathcal F=\log|\det M^{-1}|
            +\tfrac12\Re\Tr_{2m}(M\mathcal S(M)).
\]
We differentiate along a real $C^1$ path in $\eta,K,A_1,\ldots,A_k$,
with $\eta>0$. A dot denotes the path derivative, and
$\dot{\mathcal S}$ differentiates the coefficients of the map with its
argument fixed. Jacobi's formula for an invertible matrix $D$ reads
$\partial_t\log|\det D|=\Re\Tr(D^{-1}\dot D)$. Applying it to
$D=M^{-1}$ and differentiating \eqref{sp:mde-eq} gives
\[
 \partial_t\log|\det M^{-1}|
 =\Re\Tr_{2m}\bigl(M(i\dot\eta I-\dot J_K
               -\dot{\mathcal S}(M)-\mathcal S(\dot M))\bigr).
\]
By the bilinear trace symmetry \eqref{sp:covariance-trace},
\[
 \tfrac12\partial_t\Tr_{2m}(M\mathcal S(M))
 =\Tr_{2m}(M\mathcal S(\dot M))
     +\tfrac12\Tr_{2m}(M\dot{\mathcal S}(M)).
\]
Adding the real parts cancels the terms containing $\dot M$ and gives
\begin{equation}\label{sp:action-differential}
 \dot{\mathcal F}
 =\Re\Tr_{2m}\bigl(M(i\dot\eta I-\dot J_K)\bigr)
       -\tfrac12\Re\Tr_{2m}\bigl(M\dot{\mathcal S}(M)\bigr).
\end{equation}

First vary only $\eta$. By \eqref{sp:M-def} and
\eqref{sp:trace-balance}, $\partial_\eta\mathcal F=\Re\Tr_{2m}(iM)
=\Tr Y+\Tr X=2\Tr X$. Tracial differentiation
\eqref{eq:trace-log-derivative} gives
\[
 \partial_\eta(2mV_\eta)
 =2\eta(\Tr\otimes\tau)\bigl((T^*T+\eta^2I)^{-1}\bigr)=2\Tr X.
\]
To determine the constant, let $\eta\to\infty$. Since $T$ is bounded,
$2mV_\eta-2m\log\eta=(\Tr\otimes\tau)\log(I+\eta^{-2}T^*T)\to0$.
The resolvent expansion $M_\eta=(i\eta)^{-1}I+O(\eta^{-2})$ gives
$\log|\det M^{-1}|=2m\log\eta+o(1)$ and
$\Tr_{2m}(M\mathcal S(M))=O(\eta^{-2})$. Hence
$2mV_\eta=\mathcal F$, proving \eqref{sp:kargin-identity}.

To obtain the block formula, multiply the blocks in \eqref{sp:M-def}.
The two diagonal products contribute $-2\Tr(Y\Phi(X))$ to
$\Tr_{2m}(M\mathcal S(M))$, while the two off-diagonal products have
conjugate traces and contribute $2\Re\Tr(R\Theta(R))$. Thus
\begin{equation}\label{eq:coupling}
 2mV_\eta=\log|\det M^{-1}|-\Tr(Y\Phi(X))
                         +\Re\Tr(R\Theta(R)).
\end{equation}
The Schur complement of $iC$ in $M^{-1}$ is $iX^{-1}$
by \eqref{sp:schur}, so $|\det M^{-1}|=\det C/\det X$.
The upper-left block of $M^{-1}M=I$ and its real trace give
\begin{equation}\label{sp:block-trace}
 \begin{aligned}
 CY+(K-\Theta(R))R&=I,\\
 m&=\eta\Tr X+\Tr(Y\Phi(X))+\Re\Tr(KR)-\Re\Tr(R\Theta(R)).
 \end{aligned}
\end{equation}
Here $C=\eta I+\Phi(X)$ and $\Tr Y=\Tr X$. Substitution in
\eqref{eq:coupling} proves \eqref{sp:action-eq}.

Finally, vary $u_a$ at fixed $\eta$ in \eqref{sp:action-differential}.
The off-diagonal blocks of $\dot J_K$ give
$-\Re\Tr_{2m}(M\dot J_K)=2\Re\Tr(R\dot K)$.
The same block multiplication as in \eqref{eq:coupling}, with
$\mathcal S$ replaced by $\dot{\mathcal S}$, gives
\[
 -\tfrac12\Re\Tr_{2m}(M\dot{\mathcal S}(M))
 =\Tr(Y\dot\Phi(X))-\Re\Tr(R\dot\Theta(R)).
\]
Together these identities give \eqref{sp:gradient}.
\end{proof}

\subsection{Regularization and its limits}\label{sp:stability}
From now until the end of \cref{sp:analysis}, return to $K=I$ and
write $J_0=J_I$. The estimates below use this normalization.
Related estimates for Dyson equations appear in
\cite[Proposition 3.2]{AEK}, \cite[Lemmas 5.2 and 5.7]{AK}, and
\cite[Lemma 3.6]{AEKN}. We need a comparison that requires neither
flatness nor irreducibility. We first prove it for $X_\eta$ and
$C_\eta=\eta I+\Phi(X_\eta)$, and then normalize these matrices
to study the limit as $\eta\downarrow0$.

\begin{lemma}\label[lemma]{sp:bounds}
For all coefficient tuples and all $\eta>0$,
\[
 C_\eta\ge(1+\sqrt m)^{-2}X_\eta,\qquad
 V_\eta\ge-\log(1+\sqrt m)-\tfrac12.
\]
For every $L>0$, there is $c_2>0$, depending only
on $m,k,L$, with $c_1=(1+\sqrt m)^{-2}$, such that
\[
 c_1X_\eta\le C_\eta\le c_2X_\eta,\qquad
 0<\frac{\eta}{\Tr X_\eta}\le\frac{c_2}{m}
\]
whenever $\max_j\|A_j\|\le L$ and $0<\eta\le1$.
In particular $V=\lim_{\eta\downarrow0}V_\eta$ is finite.
\end{lemma}
\begin{proof}
Fix $\eta>0$ and suppress its subscript.

\textit{Step 1: the lower comparison.}
By \eqref{sp:schur}, the Schur complement of $X>0$ in the
block matrix below is $C^{-1}-R^*X^{-1}R=Y>0$. Moreover,
$X^{-1}-B^*C^{-1}B=\eta I+\Phi^*(Y)>0$; conjugating by
$X^{1/2}$ gives the second inequality:
\[
 \begin{pmatrix}X&R\\R^*&C^{-1}\end{pmatrix}\ge0,
 \qquad X^{1/2}B^*C^{-1}BX^{1/2}\le I.
\]
Put $\widehat{\mathcal A}_j=\operatorname{diag}(A_j,A_j^*)$. Applying the completely positive map $H\mapsto\sum_j\widehat{\mathcal A}_jH\widehat{\mathcal A}_j^*$ to the block matrix gives
\[
 \begin{aligned}
 0
 &\leq
 \sum_j
 \begin{pmatrix}A_j&0\\0&A_j^*\end{pmatrix}
 \begin{pmatrix}X&R\\R^*&C^{-1}\end{pmatrix}
 \begin{pmatrix}A_j^*&0\\0&A_j\end{pmatrix}=
 \begin{pmatrix}
 \Phi(X)&\Theta(R)\\
 \Theta(R)^*&\Phi^*(C^{-1})
 \end{pmatrix}.
 \end{aligned}
\]
Finally, conjugating by $D=\begin{pmatrix}C^{-1/2}&0\\0&X^{1/2}\end{pmatrix}$
preserves positivity. Thus
\[
 \begin{pmatrix}
 C^{-1/2}\Phi(X)C^{-1/2}
 &
 C^{-1/2}\Theta(R)X^{1/2}\\[1mm]
 X^{1/2}\Theta(R)^*C^{-1/2}
 &
 X^{1/2}\Phi^*(C^{-1})X^{1/2}
 \end{pmatrix}\geq0.
\]

The upper-left block is $I-\eta C^{-1}\le I$ as $C=\eta I+\Phi(X)$. The lower-right
block is positive and has trace
\[
 \Tr(X\Phi^*(C^{-1}))=\Tr(C^{-1}\Phi(X))
 =m-\eta\Tr C^{-1}\le m,
\]
so it is at most $mI$. Enlarging these diagonal blocks and taking
the Schur complement of $I$ yields
\[
 \|C^{-1/2}\Theta(R)X^{1/2}\|\le\sqrt m,
 \qquad \|C^{-1/2}BX^{1/2}\|\le1.
\]
Since $I=B+\Theta(R)$,
\[
 \|C^{-1/2}X^{1/2}\|\le1+\sqrt m,
 \qquad C\ge c_1X,\quad c_1=(1+\sqrt m)^{-2}.
\]

\textit{Step 2: the sign of the trace term.}
The inequality $|\Tr Z^2|\le\|Z\|_{\HS}^2$ gives
\[
 |\Tr(R\Theta(R))|\le\sum_j\|X^{-1/2}RA_jX^{1/2}\|_{\HS}^2
 =\sum_j\Tr(A_jXA_j^*R^*X^{-1}R)=\Tr(\Phi(X)R^*X^{-1}R).
\]
The adjoint of $R=XB^*C^{-1}$ implies
$I=CR^*X^{-1}+\Theta(R)$. Multiplication on the right by $R$ gives
\begin{align*}
 \Re\Tr R
 &=\eta\Tr(R^*X^{-1}R)+\Tr(\Phi(X)R^*X^{-1}R)
       +\Re\Tr(R\Theta(R))\\
 &\ge\eta\Tr(R^*X^{-1}R)\ge0.
\end{align*}

\textit{Step 3: the lower bound for the potential.}
Put $\Gamma=X^{-1/2}CX^{-1/2}$. Then
\[
 \Gamma\ge c_1I,\qquad \det \Gamma=\frac{\det C}{\det X}.
\]
By \eqref{sp:action-eq} and Step 2,
\[
 2mV_\eta=\log\det \Gamma+\Re\Tr R+\eta\Tr X-m
 \ge m\log c_1-m.
\]
Dividing by $2m$ proves the asserted lower bound.

\textit{Step 4: the upper comparison and finiteness.}
Under the coefficient bound and $\eta\le1$,
\[
 \|T\|\le1+2kL,
 \qquad V_\eta\le\tfrac12\log\bigl((1+2kL)^2+1\bigr).
\]
The nonnegative terms in Step 3 therefore imply
$\det \Gamma\le K_{m,k,L}$ for a fixed constant $K_{m,k,L}$.
Writing $\lambda_1,\ldots,\lambda_m$ for the eigenvalues of $\Gamma$,
\[
 \lambda_i\ge c_1,\qquad
 \lambda_i=\frac{\det \Gamma}{\prod_{j\ne i}\lambda_j}
 \le K_{m,k,L}c_1^{-(m-1)}=:c_2.
\]
Thus $\Gamma\le c_2I$, or $C\le c_2X$. Finally,
\[
 m\eta\le\Tr C\le c_2\Tr X
\]
gives the remaining bound. Since $V_\eta$ decreases as
$\eta\downarrow0$, its uniform lower bound and the displayed upper
bound show that its limit is a finite real number.
\end{proof}

When $X_\eta,Y_\eta\to0$, their relative sizes still distinguish the
invertible and critical cases (\cref{sp:limit-cases}). We normalize
$X_\eta$ to trace one and retain its scale:
\begin{equation}\label{sp:gamma-def}
 \gamma_\eta=\Tr X_\eta=\Tr Y_\eta,\qquad
 x_\eta=X_\eta/\gamma_\eta,\qquad \Tr x_\eta=1,
\end{equation}
and set
\begin{equation}\label{sp:normalize}
 \alpha_\eta=\eta/\gamma_\eta,\qquad
 C_{0,\eta}=C_\eta/\gamma_\eta=\alpha_\eta I+\Phi(x_\eta),\qquad
 \widehat Y_\eta=\gamma_\eta Y_\eta.
\end{equation}
In particular,
\begin{equation}\label{sp:normalization-traces}
 \Tr\widehat Y_\eta=\gamma_\eta^2,\qquad
 \alpha_\eta^2\Tr\widehat Y_\eta=\eta^2,
 \qquad \eta\gamma_\eta=\alpha_\eta\Tr\widehat Y_\eta.
\end{equation}
\Cref{sp:bounds} gives $c_1x_\eta\le C_{0,\eta}\le c_2x_\eta$ and
$0<\alpha_\eta\le c_2/m$. The Schur identities become
\begin{equation}\label{sp:normalized-eq}
 \begin{gathered}
 R=xB^*C_0^{-1},\qquad
 \widehat Y=C_0^{-1}-R^*x^{-1}R,\\
 x^{-1}-B^*C_0^{-1}B=\Phi^*(\widehat Y)+\alpha\Tr(\widehat Y)I.
 \end{gathered}
\end{equation}
For example, the last line follows by multiplying the last identity
of \eqref{sp:schur} by $\gamma_\eta$ and using
$\gamma_\eta\eta=\alpha_\eta\Tr\widehat Y_\eta$.

We now use irreducibility to prevent $x_\eta$ from approaching a
singular trace-one matrix, and to exclude a nonzero singular limit
of $\widehat Y_\eta$.
\begin{proposition}\label[proposition]{sp:compact}
Let $\eta_n\downarrow0$ and $A_j^{(n)}\to A_j$, where $(A_j)$ is
irreducible. The corresponding tuples
$(x_{\eta_n},\alpha_{\eta_n},\widehat Y_{\eta_n},R_{\eta_n})$
have convergent subsequences. Every limit $(x,\alpha,\widehat Y,R)$
satisfies $x>0$, $C_0=\alpha I+\Phi(x)>0$, and exactly one of
\[
 \alpha=0,\quad\widehat Y>0;\qquad\text{or}\qquad \widehat Y=0.
\]
In the first case $X=\sqrt{\Tr\widehat Y}\,x$ and
$Y=\widehat Y/\sqrt{\Tr\widehat Y}$ satisfy $X,Y>0$,
$\Tr X=\Tr Y$, and the Dyson equation \eqref{sp:mde-eq} at $\eta=0$.
\end{proposition}
\begin{proof}
\textit{Step 1: compactness and positive definiteness of $x$.}
Choose $L$ bounding every $\|A_j^{(n)}\|$, which is possible by
convergence, and discard finitely many terms so that $\eta_n\le1$.
Thus \cref{sp:bounds} supplies common comparison constants.
Since $x_\eta\ge0$, $\Tr x_\eta=1$, and $\alpha_\eta$ is bounded,
first take a convergent subsequence of $(x_\eta,\alpha_\eta)$.
Its limit satisfies
\[
 c_1x\le C_0=\alpha I+\Phi(x)\le c_2x,\qquad \Tr x=1.
\]
If $\alpha>0$, then $x\ge(\alpha/c_2)I>0$. If $\alpha=0$, the
two-sided comparison gives $\ran\Phi(x)=\ran x$. To identify this
range, observe that
\[
 \langle v,\Phi(x)v\rangle
 =\sum_j\|x^{1/2}A_j^*v\|^2.
\]
Since all summands are nonnegative,
$\ker\Phi(x)=\bigcap_j\ker(x^{1/2}A_j^*)$.
Taking orthogonal complements in finite dimension gives
\[
 \ran\Phi(x)=\Span_j\ran(A_jx^{1/2})
            =\Span_j A_j\ran x.
\]
Together with $\ran\Phi(x)=\ran x$, this makes $\ran x$ a
nonzero common invariant subspace. Irreducibility
forces $\ran x=\C^m$, so again $x,C_0>0$.

For the approximating matrices, $x_\eta\le I$ and
\eqref{sp:normalized-eq} imply
\[
 R_\eta^*R_\eta\le R_\eta^*x_\eta^{-1}R_\eta
 \le C_{0,\eta}^{-1},\qquad
 0\le\widehat Y_\eta\le C_{0,\eta}^{-1}.
\]
Since $C_{0,\eta}\to C_0>0$, these bounds allow a further convergent
subsequence of $(\widehat Y_\eta,R_\eta)$.
If $\alpha>0$, \eqref{sp:normalization-traces} forces
$\Tr\widehat Y=0$, hence $\widehat Y=0$.

\textit{Step 2: the kernel argument when $\alpha=0$.}
Set
\[
 \Sigma=x^{-1}-B^*C_0^{-1}B=\Phi^*(\widehat Y),\qquad
 N=\ker \Sigma,\qquad F=\ker\widehat Y.
\]
The two Schur complements of
$\mathcal B=\left(\begin{smallmatrix}x^{-1}&B^*\\B&C_0\end{smallmatrix}\right)$
are $\Sigma$ and $C_0-BxB^*=C_0\widehat Y C_0$. Elimination using either
invertible diagonal block gives
\[
 \operatorname{rank}\mathcal B
 =m+\operatorname{rank}\Sigma
 =m+\operatorname{rank}(C_0\widehat Y C_0).
\]
Consequently $\dim N=\dim F$. The limiting Schur identities give
\[
 \Sigma R=B^*\widehat Y,\qquad
 \widehat YB=C_0^{-1}Bx\Sigma,
\]
so $RF\subset N$ and $BN\subset F$. Moreover,
\[
 v\in N\quad\Longrightarrow\quad
 0=\langle v,\Sigma v\rangle=\sum_j\|\widehat Y^{1/2}A_jv\|^2
 \quad\Longrightarrow\quad A_jv\in F.
\]
Thus $A_jN\subset F$, $RA_jN\subset N$, and
$\Theta(R)N=\sum_jA_jRA_jN\subset F$. Since $I=B+\Theta(R)$,
we obtain $N\subset F$. Equal dimensions imply $N=F$, which is
then invariant under every $A_j$. Irreducibility gives $F=0$ or
$F=\C^m$, that is, $\widehat Y>0$ or $\widehat Y=0$.

\textit{Step 3: recover the unnormalized positive limit.}
If $\widehat Y>0$, then
$\gamma_\eta^2=\Tr\widehat Y_\eta\to\Tr\widehat Y>0$.
Hence $X_\eta=\gamma_\eta x_\eta$ and
$Y_\eta=\widehat Y_\eta/\gamma_\eta$ converge to the stated $X,Y$.
Writing $M$ for the corresponding limit of $M_\eta$, pass to the
limit in the product form of the Dyson equation:
\[
 (i\eta I-J_0-\mathcal S_\eta(M_\eta))M_\eta=I
 \quad\Longrightarrow\quad (-J_0-\mathcal S(M))M=I.
\]
Here $\mathcal S_\eta$ denotes the covariance map for the approximating
coefficient tuple. The last identity proves invertibility of $M$
and the zero-regularization Dyson equation.
\end{proof}

\subsection{Positive maps and invertibility}\label{sp:positive-maps}
We characterize invertibility of $T$ by a finite-dimensional equation
for $R$ and the strict inequality $r(\Phi_{R\mathbf A})<1$.
This will distinguish the two possibilities left open when
$\widehat Y=0$: $\alpha>0$ corresponds to invertibility, whereas
$\alpha=0$ is the critical case $r(\Phi_{R\mathbf A})=1$; see
\cref{sp:limit-cases}. The equality case of the spectral-radius
estimate below will also control uniqueness at that threshold.

For a linear map $\Psi$ on a matrix algebra, $r(\Psi)$ denotes its
spectral radius. For a tuple $B=(B_1,\ldots,B_k)\in M_m(\C)^k$,
put $\Theta_{\mathbf B}(H)=\sum_j B_jHB_j$, and define maps
$\Phi_{\mathbf B},\Phi_{\mathbf B}^*:M_m(\C)\to M_m(\C)$ by
$\Phi_{\mathbf B}(H)=\sum_jB_jHB_j^*$ and $\Phi_{\mathbf B}^*(H)=\sum_jB_j^*HB_j$.
In particular, $\Phi=\Phi_{\mathbf A}$ for $\mathbf A=(A_1,\ldots,A_k)$, and
$\Phi_u=\Phi_{\mathbf A(u)}$ when the coefficients depend on $u$.

For $R\in\mathrm{GL}_m(\C)$, write
$RA=(RA_1,\ldots,RA_k)$ and $AR=(A_1R,\ldots,A_kR)$. Then
$\Phi_{R\mathbf A}(H)=R\Phi_{\mathbf A}(H)R^*$ and
$\Phi_{\mathbf AR}^*(H)=R^*\Phi_{\mathbf A}^*(H)R$.

\begin{lemma}\label[lemma]{sp:positive-norm}
If $\Psi$ is completely positive on a matrix algebra, then
$\|\Psi\|=\|\Psi(I)\|$. If $X>0$ and $\Psi(X)\le aX$, then
$r(\Psi)\le a$. If $\Psi(X)<X$, then $r(\Psi)<1$; if
$\Psi(X)=aX$, then $r(\Psi)=a$.
\end{lemma}
\begin{proof}
For $\|H\|\le1$, apply $\Psi$ entrywise to the positive block
matrix $\left(\begin{smallmatrix}I&H\\H^*&I\end{smallmatrix}\right)$.
Block positivity gives $\|\Psi(H)\|\le\|\Psi(I)\|$; taking $H=I$
gives equality. Apply this to the similar map
$H\mapsto X^{-1/2}\Psi(X^{1/2}HX^{1/2})X^{-1/2}$.
The strict case follows by compactness, and in the equality case
$a$ is itself an eigenvalue.
\end{proof}

\begin{lemma}\label[lemma]{sp:wick-inverse}
Let $K,A_1,\ldots,A_k\in M_m(\C)$ and set
$T_K=K\otimes1+\sum_jA_j\otimes s_j$.
Suppose $R\in\mathrm{GL}_m(\C)$ satisfies
\[
 R^{-1}=K-\Theta(R),\qquad r(\Phi_{R\mathbf A})<1.
\]
For each word $w=i_1\cdots i_n$, let $U_w$ be its free semicircular
Wick polynomial, with $U_\varnothing=1$, and put
\begin{equation}\label{sp:wick-coefficients}
 C_\varnothing=R,\qquad
 C_{i_1\cdots i_n}=(-1)^nRA_{i_1}RA_{i_2}\cdots RA_{i_n}R.
\end{equation}
Then $T_K$ is boundedly invertible and
\[
 T_K^{-1}=\mathcal X=\sum_w C_w\otimes U_w
\]
with convergence in $L^2(M_m(\mathcal M),\tau_m)$. Moreover,
$(\id\otimes\tau)(T_K^{-1})=R$.
\end{lemma}
\begin{proof}
Work in $\mathcal M_0=W^*(s_1,\ldots,s_k)$ with its restricted trace.
The Wick polynomials form an orthonormal basis of
$L^2(\mathcal M_0,\tau)$; for a nonempty word $w=i_1\cdots i_n$,
\[
 s_jU_w=U_{jw}+\delta_{j,i_1}U_{i_2\cdots i_n}.
\]
For the empty word, $s_jU_\varnothing=U_j=s_j$.
These identities follow from the creation and annihilation
representation on full Fock space; see \cite[Section 2]{Miyagawa}.
The coefficients satisfy
\[
 C_{jw}=-RA_jC_w,\qquad
 \sum_{|w|=n}C_wC_w^*=\Phi_{R\mathbf A}^{\,n}(RR^*).
\]
Choose $r(\Phi_{R\mathbf A})<\theta<1$. Finite dimensionality gives
$\|\Phi_{R\mathbf A}^{\,n}\|\le C_0\theta^n$, and hence
\[
 \sum_{|w|=n}\frac1m\Tr(C_w^*C_w)
 =\frac1m\Tr\bigl(\Phi_{R\mathbf A}^{\,n}(RR^*)\bigr)
 \le C_0\|RR^*\|\theta^n.
\]
Summation and Wick orthogonality prove convergence of the series
to $\mathcal X\in L^2(M_m(\mathcal M_0),\tau_m)$.

Left multiplication by the bounded pencil $T_K$ is continuous on
matrix $L^2$. We can therefore multiply finite partial sums and
pass to the limit. Each fixed Wick coefficient involves
only finitely many summands. The constant coefficient of $T_K\mathcal X$
is
\[
 KC_\varnothing+\sum_jA_jC_j
 =(K-\Theta(R))R=I_m.
\]
For $w=i_1\cdots i_n\ne\varnothing$, write
$w'=i_2\cdots i_n$. Its coefficient is
\[
 KC_w+A_{i_1}C_{w'}+\sum_jA_jC_{jw}
 =(K-\Theta(R))C_w+A_{i_1}C_{w'}=0.
\]
Thus $T_K\mathcal X=I$ in $L^2$, and hence in the algebra of
closed operators affiliated with $M_m(\mathcal M_0)$.
The affiliated algebra is stably finite
\cite[arXiv version 2, Lemma 5.3]{MSY23}, so $\mathcal X T_K=I$ as well.
Consequently $\mathcal X$ is the affiliated inverse of $T_K$.

Let $D(s)$ be the division closure of
$\C\langle s_1,\ldots,s_k\rangle$ in the affiliated algebra of
$\mathcal M_0$. For a free semicircular family it is a skew field
\cite{MSY23}; see also \cite[Section 2]{Miyagawa}.
A division closure that is a skew field is also rationally closed:
matrices over it that are invertible in the ambient algebra have
inverses over it \cite[arXiv version 2, Proposition 4.9 and Lemma 4.10]{MSY23}.
Since $T_K\in M_m(D(s))$, this gives $\mathcal X\in M_m(D(s))$.
In addition,
\[
 \|\mathcal X\|_{2,\tau_m}^2
 =\frac1m\sum_{a,b}\|\mathcal X_{ab}\|_{2,\tau}^2<\infty.
\]
Each entry therefore lies in $D(s)\cap L^2(\mathcal M_0,\tau)$.
By \cite[Theorem 1.1]{Miyagawa}, every such entry
is bounded. Hence $\mathcal X\in M_m(\mathcal M_0)$ and the two
inverse identities hold for bounded operators.

Finally, all nonempty Wick words have trace zero. Continuity of
conditional expectation on $L^2$ gives
$(\id\otimes\tau)(\mathcal X)=C_\varnothing=R$.
\end{proof}

\begin{proposition}\label[proposition]{sp:invertibility}
Let $K,A_1,\ldots,A_k\in M_m(\C)$ and set
$T_K=K\otimes1+\sum_jA_j\otimes s_j$.
Then $T_K$ is invertible if and only if there is
$R\in\mathrm{GL}_m(\C)$ such that
\begin{equation}\label{sp:invertibility-eq}
 R^{-1}=K-\Theta(R),\qquad r(\Phi_{R\mathbf A})<1.
\end{equation}
Every such $R$ equals $(\id\otimes\tau)(T_K^{-1})$.
\end{proposition}
The criterion requires no irreducibility assumption and allows a
singular constant matrix $K$. The case $K=I_m$ applies to $T$ in \eqref{sp:pencil-def}.
\begin{proof}
Suppose first that $T_K$ is boundedly invertible. Define
\begin{equation}\label{sp:inv-XYR}
\begin{aligned}
X_\eta&=(\id\otimes\tau)[\eta(T_K^*T_K+\eta^2I)^{-1}],\\
Y_\eta&=(\id\otimes\tau)[\eta(T_KT_K^*+\eta^2I)^{-1}],\\
R_\eta&=(\id\otimes\tau)[(T_K^*T_K+\eta^2I)^{-1}T_K^*],
\end{aligned}
\end{equation}
and put
\begin{equation}\label{sp:inv-Q}
Q_{\mathrm R}=(\id\otimes\tau)[(T_K^*T_K)^{-1}],\quad
Q_{\mathrm L}=(\id\otimes\tau)[(T_KT_K^*)^{-1}],\quad
R=(\id\otimes\tau)(T_K^{-1}).
\end{equation}
The matrices $Q_{\mathrm R}$ and $Q_{\mathrm L}$ are positive definite.
Put $H=T_K^*T_K$. Since $H$ is positive and boundedly
invertible, for $\eta^2\|H^{-1}\|<1$ the Neumann series converges
in operator norm and gives
\[
 (H+\eta^2I)^{-1}
 =H^{-1}(I+\eta^2H^{-1})^{-1}
 =H^{-1}-\eta^2H^{-2}+O(\eta^4).
\]
Multiplying by $\eta$ and applying the bounded conditional
expectation yields
$X_\eta=\eta Q_{\mathrm R}+O(\eta^3)$.
The same computation for $T_KT_K^*$, and multiplication of the
first expansion by $T_K^*$, give
\[
 X_\eta=\eta Q_{\mathrm R}+O(\eta^3),\qquad
 Y_\eta=\eta Q_{\mathrm L}+O(\eta^3),\qquad
 R_\eta=R+O(\eta^2).
\]
All remainders here are in matrix operator norm. More precisely, if
$\|T_K^{-1}\|\le C$, then
\begin{equation}\label{sp:resolvent-remainders}
 \|X_\eta-\eta Q_{\mathrm R}\|,\ \|Y_\eta-\eta Q_{\mathrm L}\|
 \le\eta^3C^4,\qquad \|R_\eta-R\|\le\eta^2C^3.
\end{equation}
These bounds follow directly from the positive resolvent identity.
As in \cref{sp:dyson}, define
\begin{equation}\label{sp:inv-H}
\mathcal H_K=\begin{pmatrix}0&T_K\\T_K^*&0\end{pmatrix},\quad
J_K=\begin{pmatrix}0&K\\K^*&0\end{pmatrix},\quad
M_\eta=(\id\otimes\tau)[(i\eta I-\mathcal H_K)^{-1}].
\end{equation}
By \cref{sp:mde}, including when $K$ is singular,
\begin{equation}\label{sp:general-mean-mde}
M_\eta^{-1}=i\eta I-J_K-\mathcal S(M_\eta).
\end{equation}
Multiplying by $M_\eta$ gives the two block identities
\begin{equation}\label{sp:general-block-identities}
\begin{aligned}
(\eta I+\Phi(X_\eta))Y_\eta
 +(K-\Theta(R_\eta))R_\eta&=I,\\
(\eta I+\Phi(X_\eta))R_\eta^*
 &=(K-\Theta(R_\eta))X_\eta.
\end{aligned}
\end{equation}
The first identity tends to $(K-\Theta(R))R=I$. Thus the numerical
matrix $R$ is invertible and $R^{-1}=K-\Theta(R)$. Divide the second
identity by $\eta$ and pass to the limit. This yields
\[
(I+\Phi(Q_{\mathrm R}))R^*=R^{-1}Q_{\mathrm R},\qquad
Q_{\mathrm R}-\Phi_{R\mathbf A}(Q_{\mathrm R})=RR^*>0.
\]
Choose $0<\delta<1$ such that $RR^*\ge\delta Q_{\mathrm R}$ and define
\begin{equation}\label{sp:inv-T}
\mathcal T(H)=Q_{\mathrm R}^{-1/2}\,
\Phi_{R\mathbf A}(Q_{\mathrm R}^{1/2}HQ_{\mathrm R}^{1/2})\,
Q_{\mathrm R}^{-1/2}.
\end{equation}
The map $\mathcal T$ is completely positive, similar to $\Phi_{R\mathbf A}$,
and satisfies $\mathcal T(I)\le(1-\delta)I$. By
\cref{sp:positive-norm},
\[
r(\Phi_{R\mathbf A})=r(\mathcal T)\le\|\mathcal T\|
=\|\mathcal T(I)\|\le1-\delta<1.
\]
This proves necessity. Conversely, if $R\in\mathrm{GL}_m(\C)$ satisfies
\eqref{sp:invertibility-eq}, \cref{sp:wick-inverse} applies to
$T_K$. It proves that $T_K$ is boundedly invertible and that
$(\id\otimes\tau)(T_K^{-1})=R$. Thus every such $R$ has the asserted
value.
\end{proof}

The limits in \cref{sp:compact} allow both the coefficients and
regularization to vary: take $A_j^{(n)}\to A_j$ and $\eta_n\downarrow0$,
form the normalized quadruples using \eqref{sp:normalize}, and pass
to any convergent subsequence. The next corollary applies to every
such limit, including those with fixed coefficients.

\begin{corollary}\label[corollary]{sp:limit-cases}
For a nonzero irreducible tuple, let $(x,\alpha,\widehat Y,R)$ be any
limit from \cref{sp:compact}. If $\widehat Y=0$, then
\begin{equation}\label{sp:root-normalization}
 R^{-1}=I-\Theta(R),\qquad
 x=\alpha RR^*+\Phi_{R\mathbf A}(x),\qquad r(\Phi_{R\mathbf A})\le1.
\end{equation}
The three possibilities are
\begin{align*}
 \widehat Y>0&\ \Longrightarrow\ T\text{ is not invertible},\\
 \widehat Y=0,\ \alpha>0&\ \Longleftrightarrow\ T\text{ is invertible},\\
 \widehat Y=0,\ \alpha=0&\ \Longrightarrow\ r(\Phi_{R\mathbf A})=1
                               \text{ and }T\text{ is not invertible}.
\end{align*}
\end{corollary}
\begin{proof}
If $\widehat Y=0$, \eqref{sp:normalized-eq} gives
$x^{-1}=B^*C_0^{-1}B$ and $R=B^{-1}$, hence
$x=RC_0R^*=\alpha RR^*+\Phi_{R\mathbf A}(x)$.
Conjugating $\Phi_{R\mathbf A}$ by $H\mapsto x^{1/2}Hx^{1/2}$ gives a
positive map sending $I$ to at most $I$. Its norm is therefore
at most one by \cref{sp:positive-norm}, so $r(\Phi_{R\mathbf A})\le1$.
If $\alpha>0$, the strict inequality $\Phi_{R\mathbf A}(x)<x$ gives
$r(\Phi_{R\mathbf A})<1$ by the same argument. If $\alpha=0$, the
positive fixed matrix $x$ also gives $r(\Phi_{R\mathbf A})\ge1$, hence equality.

Conversely, if the limiting pencil $T$ is invertible, the approximating
pencils $T^{(n)}$ converge to it in norm and have uniformly bounded
inverses for large $n$. Equation \eqref{sp:resolvent-remainders}
therefore applies uniformly, with
\[
 Q_{\mathrm R}^{(n)}
 =(\id\otimes\tau)[((T^{(n)})^*T^{(n)})^{-1}]
 \longrightarrow Q_{\mathrm R}>0.
\]
Thus $\widehat Y_{\eta_n}\to0$ and
$\alpha_{\eta_n}\to1/\Tr Q_{\mathrm R}>0$.
Together with \cref{sp:invertibility}, this proves all three cases.
\end{proof}
\subsection{Uniqueness and regularity as \texorpdfstring{$\eta\downarrow0$}{eta tends to zero}}\label{sp:uniqueness}
We first recall two positive-map facts used to control the limiting
Dyson equation. Their hypotheses are stated separately, so the
spectral-radius bound is also available for reducible tuples.

For a projection $p$,
$\ran\Phi_{\mathbf B}(p)=\Span_j B_j\ran p$.
Thus irreducibility of the tuple is equivalent to irreducibility
of the positive map in the sense of \cite{EHK}.
The adjoint tuple is also irreducible, since orthogonal complements
exchange its invariant subspaces with those of the original tuple.
This verifies the hypotheses for $\Phi_{\mathbf B}$ and $\Phi_{\mathbf B}^*$ below.
\begin{lemma}[Perron--Frobenius theorem; Evans--H{\o}egh-Krohn {\cite[Theorems 2.3--2.4]{EHK}}]
\label[lemma]{sp:perron-eigenmatrices}
If $(B_j)$ is nonzero and irreducible, then $r=r(\Phi_{\mathbf B})>0$ is a
simple eigenvalue and there are $X,Y>0$, unique under the normalizations
$\Tr X=\Tr Y=1$, such that $\Phi_{\mathbf B}(X)=rX$ and $\Phi_{\mathbf B}^*(Y)=rY$.
The eigenvalue $r$ and these eigenmatrices depend real-analytically
on real-analytic coefficients.
\end{lemma}

\begin{proof}
The Perron eigenvalue, its positive eigenmatrices, and geometric
simplicity are supplied by the cited theorem. There can be no nontrivial
Jordan chain at $r$: if $(\Phi_{\mathbf B}-r\id)Z=X$, pairing with $Y$ would
give $0=\Tr(YX)>0$. The analytic-dependence assertion now follows
from the implicit-function theorem for a simple eigenvalue.
\end{proof}

\begin{lemma}\label[lemma]{sp:perron}
For arbitrary tuples $(B_j),(C_j)$ in $M_m(\C)$,
\begin{equation}\label{sp:cross-bound}
 r\left(H\mapsto\sum_jB_jHC_j\right)
       \le\sqrt{r(\Phi_{\mathbf B})r(\Phi_{\mathbf C})}.
\end{equation}
Suppose both tuples are irreducible, $r(\Phi_{\mathbf B}),r(\Phi_{\mathbf C})\le1$,
and $Z=\sum_jB_jZC_j\ne0$. Then both radii equal one and
there is $H\in\mathrm{GL}_m(\C)$ with $B_jH=HC_j^*$ for all $j$.
\end{lemma}
\begin{proof}
\textit{The spectral-radius bound.}
For $\mathcal L(H)=\sum_jB_jHC_j$, expand the iterates as
\[
 \mathcal L^n(H)=\sum_{i_1,\ldots,i_n}
 B_{i_1}\cdots B_{i_n}H C_{i_n}\cdots C_{i_1}.
\]
The row--column Cauchy--Schwarz inequality gives
\[
 \|\mathcal L^n(H)\|
 \le\|\Phi_{\mathbf B}^n(I)\|^{1/2}\|(\Phi_{\mathbf C}^*)^n(I)\|^{1/2}\|H\|.
\]
For the completely positive maps here,
$\|\Psi\|=\|\Psi(I)\|$ by \cref{sp:positive-norm}.
Taking $n$th roots and using the spectral-radius formula proves
\eqref{sp:cross-bound}.

\textit{The equality case.}
Since $\mathcal L(Z)=Z\ne0$, the bound just proved implies
\[
 1\le r(\mathcal L)\le\sqrt{r(\Phi_{\mathbf B})r(\Phi_{\mathbf C})}\le1.
\]
Thus both radii are one. Choose $X,Y>0$
with $\Phi_{\mathbf B}(X)=X$ and $\Phi_{\mathbf C}^*(Y)=Y$, and put
$P=\operatorname{diag}(X,Y)$, $\widetilde{\mathcal A}_j=P^{-1/2}\operatorname{diag}(B_j,C_j^*)P^{1/2}$.
The map $\mathcal F(U)=\sum_j\widetilde{\mathcal A}_jU\widetilde{\mathcal A}_j^*$ is unital and completely
positive, and
\begin{align*}
 N=\begin{pmatrix}0&X^{-1/2}ZY^{-1/2}\\0&0\end{pmatrix}
 \quad\text{satisfies}\quad \mathcal F(N)=N.
\end{align*}
Choose $X_*,Y_*>0$ with $\Phi_{\mathbf B}^*(X_*)=X_*$ and
$\Phi_{\mathbf C}(Y_*)=Y_*$, and define
\[
 W=P^{1/2}\operatorname{diag}(X_*,Y_*)P^{1/2},\qquad
 \varphi(U)=\frac{\Tr(WU)}{\Tr W}.
\]
The fixed-point equations give
\[
 \sum_j\widetilde{\mathcal A}_j\widetilde{\mathcal A}_j^*=I,\qquad
 \mathcal F^*(W)=P^{1/2}
 \operatorname{diag}(\Phi_{\mathbf B}^*(X_*),\Phi_{\mathbf C}(Y_*))P^{1/2}=W.
\]
Since $W>0$, the functional $\varphi$ is faithful and invariant:
$\varphi\circ\mathcal F=\varphi$.
The Kadison--Schwarz inequality \cite[Corollary 2.8]{ChoiSchwarz}
gives $\mathcal F(NN^*)-NN^*\ge0$. Invariance yields
\[
 \varphi(\mathcal F(NN^*)-NN^*)=0
 \quad\Longrightarrow\quad \mathcal F(NN^*)=NN^*.
\]
The implication uses faithfulness on positive matrices. Expanding
the following sum, using $\mathcal F(N)=N$ and
$\sum_j\widetilde{\mathcal A}_j\widetilde{\mathcal A}_j^*=I$, now gives
\begin{align*}
 0=\mathcal F(NN^*)-NN^*
   =\sum_j(\widetilde{\mathcal A}_jN-N\widetilde{\mathcal A}_j)(\widetilde{\mathcal A}_jN-N\widetilde{\mathcal A}_j)^*.
\end{align*}
Thus $\widetilde{\mathcal A}_jN=N\widetilde{\mathcal A}_j$, or $B_jH=HC_j^*$ with $H=ZY^{-1}\ne0$.
Its kernel and range are invariant under the respective irreducible
tuples, so $H$ is invertible.
\end{proof}

For the rest of this subsection the coefficient tuple is nonzero
and irreducible and satisfies the no-intertwiner condition
\eqref{sp:admissible}.

\begin{lemma}\label[lemma]{sp:root-uniqueness}
Let $(A_j)$ be nonzero and irreducible and satisfy
\eqref{sp:admissible}. There is at most one $R\in\mathrm{GL}_m(\C)$ satisfying
\begin{equation}\label{sp:root}
 R^{-1}=I-\Theta(R),\qquad r(\Phi_{R\mathbf A})\le1.
\end{equation}
At such a solution the tuple $(RA_j)$ is irreducible, and the derivative in $R$ of
$R^{-1}-I+\Theta(R)$ is invertible. Hence the solution extends
locally as a holomorphic function of holomorphic coefficient parameters.
\end{lemma}
\begin{proof}
Put $B_j=RA_j$. The equation for $R$ implies
\begin{equation}\label{sp:root-identity}
 R=I+\sum_j B_j^2.
\end{equation}
Every common invariant subspace of the $B_j$ is therefore invariant
under $R$. Since $R$ is invertible and the subspace is finite-dimensional,
it is also invariant under $R^{-1}$, and hence under $A_j=R^{-1}B_j$.
Thus $(B_j)$ is irreducible.
For two solutions $R_1,R_2$, let $Z=(R_1-R_2)R_2^{-1}$. Subtracting
their equations gives
\[
 R_1^{-1}(R_1-R_2)R_2^{-1}=\Theta(R_1-R_2),
 \qquad Z=\sum_j(R_1A_j)Z(R_2A_j).
\]
If $Z\ne0$, \cref{sp:perron} gives an invertible $H$
with $(R_1A_j)H=H(R_2A_j)^*$. Equation \eqref{sp:root-identity}
then gives $R_1H=HR_2^*$. With
$\widetilde H=R_1^{-1}H=H(R_2^*)^{-1}$, we obtain
\[
 A_j\widetilde H
 =R_1^{-1}(R_1A_j)H(R_2^*)^{-1}
 =R_1^{-1}HA_j^*=\widetilde H A_j^*,
\]
contrary to \eqref{sp:admissible}.
For the derivative, a kernel vector $\dot R$ satisfies
\[
 -R^{-1}\dot R R^{-1}+\Theta(\dot R)=0
 \quad\Longrightarrow\quad
 Z=\sum_j(RA_j)Z(RA_j),\qquad Z=\dot R R^{-1}.
\]
The same argument with $R_1=R_2=R$ gives $Z=0$, hence $\dot R=0$.
The finite-dimensional derivative is therefore invertible.
\end{proof}

\begin{lemma}\label[lemma]{sp:positive-stability}
Let $(A_j)$ be nonzero and irreducible and satisfy
\eqref{sp:admissible}. Suppose $M=\left(\begin{smallmatrix}-iY&-R^*\\-R&-iX\end{smallmatrix}\right)$,
with $X,Y>0$, solves $M^{-1}=-J_0-\mathcal S(M)$.
For $Q=\operatorname{diag}(I,-I)$ and $\Pi_0=\operatorname{diag}(Y,X)$,
\begin{align*}
 \ker\bigl(H\mapsto H-M\mathcal S(H)M\bigr)=\C Q\Pi_0.
\end{align*}
In particular this kernel has zero intersection with
$\{H:\Tr H_{11}=\Tr H_{22}\}$.
\end{lemma}
\begin{proof}
\textit{Step 1: identification of the kernel
with a fixed-point space.}
Use $\mathcal S(H)=\sum_j\mathcal A_jH\mathcal A_j$ with the matrices $\mathcal A_j$ defined
in \eqref{sp:calA-def}. The identities
$QM=-M^*Q$ and $Q\mathcal A_j=-\mathcal A_jQ$ imply
\[
 Q\bigl(M\mathcal S(H)M\bigr)=M^*\mathcal S(QH)M.
\]
Consequently $H\mapsto QH$ carries the desired kernel onto the
fixed space of the completely positive map
$\mathcal F(Z)=M^*\mathcal S(Z)M$.
Taking imaginary parts of $M^{-1}=-J_0-\mathcal S(M)$ gives
\[
 \Im M^{-1}=M^{-*}\Pi_0M^{-1}=\mathcal S(\Pi_0),
 \qquad \mathcal F(\Pi_0)=\Pi_0.
\]

\textit{Step 2: verify irreducibility.}
The completely positive map $\mathcal F$ is implemented by the matrices $M^*\mathcal A_j$. A tuple is
irreducible if and only if its adjoint tuple is irreducible, since
orthogonal complements exchange invariant subspaces. We therefore
consider a common invariant subspace $F$ for the matrices $\mathcal A_jM$.
The Dyson equation implies
\[
 J_0M=-I-\sum_j(\mathcal A_jM)^2,
 \qquad J_0MF\subset F.
\]
The matrix $J_0M$ is invertible, so its restriction maps $F$ onto $F$.
It follows that
\[
 \mathcal A_jJ_0F=(\mathcal A_jM)(J_0M)^{-1}F\subset F,
 \qquad \mathcal A_jJ_0=\operatorname{diag}(A_j,A_j^*).
\]
Let $\pi_1,\pi_2$ be the coordinate projections on $\C^m\oplus\C^m$.
The spaces $\pi_1F,\pi_2F$ and
$F\cap(\C^m\oplus0),F\cap(0\oplus\C^m)$ are invariant for the
appropriate diagonal tuples. By irreducibility, each is either zero
or a full coordinate space. If $F$ is not zero, a coordinate space,
or all of $\C^{2m}$, both projections are onto and both intersections
are zero. Thus $F$ is the graph
\[
 F=\{(Hv,v):v\in\C^m\},\qquad H\in\mathrm{GL}_m(\C).
\]
Invariance under $\operatorname{diag}(A_j,A_j^*)$ then gives
$A_jH=HA_j^*$, contradicting \eqref{sp:admissible}.
The two coordinate spaces are also impossible, since
\[
 J_0M=\begin{pmatrix}-R&-iX\\-iY&-R^*\end{pmatrix}
\]
has invertible off-diagonal blocks. Hence the tuple $(M^*\mathcal A_j)_j$ is irreducible.

\textit{Step 3: remove the remaining scalar direction.}
Since $\Pi_0>0$ and $\mathcal F(\Pi_0)=\Pi_0$, conjugation by
$\Pi_0^{1/2}$ turns $\mathcal F$ into a unital positive map.
Thus $r(\mathcal F)=1$, and the Perron--Frobenius theorem \cref{sp:perron-eigenmatrices} gives
$\ker(\id-\mathcal F)=\C \Pi_0$. Step 1 yields the stated kernel
$\C Q\Pi_0$. For $H=\lambda Q\Pi_0$, the trace constraint is
\[
 \Tr H_{11}-\Tr H_{22}
 =\lambda(\Tr Y+\Tr X)=0,
\]
which forces $\lambda=0$.
\end{proof}

The following result generalizes \cite[Theorem 3.17 and Corollary 3.18]{ZhongElliptic} to the matrix-valued setting.
\begin{proposition}\label[proposition]{sp:regularity}
Let $A_j\colon \mathcal O\to M_m(\C)$ be real-analytic on an open set
$\mathcal O\subset\R^p$, and suppose the assumptions of
this subsection hold at every $u\in \mathcal O$. Define $\mathcal O_+$ to be the set
of $u$ for which some subsequential limit of $\widehat Y_\eta(u)$
is positive definite.
\begin{enumerate}[label=(\roman*),leftmargin=*]
\item $\mathcal O_+$ is open. On $\mathcal O_+$, $(X_\eta,Y_\eta,R_\eta)$ converges
locally uniformly to a unique triple $(X,Y,R)$ with $X,Y>0$.
These limits and $V$ are real-analytic.
\item If $u_0\notin \mathcal O_+$, then $\widehat Y_\eta(u)\to0$ and
$R_\eta(u)\to R(u_0)$ jointly as $(u,\eta)\to(u_0,0)$, where
$R(u_0)$ is the solution of \eqref{sp:root}. Together with (i), this
makes $R$ continuous on $\mathcal O$ and $R_\eta\to R$ locally uniformly.
\item $V\in C^1(\mathcal O)$, and $V_\eta\to V$ locally in $C^1$. For
each coordinate $u_a$,
\begin{equation}\label{sp:limit-gradient}
 2m\partial_{u_a}V=
 \begin{cases}
 \Tr\bigl(Y(\partial_{u_a}\Phi_u)(X)\bigr)
 -\Re\Tr\bigl(R(\partial_{u_a}\Theta_u)(R)\bigr),&u\in \mathcal O_+,\\
 -\Re\Tr\bigl(R(\partial_{u_a}\Theta_u)(R)\bigr),&u\notin \mathcal O_+.
 \end{cases}
\end{equation}
\end{enumerate}
\end{proposition}

\begin{proof}
Write
$\mathcal Z(u,\eta)=(x_\eta(u),\alpha_\eta(u),\widehat Y_\eta(u),R_\eta(u))$.
For $u_0\in \mathcal O$, its set of joint limits is
\begin{equation}\label{sp:joint-limits}
 \mathscr C(u_0)=\bigcap_{\epsilon>0}
 \overline{\{\mathcal Z(u,\eta):|u-u_0|<\epsilon,\ 0<\eta<\epsilon\}},
\end{equation}
where it suffices to take small $\epsilon$ with the coefficient ball
inside $\mathcal O$. The estimates in \cref{sp:compact}, applied to any
sequence $(u,\eta)\to(u_0,0)$, give local boundedness and show that
$\mathscr C(u_0)$ is nonempty and compact. It is connected since each set inside the closure is a continuous image of a connected set, and these closures are nested compact connected sets for small $\epsilon$.

\textit{(i) Positive limits, openness, and analyticity.}
Suppose a joint limit has $\widehat Y>0$. By \cref{sp:compact} it
determines $X,Y>0$ with $\Tr X=\Tr Y$ and a solution $M$ of the
Dyson equation at $\eta=0$. At $\eta=0$ the system is invariant under
$(X,Y,R)\mapsto(\lambda X,\lambda^{-1}Y,R)$.
The trace balance $\Tr X=\Tr Y$ removes the tangent to this scaling
orbit, which is the kernel direction in \cref{sp:positive-stability}.
Accordingly, the variations lie in the real vector space
\[
 \mathcal V=\left\{
 \begin{pmatrix}-i\dot Y&-\dot R^*\\-\dot R&-i\dot X\end{pmatrix}:
 \dot X=\dot X^*,\ \dot Y=\dot Y^*,\ \Tr\dot X=\Tr\dot Y
 \right\}.
\]
The real dimension of $\mathcal V$ is $4m^2-1$.
For $\mathcal E(M)=M^{-1}+J_0+\mathcal S(M)$ its derivative is
\[
 D\mathcal E_M[H]=-M^{-1}HM^{-1}+\mathcal S(H).
\]
The structured target $\{N:QN=-N^*Q\}$, with
$Q=\operatorname{diag}(I,-I)$, has real dimension $4m^2$.
Multiplication by $M$ on both sides identifies the derivative's
kernel on $\mathcal V$ with
\[
 \mathcal V\cap\ker(H\mapsto H-M\mathcal S(H)M)=\{0\},
\]
by \cref{sp:positive-stability}. Thus we can select $4m^2-1$
real components of the full equation
$\mathcal E(M)-i\eta I_{2m}=0$ whose derivative on $\mathcal V$
is invertible at this solution. The real-analytic implicit-function
theorem gives a unique local branch solving these selected equations,
with parameters $(u,\eta)$. At this stage the branch is known to solve
only the selected equations. In particular, at $(u_0,0)$ the given
solution is isolated among trace-balanced solutions of the full system.

The normalization can be inverted near this positive solution:
\[
 \gamma=\sqrt{\Tr\widehat Y},\qquad X=\gamma x,\quad
 Y=\widehat Y/\gamma;\qquad
 x=X/\Tr X,\quad\widehat Y=(\Tr X)Y.
\]
Every sufficiently nearby point of $\mathscr C(u_0)$ has $\widehat Y>0$.
By \cref{sp:compact}, $\alpha=0$ and the displayed formulas reconstruct
a trace-balanced solution of the full equation near $M$. This solution
also satisfies the selected equations, so it equals $M$ by local
uniqueness. Hence the nearby point coincides with the given one, which
is therefore isolated in the connected set $\mathscr C(u_0)$ and is its
only point.
Compactness now gives joint convergence of $\mathcal Z(u,\eta)$ to
this point as $(u,\eta)\to(u_0,0)$. Otherwise a sequence staying away
would have a different joint limit. In particular, the same limit is
obtained with $u=u_0$ fixed.

By joint convergence and the normalization formulas, the actual
matrices $M_\eta(u)$ lie in the uniqueness neighborhood of the
implicit-function theorem for all $(u,\eta)$ sufficiently close
to $(u_0,0)$ with $\eta>0$. They solve the full Dyson system,
and hence the selected equations, so they coincide with the analytic
branch for $\eta>0$. By continuity the branch also solves the full
equation at $\eta=0$. Its $X,Y$ remain positive for $u$ near $u_0$.
Thus $u_0$ and a neighborhood of it belong to $\mathcal O_+$, proving
openness. The analytic extension gives locally uniform
convergence and real-analytic limits there. Formula \eqref{sp:action-eq}
then gives real analyticity of $V$.

\textit{(ii) Limits outside $\mathcal O_+$.}
If $u_0\notin \mathcal O_+$, part (i) excludes every positive joint limit,
not only limits with $u=u_0$ fixed. By \cref{sp:compact},
\[
 \widehat Y_\eta(u)\longrightarrow0\quad
 \text{as }(u,\eta)\to(u_0,0).
\]
By \cref{sp:limit-cases}, each joint limit of $R_\eta(u)$ satisfies
\eqref{sp:root}. Its uniqueness follows from \cref{sp:root-uniqueness}.
Hence $R_\eta(u)\to R(u_0)$ jointly, which proves (ii) and continuity
of $R$ on $\mathcal O\setminus \mathcal O_+$.

\textit{(iii) Convergence of the gradients.}
The first term in \eqref{sp:gradient} can be written without the
possibly vanishing scalar normalization:
\[
 \Tr\bigl(Y_\eta(\partial_{u_a}\Phi_u)(X_\eta)\bigr)
 =\Tr\bigl(\widehat Y_\eta(\partial_{u_a}\Phi_u)(x_\eta)\bigr).
\]
At a point of $\mathcal O_+$, part (i) determines its joint limit. At a point
outside $\mathcal O_+$, part (ii) makes it tend to zero, since $\|x_\eta\|\le1$.
The term involving $R_\eta$ also has a unique joint limit in both
cases. Thus $\nabla V_\eta$ has a unique joint limit $g(u_0)$ at
every $u_0$. The sequential criterion for continuity and compactness
give continuity of $g$ and locally uniform convergence to $g$.

To identify $g$ as the gradient of $V$, take a convex ball compactly
contained in $\mathcal O$ and a base point $u_*$. On that ball,
\[
 V_\eta(u)-V_\eta(u_*)
 =\int_0^1\nabla V_\eta(u_*+t(u-u_*))\cdot(u-u_*)\dd t.
\]
The convergence of the gradients is uniform, and
$V_\eta(u_*)\downarrow V(u_*)> -\infty$. Passing to the limit proves
locally uniform convergence of the potentials and $V\in C^1$ with
$\nabla V=g$. The two expressions for $g$ are exactly
\eqref{sp:limit-gradient}.
\end{proof}

\subsection{Harmonicity and constancy of traces}\label{sp:rigidity}
We now let the coefficients depend holomorphically on $z$.
In this subsection a prime denotes $\partial_z$.
The set $\mathcal O_+$ and the limits below are those of \cref{sp:regularity},
with $u=(\Re z,\Im z)$. By \cref{sp:limit-cases}, a noninvertible
point either lies in $\mathcal O_+$ or satisfies $r(\Phi_{R\mathbf A})=1$.
We treat these two cases separately, in
\cref{sp:positive-rigidity} and \cref{sp:critical-rigidity}.

\begin{proposition}\label[proposition]{sp:positive-rigidity}
Let $(A_j(z))_{j=1}^k$ be a holomorphic nonzero irreducible family
on a connected domain $D\subset\C$, satisfying \eqref{sp:admissible}.
If $D\subset \mathcal O_+$ and $V$ is harmonic on $D$, then at every point
there is a matrix $H$ such that $A_j'=[H,A_j]$ for all $j$.
In particular, $\Tr A_w$ is constant on $D$ for every word $w$.
\end{proposition}
\begin{proof}
\textit{Step 1: vary all coefficients independently.}
Allow variations of the constant term as well. On
$\mathrm{GL}_m(\C)\times M_m(\C)^k$, define the real-valued potential
\begin{equation}\label{sp:F-def}
 F(K,\mathbf A)=\log\Delta_m\left(K\otimes1+\sum_jA_j\otimes s_j\right)
       =m^{-1}\log|\det K|+V(K^{-1}\mathbf A).
\end{equation}
Here $K^{-1}\mathbf A=(K^{-1}A_1,\ldots,K^{-1}A_k)$; finiteness follows
from \cref{sp:bounds}. We work near each point $(I,\mathbf A(z))$.
Apply \cref{sp:regularity} to the normalized tuple $K^{-1}\mathbf A$.
Irreducibility and \eqref{sp:admissible} are open, as noted after \cref{sp:no-intertwiner},
and positivity persists by openness of $\mathcal O_+$.
Thus $F$ is real-analytic near each $(I,\mathbf A(z))$.
It is plurisubharmonic: on every complex line in coefficient space,
\cref{lem:levi} applies to the regularized logarithmic determinant,
and the decreasing limit is finite by \cref{sp:bounds}.

Using the complex partial derivatives defined in
\cref{subsec:shared-rotation}, write $\partial_{A_j}F(K,\mathbf A)[L]$
for the $\partial_\zeta$ derivative at $\zeta=0$ when $A_j$ is
replaced by $A_j+\zeta L$, and define $\partial_KF(K,\mathbf A)[L]$
similarly. Then \eqref{sp:limit-gradient} gives
\begin{equation}\label{sp:complex-gradient}
 2m\partial_{A_j}F(I,\mathbf A)[L]=\Tr\bigl((XA_j^*Y-RA_jR)L\bigr),
 \qquad 2m\partial_KF(I,\mathbf A)[L]=\Tr(RL).
\end{equation}
For the second identity, differentiate $K^{-1}A_j$ at $K=I$ in
direction $L$, obtaining $-LA_j$. Consequently
\begin{align*}
 2m\partial_KF(I,\mathbf A)[L]
 &=\Tr L-\sum_j\Tr\bigl((XA_j^*Y-RA_jR)LA_j\bigr)\\
 &=\Tr\bigl((I-\Phi(X)Y+\Theta(R)R)L\bigr)=\Tr(RL).
\end{align*}
The last equality is \eqref{sp:block-trace} at $K=I$, $\eta=0$.
There $C=\Phi(X)$; in particular the identities at $\eta=0$
used below are
\begin{equation}\label{sp:zero-height-blocks}
 (I-\Theta(R))R+\Phi(X)Y=I,\qquad
 XB^*=RC,\qquad BX=CR^*,\quad B=I-\Theta(R).
\end{equation}

\textit{Step 2: use positivity of the complex Hessian.}
For a real-valued $C^2$ function $f$ of complex coordinates
$\zeta_1,\ldots,\zeta_N$, its \emph{Levi matrix} is the complex
Hessian $(\partial_{\zeta_a}\partial_{\bar\zeta_b}f)_{a,b}$;
the associated Hermitian quadratic form is called the Levi form.
Plurisubharmonicity means that this matrix is positive semidefinite,
as follows by restricting $f$ to each complex line.

Consider $\mathscr F(z,L,E)=F(I+L,\mathbf A(z)+E)$, where the entries of
$L,E_1,\ldots,E_k$ are independent complex coordinates. This is
plurisubharmonic, being the composition of $F$ with a holomorphic map,
so its Levi matrix is positive semidefinite. At $L=E=0$, harmonicity
gives $\partial_z\partial_{\bar z}\mathscr F=0$.
For each scalar coefficient coordinate $\xi$, nonnegativity of the
corresponding two-by-two principal minor gives
\[
 |\partial_\xi\partial_{\bar z}\mathscr F|^2
 \le (\partial_\xi\partial_{\bar\xi}\mathscr F)
       (\partial_z\partial_{\bar z}\mathscr F)=0.
\]
Hence every coefficient derivative of $\mathscr F$ is holomorphic
along $z\mapsto(I,\mathbf A(z))$. Applying this first to $K$ and then to
$A_j$ in \eqref{sp:complex-gradient} gives
\[
 \partial_{\bar z}R=0,\qquad
 \partial_{\bar z}(XA_j^*Y-RA_jR)=0,
 \qquad \partial_{\bar z}(XA_j^*Y)=0.
\]
The last implication uses holomorphicity of $R$ and $A_j$.
Thus $D_j:=XA_j^*Y$ is holomorphic for every $j$.

\textit{Step 3: differentiate the Dyson identities.}
Put $B=I-\Theta(R)$, $C=\Phi(X)$, and
\begin{equation}\label{sp:HX-def}
 H_X=(\partial_z X)X^{-1},\qquad H_Y=Y^{-1}(\partial_z Y).
\end{equation}
Since $D_j=XA_j^*Y$ is holomorphic, $\partial_z(D_j^*)=\partial_z(YA_jX)=0$.
Multiplication on the left by $Y^{-1}$ and on the right by $X^{-1}$
gives
\[
 A_j'=-H_YA_j-A_jH_X.
\]
Consequently, using $\partial_z(A_j^*)=0$ and $\partial_z X=H_XX$,
\[
 \partial_z C=\sum_j(A_j'XA_j^*+A_j(\partial_z X)A_j^*)=-H_YC.
\]
The nonzero space $\sum_j\ran A_j$ is invariant under every $A_j$,
so irreducibility makes it all of $\C^m$. Since $X>0$,
$\ker C=\bigcap_j\ker A_j^*=\{0\}$, and hence $C>0$.
The identities \eqref{sp:zero-height-blocks} at $\eta=0$, valid
by \cref{sp:regularity}(i), give $XB^*=RC$ and $BX=CR^*$.
Both $B$ and $R$ are holomorphic. Differentiating the first identity
and canceling $C$ therefore gives
\[
 H_XRC=R'C-RH_YC,\qquad R'=H_XR+RH_Y.
\]
Substituting the formulas for $A_j'$ and $R'$ into the product rule,
\[
 (\Theta(R))'=\sum_j\bigl(A_j'RA_j+A_jR'A_j+A_jRA_j'\bigr)
 =-H_Y\Theta(R)-\Theta(R)H_X.
\]
We can now compute $B'$ in two ways. Differentiating
$B=I-\Theta(R)$ gives the first line below; differentiating
$BX=CR^*$ and using $\partial_z(R^*)=0$ gives the second:
\[
 \begin{aligned}
 B'&=H_Y+H_X-H_YB-BH_X,\\
 B'&=-H_YB-BH_X.
 \end{aligned}
\]
Thus $H_Y=-H_X$, and $A_j'=[H_X,A_j]$ for every $j$.
The Leibniz rule then gives
\[
 (A_w)'=[H_X,A_w],\qquad (\Tr A_w)'=\Tr[H_X,A_w]=0.
\]
Connectedness of $D$ proves the constancy assertion.
\end{proof}

\begin{lemma}\label[lemma]{sp:curvature}
Let $(B_j(z))_{j=1}^k$ be a holomorphic nonzero irreducible tuple,
and set $r(z)=r(H\mapsto\sum_jB_j(z)HB_j(z)^*)$.
Then $r$ is subharmonic.
If $r$ is harmonic on an open set, at each
point of that set there is a matrix $H$ such that
$B_j'=[H,B_j]$ for all $j$, where $[H,B_j]=HB_j-B_jH$.
\end{lemma}
\begin{proof}
\textit{Normalize the eigenmatrix at one point.}
By \cref{sp:perron-eigenmatrices}, choose a local real-analytic
$X(z)>0$ satisfying $\Phi_{\mathbf B}(X)=rX$. Fix $z_0$ and first apply a
constant similarity so that $X(z_0)=I$. Set
\[
 S(z)=I+(z-z_0)(\partial_z X)(z_0),\qquad
 \widetilde B_j=S^{-1}B_jS,\qquad
 \widetilde X=S^{-1}XS^{-*}.
\]
The matrix $S$ is invertible near $z_0$. These substitutions preserve
the eigenmatrix equation and the spectral radius. At $z_0$,
\[
 \widetilde X=I,\qquad
 \partial_z\widetilde X=\partial_z X-S'=0,\qquad
 \partial_{\bar z}\widetilde X=\partial_{\bar z}X-(S')^*=0.
\]
The last equality uses $X=X^*$.

\textit{Differentiate and pair with a left eigenmatrix.}
In the rest of this calculation all quantities are evaluated at $z_0$.
Differentiate
$\sum_j\widetilde B_j\widetilde X\widetilde B_j^*=r\widetilde X$
in $z$ and $\bar z$. Holomorphicity of the $\widetilde B_j$ and the
vanishing first derivatives of $\widetilde X$ give
\begin{align*}
 \sum_j\widetilde B_j'\widetilde B_j^*&=(\partial_z r)I,\\
 \Phi_{\widetilde{\mathbf B}}(\partial_z\partial_{\bar z}\widetilde X)
       +\sum_j\widetilde B_j'\widetilde B_j'^*
 &=r(\partial_z\partial_{\bar z}\widetilde X)+(\partial_z\partial_{\bar z}r)I.
\end{align*}
Choose $Y>0$ with $\Phi_{\widetilde{\mathbf B}}^*(Y)=rY$.
Pairing the second line with $Y$ cancels the terms in
$\partial_z\partial_{\bar z}\widetilde X$:
\[
 (\partial_z\partial_{\bar z}r)\Tr Y
 =\sum_j\Tr(Y\widetilde B_j'\widetilde B_j'^*)
 =\sum_j\|Y^{1/2}\widetilde B_j'\|_{\HS}^2\ge0.
\]
Since $z_0$ was arbitrary, $r$ is subharmonic. If $r$ is harmonic,
each $\widetilde B_j'$ vanishes. With $H=S'(z_0)$,
\[
 0=\widetilde B_j'(z_0)=B_j'(z_0)-[H,B_j(z_0)].
\]
Undoing the initial constant similarity gives the claimed common
commutator in the original coordinates.
\end{proof}

\begin{lemma}\label[lemma]{sp:holomorphic-root}
Let $(A_j(z))$ be holomorphic, nonzero, and irreducible on an open
set $D$, satisfying \eqref{sp:admissible}. If $D\cap \mathcal O_+=\varnothing$,
then the limit $R$ in \cref{sp:regularity}(ii) is holomorphic on $D$.
In this case $V$ is harmonic, by \eqref{sp:limit-gradient}.
\end{lemma}
\begin{proof}
Fix $z_0\in D$. By \cref{sp:regularity}(ii), $R$ is continuous,
and by \cref{sp:limit-cases} it solves
$R^{-1}-I+\Theta_{\mathbf A(z)}(R)=0$ with $r(\Phi_{R\mathbf A})\le1$.
The derivative of this equation in $R$ is invertible at
$(z_0,R(z_0))$ by \cref{sp:root-uniqueness}.
The holomorphic implicit-function theorem therefore supplies a
neighborhood $N$ of $z_0$, a neighborhood $\mathcal U$ of $R(z_0)$,
and a holomorphic branch $\widetilde R$ that is the unique solution
in $\mathcal U$ for each $z\in N$. Shrink $N$ so that continuity
puts $R(z)$ in $\mathcal U$. Both $R(z)$ and $\widetilde R(z)$ solve
the same equation there, so they agree. Since $z_0$ was arbitrary,
$R$ is holomorphic on $D$.
\end{proof}

\begin{proposition}\label[proposition]{sp:critical-rigidity}
Let $(A_j(z))_{j=1}^k$ be holomorphic, nonzero, and irreducible on
a connected domain $D\subset\C$, satisfying \eqref{sp:admissible}.
If $D\cap \mathcal O_+=\varnothing$ and $r(\Phi_{R(z)\mathbf A(z)})=1$ for every $z\in D$,
then $\Tr A_w$ is constant on $D$ for every word $w$.
\end{proposition}
\begin{proof}
By \cref{sp:holomorphic-root}, $R$ is holomorphic on $D$.
Set $B_j=RA_j$.
This tuple is irreducible by \cref{sp:root-uniqueness}, and
$r(\Phi_{\mathbf B})=1$. By \cref{sp:curvature}, at each point there is a
common $H$ with $B_j'=[H,B_j]$. Differentiating
$R=I+\sum_jB_j^2$ gives
\[
 R'=\sum_j\bigl([H,B_j]B_j+B_j[H,B_j]\bigr)
    =[H,\sum_jB_j^2]=[H,R].
\]
Since $A_j=R^{-1}B_j$, it follows that
\[
 A_j'=-R^{-1}[H,R]R^{-1}B_j+R^{-1}[H,B_j]=[H,A_j].
\]
The Leibniz rule now yields $A_w'=[H,A_w]$ for every word, and
$ (\Tr A_w)'=\Tr[H,A_w]=0$. No regularity of the choice of $H$
is needed: the calculation is pointwise.
\end{proof}

\begin{proof}[Proof of \cref{sp:harmonic-disk}]
Suppose $T$ is noninvertible at $z_0\in D$. Apply
\cref{sp:regularity} on the open subset of $\Omega$ where
\eqref{sp:admissible} holds; this subset contains $D$.

If $D\cap \mathcal O_+\ne\varnothing$, its openness gives a disk
$D_0\subset D\cap \mathcal O_+$. By \cref{sp:positive-rigidity}, every
$\Tr A_w$ is constant on $D_0$. Holomorphic continuation makes
these traces constant on the connected domain $\Omega$.
By \cref{sp:traces}, the coefficient tuples at any two points of
$\Omega$ are simultaneously similar, and so are their pencils.
Since $T(z_0)$ is noninvertible, $T$ is noninvertible throughout
$\Omega$.

Otherwise $D\cap \mathcal O_+=\varnothing$. Then
\cref{sp:holomorphic-root} makes $R$ holomorphic on $D$.
The tuple $(RA_j)_j$ is irreducible by \eqref{sp:root-identity},
and \cref{sp:limit-cases} gives $r(\Phi_{R\mathbf A})\le1$ throughout $D$.
At $z_0$, noninvertibility rules out the case $\alpha>0$, so
that same corollary gives $r(\Phi_{R(z_0)\mathbf A(z_0)})=1$.
By \cref{sp:curvature} this radius is subharmonic; the maximum
principle on the connected set $D$ therefore makes it identically
one. Now \cref{sp:critical-rigidity} gives constant word traces
on $D$, and holomorphic continuation and \cref{sp:traces} give
the same conclusion throughout $\Omega$ as in the first case.

Finally, if $T(z_*)$ is invertible for some $z_*\in\Omega$,
noninvertibility at any point of $D$ would contradict this conclusion.
Thus $T$ is invertible throughout $D$.
\end{proof}

The next lemma allows us to choose a disk on which
\eqref{sp:admissible} holds for every block family whose coefficients
cannot all be conjugated to a fixed tuple by one common change of basis
at each parameter. There are only finitely
many blocks, and each such family excludes a closed set with
empty interior.
\begin{lemma}\label[lemma]{sp:exceptional}
Let $(A_j(z))_{j=1}^k$ be holomorphic and irreducible on a connected
domain $\Omega\subset\C$. The set
$\{z\in\Omega:\exists H\ne0,\ A_j(z)H=HA_j(z)^*\text{ for all }j\}$
is closed in $\Omega$. This set has empty interior unless, for every
$z\in\Omega$, there is an invertible $S_z$ with
$A_j(z)=S_zA_j(z_0)S_z^{-1}$ for all $j$, where $z_0\in\Omega$ is fixed.
\end{lemma}
\begin{proof}
Closedness follows from the rank condition for
$H\mapsto(A_j(z)H-HA_j(z)^*)_j$.
Irreducibility makes every nonzero solution $H$ invertible.
At such a parameter,
$\Tr A_w=\overline{\Tr A_{w^{\mathrm{rev}}}}$, where
$w^{\mathrm{rev}}=(j_\ell,\ldots,j_1)$ for $w=(j_1,\ldots,j_\ell)$.
If the set contains a disk, every trace of a word is both holomorphic
and antiholomorphic on that disk, hence constant.
The hypotheses of \cref{sp:traces} now hold: the tuple is
holomorphic and irreducible on the connected domain $\Omega$,
and all word traces are constant on a nonempty open subset.
\Cref{sp:traces} supplies the simultaneous similarities $S_z$ to the
tuple at $z_0$, proving the assertion.
\end{proof}

\section{Support of the Brown measure and spectral convergence}\label{sp:support}
We first prove equality of the spectrum and the support of the Brown
measure for polynomials in elliptic elements. We then adapt the
complex-rotation argument of \cref{sec:weak} to prove \cref{sp:main}, and
establish Hausdorff convergence of the spectra of GUE polynomials.
Throughout this section, the ambient finite von Neumann algebra is
equipped with a faithful normal tracial state $\tau$. Recall that 
$U_B(z)=\log\Delta_\tau(B-z)$ and the normalized trace is denoted by
$\tau_m=m^{-1}(\Tr\otimes\tau)$ in the definition of $\Delta_m$.

\subsection{Elliptic polynomials: the nondegenerate case}\label{sp:elliptic}
Let $s_1,\ldots,s_k,t_1,\ldots,t_k$ be a jointly free standard
semicircular family. For $0\le\rho\le1$, put
\begin{equation}\label{sp:elliptic-input}
 a_\rho=\sqrt{(1+\rho)/2},\qquad
 b_\rho=\sqrt{(1-\rho)/2},\qquad
 e_j=a_\rho s_j+i b_\rho t_j.
\end{equation}
The rotation dictionary is
\begin{equation}\label{sp:rotation-dictionary}
 \rho_y=\frac1{\cosh(2y)},\qquad
 \cosh(y)s_j+i\sinh(y)t_j=\sqrt{\cosh(2y)}\,e_j(\rho_y).
\end{equation}
Then $\tau(e_i e_j^*)=\delta_{ij}$ and
$\tau(e_i e_j)=\rho\delta_{ij}$. In particular, $\rho=0$ gives
a circular family and $\rho=1$ gives the original semicircular family.
For a matrix tuple $\mathbf A=(A_j)$, write
\begin{equation}\label{sp:PhiA-def}
 \Phi_{\mathbf A}(H)=\sum_j A_jHA_j^*,\qquad
 \Theta_{\mathbf A}(H)=\sum_j A_jHA_j.
\end{equation}

\begin{lemma}\label[lemma]{sp:elliptic-reduction}
For $0\le\rho\le1$ and $A_1,\ldots,A_k\in M_m(\C)$, set
\begin{equation}\label{sp:tilde-A}
 \widetilde{\mathbf A}=(a_\rho A_1,\ldots,a_\rho A_k,
                    i b_\rho A_1,\ldots,i b_\rho A_k).
\end{equation}
The pencil $I+\sum_j A_j\otimes e_j$ is the semicircular pencil
with coefficients $\widetilde{\mathbf A}$, and
\begin{equation}\label{sp:elliptic-covariance}
 \Phi_{\widetilde{\mathbf A}}=\Phi_{\mathbf A},\qquad
 \Theta_{\widetilde{\mathbf A}}=\rho\Theta_{\mathbf A}.
\end{equation}
The tuples $\mathbf A$ and $\widetilde{\mathbf A}$ have the same common invariant
subspaces. If $\mathbf A$ is nonzero and irreducible and $\rho<1$, then
$\widetilde{\mathbf A}$ satisfies the no-intertwiner condition
\eqref{sp:admissible}:
\[
 \{H:\widetilde A_jH=H\widetilde A_j^*\text{ for all }j\}
 =\{0\}.
\]
\end{lemma}
\begin{proof}
The identities follow from $a_\rho^2+b_\rho^2=1$ and
$a_\rho^2-b_\rho^2=\rho$. Since $a_\rho>0$, the invariant
subspaces agree. When $\rho<1$, an intertwiner $H$ satisfies both
$A_jH=HA_j^*$ and $A_jH=-HA_j^*$ for every $j$. Hence
$\ran H\subseteq\bigcap_j\ker A_j$. This intersection is a
common invariant subspace. Since $\mathbf A$ is irreducible and nonzero,
it must be $\{0\}$.
\end{proof}

\begin{corollary}\label[corollary]{sp:elliptic-invertibility}
Let $e_1,\ldots,e_k$ be the free elliptic family of
\eqref{sp:elliptic-input}, with $0\le\rho\le1$. For arbitrary
$K\in M_m(\C)$ and $\mathbf B=(B_1,\ldots,B_k)\in M_m(\C)^k$, the pencil
$L=K\otimes1+\sum_j B_j\otimes e_j$ is invertible if and only if
there exists $R\in\mathrm{GL}_m(\C)$ such that
\begin{equation}\label{sp:elliptic-invertibility-eq}
 R^{-1}=K-\rho\Theta_{\mathbf B}(R),\qquad
 r(\Phi_{R\mathbf B})<1,\qquad
 \Phi_{R\mathbf B}(H)=R\Phi_{\mathbf B}(H)R^*.
\end{equation}
Every such $R$ equals $(\id\otimes\tau)(L^{-1})$.
In particular, for the linearization $L_z$ of a polynomial $P$ in
\cref{lem:linearization}, $z\notin\spec P(e)$ if and only if
\eqref{sp:elliptic-invertibility-eq} holds with $K=K_z$ for some
$R\in\mathrm{GL}_m(\C)$.
\end{corollary}
\begin{proof}
By \cref{sp:elliptic-reduction}, $L$ is a semicircular pencil with
constant term $K$ and doubled coefficient tuple $\widetilde{\mathbf B}$.
Its covariance maps are $\Phi_{\mathbf B}$ and $\rho\Theta_{\mathbf B}$
by \eqref{sp:elliptic-covariance}. Applying \cref{sp:invertibility}
gives the criterion and the formula for $R$. For $L_z$, the Schur
factorization \eqref{eq:schur-factorization} gives invertibility of
$L_z(e)$ exactly when $P(e)-z$ is invertible.
\end{proof}

We will also use the following elementary observation.
\begin{lemma}\label[lemma]{sp:riesz}
Every isolated point of the spectrum of an element $B$ in a finite
von Neumann algebra with faithful trace belongs to $\supp\mu_B$.
Consequently, if $\spec B\setminus\supp\mu_B$ is finite, it is empty.
\end{lemma}
\begin{proof}
Suppose that $z_0$ is an isolated spectral point outside
$\supp\mu_B$. Choose a disk $D$ about $z_0$ whose closure is
disjoint from $\supp\mu_B$ and contains no other spectral point.
The potential $U_B$ is harmonic near $\overline D$, and on
$\partial D$, \eqref{eq:brown-resolvent} gives
$2\partial_zU_B(z)=\tau((z-B)^{-1})$. Thus the nonzero Riesz
idempotent
\[
 q=\frac{1}{2\pi i}\int_{\partial D}(z-B)^{-1}\dd z
\]
has trace zero. If $p$ is its range projection, then $pq=q$ and
$qp=p$, so $\tau(q)=\tau(p)>0$, a contradiction.
If $\spec B\setminus\supp\mu_B$ is finite, each of its points
is isolated in $\spec B$, since $\supp\mu_B$ is closed.
The first assertion therefore applies.
\end{proof}

\begin{lemma}\label[lemma]{sp:block-disk}
Let $P\in\C\langle x_1,\ldots,x_k\rangle$ be a nonconstant
scalar polynomial without adjoints, fix $0\le\rho\le1$, and put
$p_0=P(0)$. Then there exists a finite set $\Lambda\subset\C$, containing
$p_0$, such that on every simply connected domain
$\Omega\subset\C\setminus \Lambda$ there are holomorphic irreducible
coefficient tuples $\mathbf A^{(\ell)}=(A_j^{(\ell)})_{j=1}^k$ with
\begin{equation}\label{sp:scalar-block-potential}
 U_{P(e)}(z)=\log|p_0-z|+\sum_{\ell=1}^q m_\ell V_\ell(z),
 \qquad V_\ell(z)=\log\Delta_{m_\ell}(T_\ell(z)),
\end{equation}
where
\[
 T_\ell(z)=I_{m_\ell}+\sum_j A_j^{(\ell)}(z)\otimes e_j.
\]
Each tuple is either nonzero throughout $\Omega$ or identically
zero with $m_\ell=1$. Moreover, $P(e)-z$ is invertible if and
only if every $T_\ell(z)$ is invertible.
For every $0\le\rho\le1$, the functions $V_\ell$ are finite and
subharmonic. Wherever $U_{P(e)}$ is harmonic, all the $V_\ell$
are harmonic as well.
\end{lemma}
\begin{proof}
The decomposition follows from \cref{sp:linear-det,sp:blocks}.
By \cref{sp:elliptic-reduction,sp:bounds}, each $V_\ell$ is finite,
and it is subharmonic by \cref{lem:levi}.
Where $U_{P(e)}$ is harmonic, \eqref{sp:scalar-block-potential}
shows that
\[
 \sum_{\ell=1}^q m_\ell\Lap V_\ell=0
\]
in the sense of distributions. Since each $\Lap V_\ell$ is a
nonnegative measure and $m_\ell>0$, every $V_\ell$ is harmonic there.
\end{proof}

\begin{theorem}
\label[theorem]{sp:elliptic-theorem}
For $0\le\rho<1$ and every scalar noncommutative polynomial $P$
without adjoints,
\[
 \spec P(e)=\supp\mu_{P(e)}.
\]
\end{theorem}
\begin{proof}
The assertion is immediate for constant $P$, so assume $P$ is
nonconstant. Let $\Lambda$ be as in \cref{sp:block-disk}.
Choose a disk $D$ with
\[
 \overline D\subset\C\setminus(\supp\mu_{P(e)}\cup \Lambda).
\]
Choose $z_*\notin \Lambda\cup\spec P(e)$, which is possible for
$|z_*|$ sufficiently large. Choose a simply connected domain \(\Omega\subset\mathbb C\setminus\Lambda\) containing \(\overline D\) and \(z_*\) as in \cref{sp:joining-domain}. By \cref{sp:block-disk},
all the block pencils $T_\ell(z_*)$ are invertible, and their
potentials $V_\ell$ are harmonic on $D$.

For each nonzero block, \cref{sp:elliptic-reduction} gives an
irreducible semicircular coefficient tuple satisfying the
no-intertwiner condition \eqref{sp:admissible}. By \cref{sp:harmonic-disk}, \(T_\ell(z)\) is invertible for every \(z\in D\). If all coefficients of a block vanish, then \(T_\ell(z)=1\). Thus $P(e)-z$ is invertible
on $D$, and
$\spec P(e)\setminus\supp\mu_{P(e)}\subseteq \Lambda$.
Apply \cref{sp:riesz} to finish the proof.
\end{proof}

\subsection{The semicircular endpoint}
\label{sp:rotation}\label{sp:endpoint}
We first show that the boundary of $\spec P(s)$ belongs to the
support of its Brown measure. The argument uses the following
consequence of \cref{prop:subharm}.

\begin{lemma}\label[lemma]{sp:moment-roots}
Let $B(w)$ be a norm-holomorphic family in a finite tracial von
Neumann algebra. For every integer $\ell\ge1$,
\[
 M_\ell(w)=\left(\int_\C |z|^{2\ell}\dd\mu_{B(w)}(z)\right)^{1/(2\ell)}
\]
is subharmonic, and $\log M_\ell$ is also subharmonic or identically
$-\infty$ on each connected component.
Furthermore,
\[
 M_\ell(w)\uparrow r_{\mathrm B}(B(w))
 :=\max\{|z|:z\in\supp\mu_{B(w)}\}.
\]
\end{lemma}
\begin{proof}
By \cref{sp:det-facts},
$u(w):=M_\ell(w)^{2\ell}=b(B(w)^\ell)$ is subharmonic, where
$b(C)=\int|z|^2\dd\mu_C(z)$.
For every holomorphic scalar function $f$,
\[
 w\longmapsto e^{\Re f(w)}u(w)
   =b(e^{f(w)/2}B(w)^\ell)
\]
is also subharmonic by \cref{prop:subharm}.
Let $h$ be harmonic near a closed disk, and suppose
$\log u\le h$ on its boundary. Locally $h$ is the real part
of a holomorphic function, so $e^{-h}u$ is subharmonic.
By the maximum principle, $e^{-h}u\le1$ throughout the disk.
It follows that $\log u$ is subharmonic (or identically $-\infty$).
Since exponentiation preserves subharmonicity,
$M_\ell=\exp((2\ell)^{-1}\log u)$ is subharmonic wherever $u$ is not
identically zero; on a component where $u\equiv0$, $M_\ell\equiv0$
is subharmonic as well. Here and below $\log0=-\infty$.
Finally, the convergence of $L^p$ norms to the $L^\infty$ norm
gives the last assertion.
\end{proof}

\begin{proposition}
\label[proposition]{sp:boundary}
For every $z_0\notin\spec P(s)$,
\begin{equation}\label{sp:distance-transfer}
 \operatorname{dist}(z_0,\spec P(s))
 =\operatorname{dist}(z_0,\supp\mu_{P(s)}).
\end{equation}
In particular,
$\partial\spec P(s)\subseteq\supp\mu_{P(s)}$.
\end{proposition}
\begin{proof}
Take a free standard semicircular family $(s,t)$ and put
\begin{equation}\label{sp:X-rotation}
 X(w)=P(\cos(w)s+\sin(w)t).
\end{equation}
By \eqref{eq:rotation-law}, for real $x$ the pair
$(\cos(x)s+\sin(x)t,-\sin(x)s+\cos(x)t)$ has the same joint
$*$-distribution as $(s,t)$. Hence the $*$-distribution of
$X(x+iy)$ depends only on $y$, as do its spectrum and
resolvent norms, by faithfulness of the trace.

Fix $z_0\notin\spec P(s)$. For all real $x$,
\[
 \|(X(x)-z_0)^{-1}\|=\|(P(s)-z_0)^{-1}\|.
\]
Moreover, periodicity gives
\[
 \sup_{x\in\R}\|X(x+iy)-X(x)\|\longrightarrow0
 \qquad(y\to0).
\]
By the Neumann series, there is $a>0$ such that
$\mathcal R(w)=(X(w)-z_0)^{-1}$ is norm-holomorphic and uniformly
bounded on the strip $|\Im w|<a$.

Apply \cref{sp:moment-roots} to $\mathcal R(w)$. Each $M_\ell(x+iy)$
depends only on $y$ and is subharmonic, so
$y\mapsto M_\ell(iy)$ is convex by \cref{lem:rotation-strip}. Consequently,
\begin{equation}\label{sp:Fy-def}
 F(y):=r_{\mathrm B}(\mathcal R(iy))=\sup_\ell M_\ell(iy)
\end{equation}
is convex. The spectral radius of a norm-holomorphic family is
also subharmonic: this follows, for example, by exponentiating
the decreasing limit of the subharmonic functions
$2^{-j}\log\|\mathcal R(w)^{2^j}\|$; see also \cite{Vesentini1968}.
Thus
\begin{equation}\label{sp:Gy-def}
 G(y):=r_{\mathrm{op}}(\mathcal R(iy)),\qquad
 r_{\mathrm{op}}(B)=\max\{|z|:z\in\spec B\},
\end{equation}
is convex as well. Both $F$ and $G$ are finite on $(-a,a)$,
and hence continuous at zero.

For $0<y<a$, the family
\[
 \frac{\cosh(y)s+i\sinh(y)t}{\sqrt{\cosh(2y)}}
\]
is elliptic with parameter $\rho=1/\cosh(2y)<1$ by
\eqref{sp:rotation-dictionary}.
Applying \cref{sp:elliptic-theorem} to
$Q_y(x)=P(\sqrt{\cosh(2y)}\,x)$ gives
$\spec X(iy)=\supp\mu_{X(iy)}$.
The spectral mapping theorem and the corresponding property
of the Brown measure (\cref{sp:det-facts}) then give
$F(y)=G(y)$. By continuity, $F(0)=G(0)$, while
\[
 F(0)=\frac1{\operatorname{dist}(z_0,\supp\mu_{P(s)})},
 \qquad
 G(0)=\frac1{\operatorname{dist}(z_0,\spec P(s))}.
\]
This proves \eqref{sp:distance-transfer}.
If $z\in\partial\spec P(s)$, take $z_n\notin\spec P(s)$
with $z_n\to z$. Then
\[
 \operatorname{dist}(z_n,\supp\mu_{P(s)})
 =\operatorname{dist}(z_n,\spec P(s))
 \le|z_n-z|\longrightarrow0,
\]
so $z\in\supp\mu_{P(s)}$ by closedness.
\end{proof}

\begin{proof}[Proof of \cref{sp:main}]
The assertion is immediate for constant $P$, so assume $P$ is
nonconstant. Suppose
\[
 \spec P(s)\setminus\supp\mu_{P(s)}\ne\varnothing.
\]
Since $\partial\spec P(s)\subseteq\supp\mu_{P(s)}$
by \cref{sp:boundary}, this set is open. Use
\cref{sp:block-disk} with $\rho=1$. Choose a disk
$D$ with
\[
 \overline D\subset
 \spec P(s)\setminus\supp\mu_{P(s)}
\]
avoiding the finite exceptional set $\Lambda$ of
\cref{sp:block-disk}. Choose $z_*\notin \Lambda\cup\spec P(s)$ and
choose a joining domain $\Omega\subset\C\setminus \Lambda$ as in
\cref{sp:joining-domain}.
By \cref{sp:block-disk}, all block potentials are harmonic
on $D$, and all block pencils are invertible at $z_*$.

Consider a nonzero irreducible block tuple $\mathbf A^{(\ell)}$.
If it is simultaneously similar to a fixed tuple
throughout $\Omega$, its pencil is invertible throughout
$\Omega$, since it is invertible at $z_*$.
Otherwise, \cref{sp:exceptional} shows that
\begin{equation}\label{sp:El-def}
 E_\ell=\left\{z\in\Omega:
 \begin{array}{l}
 \exists H\ne0 \text{ such that}\\
 A_j^{(\ell)}(z)H=H A_j^{(\ell)}(z)^*
 \quad\text{for all }j
 \end{array}
 \right\}
\end{equation}
is closed in $\Omega$ with empty interior.
Since there are finitely many blocks, we can choose
a nonempty disk $\overline{D'}\subset D$ such that
\[
 D'\cap\bigcup_\ell E_\ell=\varnothing,
\]
where the union is over the blocks that are not
simultaneously conjugate to a fixed tuple.
On $D'$, these blocks have harmonic potentials
and satisfy \eqref{sp:admissible}, so their pencils
are invertible by \cref{sp:harmonic-disk}.
The remaining blocks are also invertible, as
observed above. Hence $P(s)-z$ is invertible
on $D'$, contradicting $D'\subset\spec P(s)$.
\end{proof}

\begin{example}[An element of $C^*(s_1,s_2)$ with a smaller support]\label[example]{sp:nonpolynomial-example}
The conclusion of \cref{sp:main} fails for general elements of
$C^*(s_1,s_2)$. Let $F$ be the semicircle distribution function and set
$u=\exp(2\pi iF(s_1))$, $h=|s_2|^{1/4}$, and $T=uh$.
Then $u$ is a Haar unitary free from $h$. The formulas \cite[Theorem 4.4(i) and Proposition 4.6]{HL} give
\[
 \supp\mu_T=\{z:r_-\le|z|\le r_+\},\qquad
 \spec T=\{z:|z|\le r_+\},
\]
where $r_-=(\tau(|s_2|^{-1/2}))^{-1/2}>0$ and
$r_+=(\tau(|s_2|^{1/2}))^{1/2}$. The negative moment is finite
because the semicircle density is bounded near zero, whereas $h$
is not invertible. Thus $0\in\spec T\setminus\supp\mu_T$.
\end{example}

\subsection{Hausdorff convergence}
\label{sp:matrix}
Let $G_1^{(n)},\ldots,G_k^{(n)}$ be independent normalized GUE
matrices, and put $A_n=P(G^{(n)})$, $A=P(s)$.
Strong convergence \cite[Theorem A]{HT},
\cref{thm:main}, and the equality
$\supp\mu_A=\spec A$ imply
almost-sure Hausdorff convergence by
\cite[Proposition 7.3(3)]{BCKP}.
We prove the quantitative estimates below.

\paragraph{Strong convergence and the absence of outliers.}
\begin{theorem}[Haagerup--Thorbj{\o}rnsen {\cite[Theorem A]{HT}}]\label[theorem]{sp:HT}
For a fixed independent normalized GUE family and free standard
semicircular family, every fixed scalar $*$-polynomial $Q$ satisfies
\[
 \|Q(G^{(n)})\|\longrightarrow\|Q(s)\|
 \quad\text{almost surely}.
\]
\end{theorem}
We use \cref{sp:HT} to prove that $\sg(A_n-zI_n)$ converges almost
surely to $\sg(A-z)$, locally uniformly in $z$.
For the spectral inclusion, see also \cite[Proposition 7.3(1)]{BCKP}.
For related upper bounds on the spectral radius, see
\cite[Lemma 2.12]{vH2026} and \cite[Remark 3.5]{Miyagawa2023}.

For a compact set $K\subset\C$, put
$K_\epsilon=\{z:\operatorname{dist}(z,K)<\epsilon\}$.
\begin{lemma}\label[lemma]{sp:exclusion}
Let $A_n=P(G^{(n)})$ and $A=P(s)$. Almost surely, for every
compact $K\subset\C$,
\begin{equation}\label{sp:uniform-lsv}
 \sup_{z\in K}|\sg(A_n-zI_n)-g(z)|\longrightarrow0,
 \qquad g(z)=\sg(A-z).
\end{equation}
Here $g$ is $1$-Lipschitz, vanishes precisely on $\spec A$, and
equals $\|(A-z)^{-1}\|^{-1}$ off that spectrum. Consequently,
for every $\epsilon>0$, almost surely
$\spec A_n\subset(\spec A)_\epsilon$ for all sufficiently large $n$.
\end{lemma}
\begin{proof}
For fixed $z$, put $Q_z=(P-z)^*(P-z)$ and choose
$M>\|Q_z(s)\|$. By \cref{sp:HT}, eventually
$0\le Q_z(G^{(n)})<MI_n$ and
\[
 \lambda_{\min}(Q_z(G^{(n)}))
 =M-\|M-Q_z(G^{(n)})\|
 \longrightarrow M-\|M-Q_z(s)\|=g(z)^2.
\]
This argument also applies when $z\in\spec A$. Take square roots
and intersect the probability-one events for a countable dense set
of $z$. The $1$-Lipschitz dependence on $z$ on both sides gives
uniform convergence on every compact set by finite nets. Strong
convergence also bounds $\|A_n\|$ eventually, almost surely.
Inside a fixed containing disk, the compact set at distance at
least $\epsilon$ from $\spec A$ has a positive minimum of $g$.
Equation \eqref{sp:uniform-lsv} excludes eigenvalues there.
\end{proof}

\begin{proof}[Proof of \cref{thm:hausdorff}]
\textit{Step 1: no outliers.}
We first obtain a uniform resolvent bound on a compact set
$K\subset\C\setminus\spec A$, assumed nonempty. Set
\begin{equation}\label{sp:gK-def}
 g_K=\min_{z\in K}\sg(A-z)>0,\qquad
 F_{n,K}=\min_{z\in K}\sg(A_n-zI_n).
\end{equation}
By \cref{sp:exclusion}, $F_{n,K}\to g_K$ in probability; only
this consequence of strong convergence is needed here.
Let $\mathcal N_n=\{\max_j\|G_j^{(n)}\|\le8\}$.
Then $\Pp(\mathcal N_n^c)\le Ce^{-cn}$ by \cref{lem:norm}.
By \eqref{eq:singular-calculus}, taking a minimum over $K$
preserves the operator-norm Lipschitz bound. Thus
\cref{lem:gaussian-lipschitz} gives Gaussian-coordinate Lipschitz
constant $L_P/\sqrt n$ on $\mathcal N_n$, uniformly in $K$.

Extend $F_{n,K}$ from $\mathcal N_n$ to all Gaussian coordinates
with the same Lipschitz constant, and truncate the extension to
$[0,B_K]$, where $B_K$ is a fixed bound for $F_{n,K}$ on
$\mathcal N_n$. Denote the resulting function by $f_{n,K}$.
Since it agrees with $F_{n,K}$ on $\mathcal N_n$ and is uniformly
bounded, $f_{n,K}\to g_K$ in probability and in mean.
For large $n$, $\E f_{n,K}\ge3g_K/4$. Gaussian concentration
\cite[Corollary 4.4.19]{AGZ} therefore gives
\begin{equation}\label{sp:exp-lsv}
 \begin{aligned}
 \Pp\{F_{n,K}\le g_K/2\}
 &\le\Pp(\mathcal N_n^c)
   +\Pp\{f_{n,K}-\E f_{n,K}\le-g_K/4\}\\
 &\le Ce^{-cn}.
 \end{aligned}
\end{equation}
In particular, for constants depending on $P,K$,
\begin{equation}\label{sp:exp-resolvent}
 \Pp\left\{\sup_{z\in K}\|(A_n-zI_n)^{-1}\|>2/g_K\right\}
 \le Ce^{-cn}
\end{equation}
for all sufficiently large $n$. The inverse norm is understood
to be $+\infty$ at a singular matrix.

Choose $R_P>\|A\|$ so that $\|A_n\|\le R_P$ on
$\mathcal N_n$. For $\epsilon>0$, put
\begin{equation}\label{sp:F-eps}
 F_\epsilon=\{z:|z|\le R_P,\
                   \operatorname{dist}(z,\spec A)\ge\epsilon\}.
\end{equation}
If $F_\epsilon\ne\varnothing$, an eigenvalue of $A_n$ in
$F_\epsilon$ forces $F_{n,F_\epsilon}=0$.
Using \eqref{sp:exp-lsv} with $K=F_\epsilon$ and the bound on
$\Pp(\mathcal N_n^c)$, we obtain
\begin{equation}\label{sp:outer-estimate}
 \Pp\left\{\sup_{\lambda\in\spec A_n}
       \operatorname{dist}(\lambda,\spec A)\ge\epsilon\right\}
 \le C_\epsilon e^{-c_\epsilon n}.
\end{equation}
If $F_\epsilon$ is empty, the norm bound alone gives the same
conclusion. This proves \eqref{eq:main-outer}.
It also follows that, for $h\in C_c(\C\setminus\supp\mu_A)$,
\begin{equation}\label{sp:exp-test-function}
 \E\left|\int h\dd\mu_{A_n}\right|
 \le\|h\|_\infty\Pp\{\spec A_n\cap\supp h\ne\varnothing\}
 \le C_he^{-c_hn},
\end{equation}
since \cref{sp:main} identifies $\supp\mu_A$ with $\spec A$.

\textit{Step 2: eigenvalues near the limiting support.}
Choose finitely many
$z_j\in\supp\mu_A$ such that the balls
$B(z_j,\epsilon/2)$ cover $\supp\mu_A$.
For each $j$, choose a nonnegative smooth function $h_j$
supported in $B(z_j,\epsilon/2)$ and positive near $z_j$.
Then $a_j:=\int h_j\dd\mu_A>0$. By \eqref{eq:green},
\[
 \left|\int h_j\dd\mu_{A_n}-a_j\right|
 \le\frac{\|\Lap h_j\|_\infty}{2\pi}
      \|U_{A_n}-U_{A}\|_{L^1(\supp h_j)}.
\]
If $\int h_j\dd\mu_{A_n}=0$, the left side equals $a_j$.
Consequently \cref{as:integrated} gives
\[
 \Pp\left\{\int h_j\dd\mu_{A_n}=0\right\}
 =O_{P,\epsilon,j,q}(n^{-q})\qquad(q>0).
\]
The same argument with a smooth function supported in any open ball
$B$ with $\mu_A(B)>0$ gives
$\Pp\{\spec A_n\cap B=\varnothing\}=O_{P,B,q}(n^{-q})$.
If all the integrals $\int h_j\dd\mu_{A_n}$ are positive, every
point of $\supp\mu_A$ is within $\epsilon$ of an eigenvalue of
$A_n$. Thus, for every fixed $\epsilon,q>0$,
\begin{equation}\label{sp:inner-estimate}
 \begin{aligned}
 \Pp\left\{\sup_{z\in\supp\mu_A}\operatorname{dist}(z,\spec A_n)
                  \ge\epsilon\right\}
 &\le\sum_j\Pp\left\{\int h_j\dd\mu_{A_n}=0\right\}\\
 &=O_{P,\epsilon,q}(n^{-q}).
 \end{aligned}
\end{equation}
\textit{Step 3: Hausdorff convergence.}
Together with \eqref{sp:outer-estimate} and \(\supp\mu_A=\spec A\),
this proves \eqref{eq:main-hausdorff-tail}. Taking $q=2$ and applying
the Borel--Cantelli lemma for $\epsilon=1/m$, $m=1,2,\ldots$, gives the
almost-sure limit in \eqref{eq:main-hausdorff}. Finally,
$|r(A_n)-r(A)|\le d_H(\spec A_n,\spec A)$ proves spectral-radius convergence.
\end{proof}

\section{An algorithm for computing Brown measures}
\label{alg:algorithm}
Let $P\in\C\langle x_1,\ldots,x_k\rangle$ and let
$s_1,\ldots,s_k$ be free standard semicircular elements in
$(\mathcal M,\tau)$. This section turns the determinant formula of
\cref{sp:action} into a procedure that computes
\[
 U_P(z)=\log\Delta_\tau(z1-P(s)),\qquad
 \mu_{P(s)}=\frac1{2\pi}\Lap U_P.
\]
\Cref{alg:determinant} expresses the regularized potential of an
affine pencil for $P$ through the solution of a matrix Dyson equation.
\Cref{alg:solver} gives an iteration that computes this solution, and
two linear systems for the Laplacian of the regularized potential.
\Cref{alg:procedure} states the algorithm and proves that, as the
regularization is removed, it recovers $U_P$ and $\mu_{P(s)}$, even
when that measure has no planar density.
\Cref{alg:adjoint-example} gives the calculation for the polynomial in
\cref{fig:adjoint-example}.

The algorithm specializes the Hermitization and linearization
framework of Belinschi--\'Sniady--Speicher \cite[Sections 3--4]{BSS}
and Helton--Mai--Speicher
\cite[Sections 6.4--6.6 and Algorithm 6.15]{HMS}, see also \cite{BMS},
to semicircular inputs. Its determinant
formulation evaluates the normalized potential directly from the
Dyson solution; the matrix unknown itself remains an operator-valued
Cauchy transform. The determinant identity \eqref{sp:kargin-identity}
also appears in \cite[Theorem 3.1]{Kargin}, and the iteration is
supplied by \cite[Theorem 2.1]{HRS}.
For a variational characterization of the potential itself, see
\cite[Corollary 1.5]{Kargin}; here we use regularized determinants
and distributional Laplacians.

\subsection{Affine linearization and the determinant formula}
\label{alg:determinant}
An admissible pencil for $P$ is an affine pencil
\begin{equation}\label{alg:direct-pencil}
 L_z(s)=K_z\otimes1+\sum_j B_j\otimes s_j,\qquad
 K_z=K_0-zE_{11},\quad K_0,B_j\in M_m(\C),
\end{equation}
with constant $B_j$ and the exact determinant identity
\begin{equation}\label{alg:det-normalization}
 U_P(z)=m\log\Delta_{\tau_m}(L_z(s))\qquad(z\in\C).
\end{equation}
The pencil of \cref{lem:linearization} is admissible,
with $K_0=\operatorname{diag}(p_0,I_{m-1})$, $p_0=P(0)$, by
\eqref{eq:schur-factorization} and \cref{sp:linear-det}.
A smaller factored pencil is equally suitable when its Schur
factorization proves \eqref{alg:det-normalization}; one must also
account for the determinant of its auxiliary block. We work directly
with $K_z$, including where it is singular. All matrix traces
$\Tr$ below are unnormalized.

\begin{proposition}\label[proposition]{alg:potential}
For \eqref{alg:direct-pencil}, use $K=K_z$ and $A_j=B_j$ in
\cref{sp:mde,sp:action}. In particular, $M=M_\eta(z)$ is the
unique solution with $\Im M<0$ of
\begin{equation}\label{alg:dyson}
 M^{-1}=i\eta I-J_{K_z}-\mathcal S(M),\qquad
 J_{K_z}=\begin{pmatrix}0&K_z\\K_z^*&0\end{pmatrix}.
\end{equation}
Its blocks satisfy \eqref{sp:schur}, with $B=K_z-\Theta(R)$, and
\begin{equation}\label{alg:action}
 U_\eta(z):=mV_\eta(L_z(s))
 =-\tfrac12\log|\det M|
       +\tfrac14\Re\Tr_{2m}(M\mathcal S(M)).
\end{equation}
\end{proposition}
\begin{proof}
This follows from \cref{sp:mde,sp:action},
with the factor $m$ required by
\eqref{alg:det-normalization}.
\end{proof}

Equivalently, in block form,
\begin{equation}\label{alg:potential-limit}
 U_\eta(z)=mV_\eta(L_z(s))
 =\tfrac12\left[\log\frac{\det C_\eta}{\det X_\eta}
       +\Re\Tr(K_zR_\eta)+\eta\Tr X_\eta-m\right].
\end{equation}
This is the pencil regularization; it need not equal
$\frac12\tau\log((P(s)-z)^*(P(s)-z)+\eta^2)$ at the same $\eta$.
Its zero-regularization limit is nevertheless exactly $U_P$, as
proved in \cref{alg:recovery}.

\subsection{A convergent iteration and differentiation}
\label{alg:solver}

For fixed $\eta>0$, the sign condition $\Im M<0$ selects the
solution corresponding to the operator-valued Cauchy transform.
The following iteration implements this choice. Its convergence
is an application of Helton--Rashidi Far--Speicher
\cite[Theorem 2.1]{HRS}. Belinschi, Mai, and Speicher give a related
fixed-point iteration for operator-valued free additive subordination
\cite[Theorem 2.2]{BMS}.

\begin{proposition}
\label[proposition]{alg:iteration}
For any $K,A_1,\ldots,A_k\in M_m(\C)$ and $\eta>0$, with
$J_K$ and $\mathcal S$ as in \cref{sp:dyson}, set
$W^{(0)}=I_{2m}$ and iterate
\begin{equation}\label{alg:iteration-eq}
 W^{(\ell+1)}=
 \bigl(\eta I_{2m}+iJ_K+\mathcal S(W^{(\ell)})\bigr)^{-1}.
\end{equation}
Every iterate has positive definite real part, and
$W^{(\ell)}\to iM_\eta$ in matrix norm as $\ell\to\infty$.
At this solution the linear map
\begin{equation}\label{alg:linear-map}
 \mathcal L_\eta(H)=M_\eta\mathcal S(H)M_\eta
\end{equation}
has spectral radius strictly less than one.
Consequently $I-\mathcal L_\eta$ is invertible. More precisely,
\begin{equation}\label{alg:iteration-rate}
 \limsup_{\ell\to\infty}\|W^{(\ell)}-iM_\eta\|^{1/\ell}
 \le r(\mathcal L_\eta)<1.
\end{equation}
Newton's method for the Dyson equation converges quadratically
from a sufficiently close initial matrix.
\end{proposition}
\begin{proof}
With $W=iM$, \eqref{alg:dyson} becomes
\[
 W^{-1}=\eta I+iJ_K+\mathcal S(W).
\]
Apply \cite[Theorem 2.1]{HRS} with constant term
$\eta I/2+iJ_K$ and analytic map
$W\mapsto\eta I/2+\mathcal S(W)$. The constant term has strictly
positive real part. The map is bounded on bounded sets and sends
$\Re W>0$ to matrices with real part at least $\eta I/2$,
since $\mathcal S$ is completely positive. The theorem therefore
gives the asserted convergence and uniqueness. The inverse of a
matrix with positive definite real part again has positive definite
real part, so every iterate has the asserted sign.

For linear stability, write $M=M_\eta$, $Q=-\Im M>0$, and
$\mathcal A_j=\left(\begin{smallmatrix}0&A_j\\A_j^*&0\end{smallmatrix}\right)$,
so $\mathcal S(H)=\sum_j\mathcal A_jH\mathcal A_j$.
Taking imaginary parts of the Dyson equation and multiplying on
the two sides in either order gives
\[
 Q=\eta MM^*+M\mathcal S(Q)M^*,\qquad
 Q=\eta M^*M+M^*\mathcal S(Q)M.
\]
Thus the positive maps
$H\mapsto M\mathcal S(H)M^*$ and
$H\mapsto M^*\mathcal S(H)M$ send $Q$ to a matrix strictly
smaller than $Q$. The conjugation by $Q^{1/2}$ turns each into a
positive map sending $I$ to a matrix strictly smaller than $I$;
therefore, its spectral radius is less than one.
Apply \eqref{sp:cross-bound} to the tuples
$B_j=M\mathcal A_j$ and $C_j=\mathcal A_jM$.
The two positive maps just considered are $\Phi_{\mathbf B}$ and $\Phi_{\mathbf C}^*$,
and $r(\Phi_{\mathbf C}^*)=r(\Phi_{\mathbf C})$. Hence
\[
 r(\mathcal L_\eta)
 \le\sqrt{r(\Phi_{\mathbf B})r(\Phi_{\mathbf C})}<1.
\]
This proves the invertibility of $I-\mathcal L_\eta$.
The derivative of the fixed-point map at $iM_\eta$ is
$H\mapsto M_\eta\mathcal S(H)M_\eta$. For each
$r(\mathcal L_\eta)<r_0<1$, an equivalent matrix norm makes this
derivative smaller than $r_0$; the eventual local contraction gives
\eqref{alg:iteration-rate}. The Dyson residual is analytic near
$M_\eta$, with invertible derivative, so the usual local Newton
estimate gives quadratic convergence.
\end{proof}

\begin{proposition}\label[proposition]{alg:derivatives}
For an admissible pencil, write $z=x+iy$, $E=E_{11}$, and
\begin{equation}\label{alg:Jxy}
 J_x=\begin{pmatrix}0&-E\\-E&0\end{pmatrix},\qquad
 J_y=\begin{pmatrix}0&-iE\\iE&0\end{pmatrix}.
\end{equation}
At each $\eta>0$, the two unique solutions $N_x,N_y$ of
\begin{equation}\label{alg:derivative-system}
 (I-\mathcal L_\eta)N_a=M_\eta J_aM_\eta,\qquad a=x,y,
\end{equation}
are $N_a=\partial_aM_\eta$. Set $R_a=-(N_a)_{21}$. Then
\begin{equation}\label{alg:potential-gradient}
 \partial_xU_\eta=-\Re\Tr(ER_\eta),\qquad
 \partial_yU_\eta=\Im\Tr(ER_\eta),
\end{equation}
and
\begin{equation}\label{alg:density-formula}
 \rho_\eta(z)=\frac1{2\pi}\Lap U_\eta(z)
 =\frac{-\Re\Tr(E R_x)+\Im\Tr(E R_y)}{2\pi}.
\end{equation}
\end{proposition}
\begin{proof}
The regularized resolvent is smooth in $x,y$ and $\mathcal S$ is
constant in $z$. Differentiate \eqref{alg:dyson} and multiply by
$M_\eta$ on both sides to obtain \eqref{alg:derivative-system}.
The uniqueness follows from \cref{alg:iteration}. Formula
\eqref{sp:gradient}, with $\partial_xK_z=-E$ and
$\partial_yK_z=-iE$, gives \eqref{alg:potential-gradient}.
One more differentiation proves \eqref{alg:density-formula}.
\end{proof}
Thus, the regularized density requires two linear solves after the
Dyson solution has been computed.

\subsection{The algorithm and its convergence}
\label{alg:procedure}
The input is an admissible pencil \eqref{alg:direct-pencil} for $P$,
for example, the linearization of \cref{lem:linearization}; its size $m$
depends only on $P$.

\begin{center}
\begin{minipage}{0.96\textwidth}
\hrule\vspace{4pt}
\textbf{Algorithm 1.} Regularized potential and density at a point.
\vspace{4pt}\hrule\vspace{5pt}
\textbf{Input:} $K_0,B_1,\ldots,B_k\in M_m(\C)$; $z\in\C$; $\eta>0$;
a tolerance $\delta>0$.\\
\textbf{Output:} approximations $\widehat U_\eta(z)$ and
$\widehat\rho_\eta(z)$ to $U_\eta(z)$ and $\rho_\eta(z)$.
\begin{enumerate}[leftmargin=2em,itemsep=1pt,topsep=3pt]
\item Set $K=K_0-zE_{11}$,
$J=\begin{psmallmatrix}0&K\\K^*&0\end{psmallmatrix}$,
$\mathcal A_j=\begin{psmallmatrix}0&B_j\\B_j^*&0\end{psmallmatrix}$
and $\mathcal S(H)=\sum_j\mathcal A_jH\mathcal A_j$.
\item Set $W=I_{2m}$ and repeat
$W\leftarrow(\eta I_{2m}+iJ+\mathcal S(W))^{-1}$
until $W$ changes by at most $\delta$ in norm. Set $M=-iW$.
\item Set $\widehat U_\eta(z)=-\tfrac12\log|\det M|
+\tfrac14\Re\Tr_{2m}(M\mathcal S(M))$.
\item For $a=x,y$, solve the linear system
$N_a-M\mathcal S(N_a)M=MJ_aM$, with $J_x,J_y$ as in \eqref{alg:Jxy},
and let $R_a$ be minus the lower left $m\times m$ block of $N_a$.
If a system is singular, replace $\delta$ by $\delta/2$ and resume step 2.
\item Set $\widehat\rho_\eta(z)=\bigl(-\Re\Tr(ER_x)+\Im\Tr(ER_y)\bigr)/(2\pi)$,
where $E=E_{11}$.
\end{enumerate}
\vspace{2pt}\hrule
\end{minipage}
\end{center}

Steps 2--5 implement \eqref{alg:iteration-eq}, \eqref{alg:action},
\eqref{alg:derivative-system} and \eqref{alg:density-formula}.
The tolerance $\delta$ controls the change between successive
iterates. For each fixed $\eta>0$, the linear systems of step 4 are
invertible once the iterate is sufficiently close to the solution of
the Dyson equation, by \cref{alg:iteration}. The accuracy of the
outputs is then governed by $\delta$ and the norm of
$(I-\mathcal L_\eta)^{-1}$.

The convergence statement below is for the exact solutions at fixed
$\eta$: one first lets the iteration converge, and then lets
$\eta\downarrow0$. In this order, the potentials converge locally
in $L^1$, and the integrals of $\rho_\eta$ against compactly supported
continuous test functions converge to their integrals against
$\mu_{P(s)}$, as established in the following.

\begin{remark}\label[remark]{alg:newton-remark}
In step 2 one may instead apply Newton's method to the Dyson equation
\eqref{alg:dyson}, starting from the solution at a nearby value of $z$
or at a larger value of $\eta$. It converges quadratically from a
sufficiently close initial matrix (\cref{alg:iteration}), and the
condition $\Im M<0$ identifies the required solution.
\Cref{fig:quadratic-density} uses this method at positive $\eta$.
For \cref{fig:adjoint-example}, the calculation continues the interior
solution in $z$ directly at $\eta=0$ as justified below.
\end{remark}

\begin{proposition}
\label[proposition]{alg:recovery}
For every $z$, $U_\eta(z)\downarrow U_P(z)$ as $\eta\downarrow0$,
and this convergence holds in $L^1_{\mathrm{loc}}(\C)$.
For $\eta>0$, $U_\eta$ is smooth and subharmonic.
Consequently
\begin{equation}\label{alg:density}
 \rho_\eta=\frac1{2\pi}\Lap U_\eta
          \longrightarrow\mu_{P(s)}
\end{equation}
in distributions, and vaguely as positive measures
$\rho_\eta(z)\dd^2z$.
\end{proposition}
\begin{proof}
Monotone spectral calculus and \eqref{alg:det-normalization} give the pointwise
limit. For subharmonicity, put $T=L_z(s)$ and $Q_\eta=T^*T+\eta^2I$.
Tracial differentiation gives
\[
 \Lap U_\eta
 =2(\Tr\otimes\tau)\bigl[Q_\eta^{-1}(T')^*
       (I-TQ_\eta^{-1}T^*)T'\bigr]\ge0,
 \qquad T'=\partial_zT=-E_{11}\otimes1,
\]
since $I-TQ_\eta^{-1}T^*=\eta^2(TT^*+\eta^2I)^{-1}>0$.
The same regularized functional calculus gives smoothness.
The decreasing limit is the locally integrable logarithmic potential of
\cref{thm:brown}, therefore, it is not identically $-\infty$.
The decreasing convergence theorem for subharmonic functions gives
$L^1_{\mathrm{loc}}$ convergence. Distributional differentiation and
positivity then give the stated distributional and vague convergence.
\end{proof}

\begin{proof}[Proof of \cref{thm:algorithm}]
The construction in \cref{lem:linearization} gives
an admissible pencil explicitly; the argument applies to every
pencil satisfying \eqref{alg:det-normalization}. \Cref{alg:iteration} proves convergence
of the Dyson iteration and invertibility of the two derivative
systems. \Cref{alg:potential} gives the normalized determinant
formula, and \eqref{alg:density-formula} gives its Laplacian.
Finally, \cref{alg:recovery} proves nonnegativity and the asserted
pointwise, local $L^1$, and vague limits. The construction uses
$K_z$ without inverting it, so it also applies at $z=P(0)$.
\end{proof}

\subsection{A circular polynomial with an adjoint}
\label{alg:adjoint-example}
Let $c=(s_1+is_2)/\sqrt2$ be standard circular and $P=c^2+c^*$.
The matrix approximation in \cref{fig:adjoint-example} involves a
Ginibre matrix and its own adjoint. We now determine the boundary and
compute the density inside it using a pencil of size two.

\begin{proposition}\label[proposition]{alg:adjoint-measure}
The Brown measure of $P=c^2+c^*$ is absolutely continuous and has
real-analytic, strictly positive density on the interior of
$\supp\mu_P=\spec P$. The density is given by \eqref{alg:adjoint-density} below
using a uniquely specified solution of the unregularized Dyson equation. The support boundary is
\[
 \left\{g^{-1}+3g^2:g=u_*(\theta)e^{i\theta},\ 0\le\theta\le2\pi\right\},
\]
where $u_*(\theta)$ is the unique root in $(0,1/\sqrt2)$ of
\[
 d_\theta(u)=1-3u^2-2u^4-4u^6-8u^5\cos(3\theta).
\]
Moreover, $\mu_P(\partial\spec P)=0$.
\end{proposition}

The pencil
\begin{equation}\label{alg:adjoint-pencil}
 \begin{gathered}
 L_z=\begin{pmatrix}c^*-z&c\\c&-1\end{pmatrix}
 =K_z\otimes1+B\otimes c+B'\otimes c^*,\\
 K_z=\operatorname{diag}(-z,-1),\qquad
 B=\begin{pmatrix}0&1\\1&0\end{pmatrix},\qquad B'=E_{11}.
 \end{gathered}
\end{equation}
has Schur complement $P-z$. Since its auxiliary determinant has modulus
one, $U_P(z)=2\log\Delta_{\tau_2}(L_z)$ by determinant multiplicativity.
Write $H_z,J_z$ for the hermitizations of $L_z,K_z$. The coefficients in
$L_z$ of $s_1,s_2$ are $(B+B')/\sqrt2$ and $i(B-B')/\sqrt2$. Thus
\[
 \Phi(X)=BXB+B'XB',\qquad
 \Theta(X)=BXB'+B'XB,\qquad
 \mathcal S(N)=
 \begin{pmatrix}\Phi(N_{22})&\Theta(N_{21})\\
                 \Theta(N_{12})&\Phi(N_{11})\end{pmatrix}.
\]
Here $N_{ij}$ are $2\times2$ blocks. The expected resolvent
$M_\eta=(\id\otimes\tau)((i\eta I_4-H_z)^{-1})$ solves \eqref{alg:dyson},
and \eqref{alg:action} gives
$U_\eta(z)=\tau_2\log(L_z^*L_z+\eta^2 I_2)\downarrow U_P(z)$.

The off-diagonal Hermitization gives $M_\eta^*=-\Gamma M_\eta\Gamma$
for $\Gamma=\operatorname{diag}(I_2,-I_2)$. Uniqueness also makes
$M_\eta$ invariant under transposition followed by block interchange,
which preserves the Dyson equation. Thus $M_\eta$ and its limits have
the form
\begin{equation}\label{alg:adjoint-symmetry}
 M=\begin{pmatrix}
 -iq&-i\kappa&\bar g&\bar h\\
 -i\bar\kappa&-it&\bar h&\bar\ell\\
 g&h&-iq&-i\bar\kappa\\
 h&\ell&-i\kappa&-it
 \end{pmatrix},\qquad q,t\in\mathbb R,\quad \kappa,g,h,\ell\in\C.
\end{equation}
With $M=M_\eta$ and $Q=-\Im M$, the imaginary part of the Dyson
equation gives $Q=\eta MM^*+M\mathcal S(Q)M^*$.
Write $C=\left(\begin{smallmatrix}0&B'\\ B&0\end{smallmatrix}\right)$,
so $\mathcal S(N)=CNC^*+C^*NC$.

\begin{lemma}
\label[lemma]{lem:adjoint-compact}
The matrices $M_\eta$ stay bounded as $\eta\downarrow0$.
Every cluster limit is invertible and solves the Dyson equation at
$\eta=0$. If $q=0$ in a cluster limit, that limit is Hermitian.
\end{lemma}

\begin{proof}
Since $\mathcal S(Q)=CQC^*+C^*QC$, the imaginary-part
identity becomes
\[
Q=\eta MM^*
  +(MC)Q(MC)^*
  +(MC^*)Q(MC^*)^*.
\]
Set $T_1=Q^{-1/2}MCQ^{1/2}$ and
$T_2=Q^{-1/2}MC^*Q^{1/2}$.
Multiplying the preceding identity on both sides by
$Q^{-1/2}$ gives
\[
T_1T_1^*+T_2T_2^*
=I_4-\eta Q^{-1/2}MM^*Q^{-1/2}
\le I_4.
\]
Thus $\lVert T_1\rVert,\lVert T_2\rVert\le1$, and hence
$\lVert T_1T_2\rVert\le1$.
Similarity invariance of the trace yields
\[
\operatorname{Tr}(MC)=\operatorname{Tr}(T_1),
\qquad
\operatorname{Tr}(MCMC^*)=\operatorname{Tr}(T_1T_2).
\]
Using $|\operatorname{Tr}(T)|\le4\lVert T\rVert$
for $T\in M_4(\mathbb{C})$, we obtain the uniform bounds $|\Tr(MC)|\leq 4$ and $|\Tr(MCMC^*)|\le4$.
Expanding the traces using \eqref{alg:adjoint-symmetry} gives
\[
 |g+2\bar h|\le4,\qquad
 \tfrac12|g-2\bar h|^2+2|\kappa|^2+q^2+2qt
 =\tfrac12|g+2\bar h|^2-\Tr(MCMC^*)\le 12.
\]
Since $q,t\ge0$, this bounds $q,\kappa,g,h$. The auxiliary block of $M^{-1}$ in coordinates
$(2,4)$ is $\left(\begin{smallmatrix}i(q+\eta)&1\\1&i(q+\eta)\end{smallmatrix}\right)$,
with inverse norm at most one. The coupling blocks have entries among
$i\kappa,i\bar\kappa,-g,-\bar g$. The block inverse formula therefore
bounds $t,\ell$ by the already bounded primary block of $M$ in $(1,3)$.
Every cluster limit satisfies $M(-J_z-\mathcal S(M))=I_4$ and is invertible.
If $q=0$ at such a limit, then $\kappa=0$. If $t>0$, then
$\mathcal S(Q)=t\operatorname{diag}(E_{11},E_{11})$, and the imaginary-part
identity forces the primary block of $M$ to vanish. That block is the
inverse of the primary Schur complement, so $t=0$ and $Q=0$.
\end{proof}

\begin{lemma}
\label[lemma]{lem:adjoint-stability}
A solution of the Dyson equation at $\eta=0$ of the form
\eqref{alg:adjoint-symmetry} with $-\Im M\ge0$ and $q>0$
has $-\Im M>0$. It extends analytically in $(x,y,\eta)$ through
$\eta=0$ and is the unique limit of $M_\eta$ there.
\end{lemma}

\begin{proof}
Put $Q=-\Im M$. The Hermitized semicircular coefficients are
\[
 C_1=\tfrac1{\sqrt2}(C+C^*),\qquad C_2=\tfrac{i}{\sqrt2}(C^*-C).
\]
\noindent The matrix $C_1$ is invertible. In the standard basis,
\[
 \begin{gathered}
 (MC_1)^{-1}(MC_2)=
 \operatorname{diag}\!\left(
   -i\begin{pmatrix}1&0\\-2&1\end{pmatrix},
    i\begin{pmatrix}1&0\\-2&1\end{pmatrix}\right),\\
 (MC_1e_2)_3=(MC_1e_4)_1=-iq/\sqrt2\ne0.
 \end{gathered}
\]
A common invariant subspace for $MC_1,MC_2$ is invariant under
$(MC_1)^{-1}MC_2$, since $MC_1$ is invertible. Its two Jordan
blocks have distinct eigenvalues, so the subspace splits into its
parts in the two coordinate blocks. Each nonzero part contains the second coordinate axis.
If either part is nonzero, the second line forces the opposite block
to be full, and then so is the initial block.

The imaginary-part identity when $\eta=0$ is
\[
 Q=M\mathcal S(Q)M^*=\sum_{j=1}^2(MC_j)Q(MC_j)^*.
\]
The range of $Q$ is common invariant and nonzero, so $Q>0$.
Conjugation by $X\mapsto Q^{1/2}XQ^{1/2}$ makes the map unital.
Subtracting the least eigenvalue times $I_4$ from a Hermitian fixed
matrix shows that it must be scalar, since a nonzero positive fixed
matrix would have full range. Consequently
\[
 \begin{gathered}
 \{X:M\mathcal S(X)M^*=X\}=\C Q,\qquad
 M\mathcal S(X\Gamma)M=M\mathcal S(X)M^*\Gamma,\\
 \ker\bigl[N\mapsto N-M\mathcal S(N)M\bigr]=\C Q\Gamma .
 \end{gathered}
\]
The middle identity follows from $C_j\Gamma=-\Gamma C_j$ and $M^*=-\Gamma M\Gamma$. The diagonal blocks of
$Q\Gamma$ have opposite nonzero traces, so the residual derivative
$-M^{-1}(N-M\mathcal S(N)M)M^{-1}$ is invertible on the space
\eqref{alg:adjoint-symmetry}, which it preserves. The implicit function theorem
gives a real-analytic branch in $(x,y,\eta)$ through $\eta=0$. Its
imaginary part remains strictly negative, so uniqueness for $\eta>0$
identifies it with $M_\eta$ there.
\end{proof}

At $z=0$, substitution in the Dyson equation gives $q=\sqrt{k_0}$,
$t=q/(1+k_0)$ and $\ell=1/(1+k_0)$, where $k_0=(\sqrt5-1)/2$ and
$\kappa=g=h=0$. This solution has $-\Im M>0$ and therefore extends
to a positive branch near zero by \cref{lem:adjoint-stability}.

\begin{lemma}
\label[lemma]{lem:adjoint-boundary}
The boundary of $\spec P$ is the curve stated in
\cref{alg:adjoint-measure}. It is a simple closed real-analytic
curve and avoids the origin. The matrices $M_\eta$ converge at every
$z$, with a positive limit $-\Im M>0$ throughout the interior.
\end{lemma}

\begin{proof}
For a Hermitian cluster limit the bottom-left block $G$ is invertible
and symmetric. The Dyson equation and $G^{-1}G=I_2$ give
\[
 \begin{aligned}
 G^{-1}&=-K_z-\Theta(G)=\begin{pmatrix}z-2h&-g\\-g&1\end{pmatrix},\\
 G&=\begin{pmatrix}g&g^2\\g^2&1+g^3\end{pmatrix},
 \qquad z=g^{-1}+3g^2.
 \end{aligned}
\]
At positive height $Q=-\Im M_\eta$ satisfies
\[
 Q-M_\eta\mathcal S(Q)M_\eta^*=\eta M_\eta M_\eta^*>0.
\]
After conjugation by $X\mapsto Q^{1/2}XQ^{1/2}$, the positive map
$X\mapsto M_\eta\mathcal S(X)M_\eta^*$ sends $I_4$ strictly below
$I_4$, so its spectral radius is less than one. At a Hermitian
limit it restricts to $X\mapsto G\Phi(X)G^*$ on the lower diagonal
block. This map has matrix $T_g$ with $r(T_g)\le1$. For $g=u e^{i\theta}$,
\[
 \begin{gathered}
 T_g=\overline{GB}\otimes GB+\overline{GB'}\otimes GB',\\
 \det(I_4-T_g)=d_\theta(u),\qquad \Tr T_g=u^2+4u^4.
 \end{gathered}
\]
The full map dominates $X\mapsto GB'X(GB')^*$ on positive matrices,
so $r(T_g)\ge u^2$.
The trace excludes $u=1$ if $r(T_g)\le1$. The polynomial in
\cref{alg:adjoint-measure} satisfies
\[
 \begin{aligned}
 d_\theta(u)&\le (1-u)\bigl[4u^4(u-1)+(u+1)(1-2u^2)\bigr]<0
      &&(1/\sqrt2\le u<1),\\
 d_\theta'(u)&\le-2u(3-11u^2+12u^4)<0
      &&(0<u\le1/\sqrt2).
 \end{aligned}
\]
Since $r(T_g)\le1$ implies $\det(I_4-T_g)\ge0$, the first inequality
excludes $1/\sqrt2\le u<1$. The second gives the unique root $u_*(\theta)$
in $(0,1/\sqrt2)$ because $d_\theta(0)=1$. At $u=0$ the map is
nilpotent. Its radius is a continuous eigenvalue, so
\[
 r(T_g)<1\quad(0\le u<u_*),\qquad r(T_g)=1\quad(u=u_*),
 \qquad r(T_g)>1\quad(u>u_*).
\]

The trial inverse satisfies $(-G)^{-1}=K_z-\Theta(-G)$. The criterion
\cref{sp:invertibility} gives
\[
 0<u<u_*\ \Longrightarrow\ g^{-1}+3g^2\notin\spec P,\qquad
 u=u_*\ \Longrightarrow\ g^{-1}+3g^2
    \text{ is a limit of resolvent points}.
\]
Thus a Hermitian cluster limit can occur only in the closure of the
resolvent set. At an interior point compactness gives a limit with
$q>0$, and uniqueness follows from
\cref{lem:adjoint-stability}. At a resolvent point the
operator resolvent converges in norm to a selfadjoint inverse.
A limit with $q>0$ would persist on a neighborhood that meets the
resolvent set, so boundary limits are Hermitian and must have $u=u_*$. The positive branch
near zero puts zero in the interior. In particular $g^{-1}+3g^2\ne0$ whenever $u=u_*(\theta)$.

For $z\ne0$, the semicircular coefficients $A_1,A_2$ of $K_z^{-1}L_z$
span the same space as $K_z^{-1}B$ and $K_z^{-1}B'$. The second
has the coordinate eigenlines, which the first exchanges. The
tuple is therefore irreducible, and
\[
 A_jS=SA_j^*\ (j=1,2)\ \Longrightarrow\
 (K_z^{-1}B)S=S(K_z^{-1}B')^*\ \Longrightarrow\
 \operatorname{rank}S\le1 .
\]
The range would be a common invariant line if $S\ne0$, so $S=0$.
Normalization replaces a trial inverse $R$ by $RK_z$ and preserves
the matrices in the spectral-radius test because
\[
 (RK_z)A_j=R(K_zA_j).
\]
By \cref{sp:root-uniqueness}, there is at most one
solution of the normalized equation for which the
associated covariance map has spectral radius
at most one.
When $u=u_*(\theta)$, the solution constructed above
has spectral radius one for its covariance map.
By uniqueness, there is no solution for which this
spectral radius is strictly less than one.
Hence, by \cref{sp:invertibility},
$z=g^{-1}+3g^2$ belongs to the spectrum,
and therefore to its boundary.
The same uniqueness gives convergence at the boundary
and injectivity modulo $2\pi$, since $G_{11}=g$.
Since $d_\theta'(u_*)<0$, the implicit function theorem
implies that $u_*$ is real analytic.
The derivative of $u_*(\theta)e^{i\theta}$ is nonzero
since $u_*>0$.
The map $g\mapsto g^{-1}+3g^2$ has zero derivative
only at $g^3=1/6$.
In those directions $d_\theta(1/2)=-3/16$ gives
$u_*<1/2<6^{-1/3}$, so the boundary parametrization
is regular.
\end{proof}

\begin{lemma}
\label[lemma]{lem:adjoint-density}
The density on the interior of $\spec P$ is real analytic,
strictly positive, and given by \eqref{alg:adjoint-density}.
\end{lemma}

\begin{proof}
Along the analytic branch $M$, \eqref{alg:action}
extends real analytically across $\eta=0$,
with value $U_P$ at $\eta=0$.
Differentiating the regularized logarithm gives
$2\partial_z U_\eta=(M_\eta)_{31}$.
At $\eta=0$, this becomes
$2\partial_z U_P=g=M_{31}$, so
\begin{equation}\label{alg:adjoint-density}
 \rho_P(z)=\frac1{2\pi}\Re(g_x+i g_y).
\end{equation}
For $a=x,y$, compute $g_a=(\partial_a M)_{31}$ from
$-M^{-1}(\partial_a M)M^{-1}+\mathcal S(\partial_a M)=-\partial_aJ_z$.
\Cref{lem:adjoint-stability} gives a unique solution in
\eqref{alg:adjoint-symmetry}.

For interior $z\ne0$, the proof of \cref{lem:adjoint-boundary} verifies the
hypotheses of the coefficient regularity theorem for $K_z^{-1}L_z$
(\cref{sp:regularity}). If $a>0$ is the least eigenvalue
of $K_z^*K_z$, then $L_z^*L_z\ge a(K_z^{-1}L_z)^*(K_z^{-1}L_z)$.
Inverting after adding $\eta^2 I_2$ shows that the normalized lower
diagonal block of the negative imaginary part at height
$\eta/\sqrt a$, multiplied by $a^{-1/2}$, dominates the original
one. Thus the normalized limit lies in the positive regime.

For the normalized tuple at a nonzero interior point $z_0$, consider
$\mathscr F(z,L,E)=F(I+L,\mathbf A(z)+E)$ as in the proof of
\cref{sp:positive-rigidity}. This is real-analytic and plurisubharmonic
near $(z_0,0,0)$. There $2\mathscr F(z,0,0)=U_P(z)-\log|\det K_z|$,
so a zero of $\rho_P$ would give $\mathscr F_{z\bar z}=0$ at the point.
Positivity of its Levi matrix forces
$\mathscr F_{\xi\bar z}=0$ there for every coefficient or mean
coordinate $\xi$. By \eqref{sp:complex-gradient}, this gives
$\partial_{\bar z}R=0$ and $\partial_{\bar z}(XA_j^*Y)=0$ at $z_0$.
Step 3 of the proof of \cref{sp:positive-rigidity} uses only these
first derivatives at the point, and yields
$\partial_zA_j=[(\partial_zX)X^{-1},A_j]$. Taking the trace contradicts
$\Tr\partial_zA_1(z_0)=1/(\sqrt2 z_0^2)$, since
$A_1=K_z^{-1}(B+B')/\sqrt2$.

At the origin, the first-order auxiliary and primary Schur equations
at the explicit solution above give
$d h=-2k_0\,d\bar g/(1+2k_0)$ and $d z=k_0^{-1}d\bar g+2\,d h$.
Substituting into \eqref{alg:adjoint-density} gives
\[
 \rho_P(0)=\frac{k_0(1+2k_0)}{\pi(1+2k_0-4k_0^2)}
          =\frac{5+3\sqrt5}{6\pi}>0.
\]
\end{proof}

\begin{proof}[Proof of \cref{alg:adjoint-measure}]
By \cref{sp:main}, the support of the Brown measure coincides with the spectrum.
\Cref{lem:adjoint-boundary,lem:adjoint-density} give its boundary
and the interior density. Near the boundary $z\ne0$. With $K_z,B,B'$
as in \eqref{alg:adjoint-pencil}, put
\[
 A_1(z)=K_z^{-1}(B+B')/\sqrt2,\qquad
 A_2(z)=iK_z^{-1}(B-B')/\sqrt2 .
\]
Factoring $K_z$ from $L_z$ gives
\[
 U_P(z)=\log|z|+2\log\Delta_{\tau_2}\bigl(
 I_2\otimes1+A_1(z)\otimes s_1+A_2(z)\otimes s_2\bigr).
\]
The normalized tuple satisfies the hypotheses of \cref{sp:regularity}
near the boundary, so $U_P$ is $C^1$ there.
The boundary is a compact smooth curve and hence carries no mass by
\cref{sp:no-curve-mass}. This proves absolute continuity.
\end{proof}

The left panel of \cref{fig:adjoint-example} evaluates
\eqref{alg:adjoint-density} on the interior branch, starting at the
explicit solution at the origin. For the Ginibre samples on the right,
\cref{thm:main,thm:hausdorff,sp:main} give almost-sure weak convergence
to $\mu_P$ and almost-sure Hausdorff convergence to its support.

\appendix
\section{Regularity and absolute continuity for circular polynomials}
\label[appendix]{sp:circular}

\begin{theorem}
\label[theorem]{sp:circular-main}
Let $c_1,\ldots,c_k$ be a free standard circular family, and let
$P\in\C\langle x_1,\ldots,x_k\rangle$ be a nonconstant polynomial.
Set $a=P(c)$. Then the Brown measure of $a$ is absolutely continuous
with respect to planar Lebesgue measure. Its logarithmic potential
$U_a$ is finite everywhere and belongs to $C^1(\C\setminus E)$
for some finite set $E$. Moreover, for every $z_0\notin E$,
there are a neighborhood $D$ of $z_0$ and a nonzero polynomial
$p\in\R[x,y]$ such that $U_a$ is real-analytic on
$D\setminus\{z:p(\Re z,\Im z)=0\}$.
\end{theorem}

The exclusion of adjoints is essential. The Brown measure of the self-adjoint element
$(c_1+c_1^*)/\sqrt2$ is supported on the real line.

\subsection{The circular pencil and regularity of the potential}
\label{sp:circular-proofs}

For circular inputs ($\rho=0$), the doubled tuple has $\Theta=0$ by
\eqref{sp:elliptic-covariance}. The equation for $R$ in
\cref{sp:elliptic-invertibility} reduces to $R^{-1}=K$, giving
\begin{equation}\label{sp:circular-invertibility}
 K\otimes1+\sum_jB_j\otimes c_j\text{ is invertible}
 \quad\Longleftrightarrow\quad
 K\text{ is invertible and }r(\Phi_{K^{-1}\mathbf B})<1 .
\end{equation}

\begin{proposition}\label[proposition]{sp:circular-regularity}
Let $\mathbf A(z)=(A_1(z),\ldots,A_k(z))$ be a real-analytic family of
nonzero irreducible tuples on an open $\Omega\subset\C$. Put
\[
 V(z)=\log\Delta_m\Big(I+\sum_j A_j(z)\otimes c_j\Big),\qquad
 D=\{z\in\Omega:r(\Phi_{\mathbf A(z)})>1\}.
\]
Then $V$ is $C^1$ on $\Omega$, positive and real-analytic on $D$, and
zero on $\Omega\setminus D$.
\end{proposition}

\begin{proof}
By \cref{sp:elliptic-reduction}, the associated semicircular tuple
satisfies the hypotheses of \cref{sp:regularity}. Thus $V\in C^1(\Omega)$
and $V$ is real-analytic on the set $\mathcal O_+$ defined there.
For $z\in\mathcal O_+$, the matrices $X,Y$ are positive definite, and
\eqref{sp:schur} at $\eta=0$ gives
\[
 \Phi(X)-X=\Phi(X)Y\Phi(X)>0,
 \qquad \Phi=\Phi_{\mathbf A(z)}.
\]
Pairing this inequality with a positive definite Perron eigenmatrix
of $\Phi^*$ (\cref{sp:perron-eigenmatrices}) shows that $r(\Phi)>1$.
For $z\notin\mathcal O_+$, \cref{sp:regularity}(ii) gives $R=I$ and
$r(\Phi)\le1$. Hence $\mathcal O_+=D$.

Now fix $z\in\Omega$ and write $\Phi=\Phi_{\mathbf A(z)}$. 
For the scaled tuple $t\mathbf A(z)$, $t>0$, the covariance
map is $\Phi_t=t^2\Phi$. Thus $r(\Phi_t)>1$ if and only if
$t>t_0:=r(\Phi)^{-1/2}$.
Set
\[
 f(t)=\log\Delta_m\Big(I+t\sum_j A_j(z)\otimes c_j\Big).
\]
By \eqref{sp:limit-gradient}, $f'(t)=0$ for $0<t\le t_0$, whereas
\[
 f'(t)=\frac1{mt}\Tr(Y_t\Phi_t(X_t))>0
 \qquad (t>t_0).
\]
Since $f(t)\to0$ as $t\downarrow0$, it follows that $f(1)=0$ if
$r(\Phi)\le1$ and $f(1)>0$ otherwise.
\end{proof}

For any constant coefficient tuple, including a reducible one, the
potential $\log\Delta_m(I+\sum_j A_j\otimes c_j)$ is nonnegative:
$w\mapsto\log\Delta_m(I+w\sum_j A_j\otimes c_j)$ is radial and
subharmonic and takes the value zero at the origin.

\begin{corollary}
\label[corollary]{sp:circular-rational}
If $A_j\in M_m(\C(z))$, the potential
$V(z)=\log\Delta_m(I+\sum_j A_j(z)\otimes c_j)$ is $C^1$
off a finite set $E$ containing the poles. Locally off $E$ it is
real-analytic away from the zeros of a nonzero polynomial in
$x=\Re z$ and $y=\Im z$.
\end{corollary}

\begin{proof}
By \cref{sp:blocks}, there is a finite set $E$ such that, on a
sufficiently small disk in $\C\setminus E$, the coefficient tuple is
simultaneously block triangularizable, with holomorphic entries
algebraic over $\C(z)$. By \eqref{sp:block-potential}, $V$ is a
weighted sum of the potentials associated with the diagonal blocks.
A zero block contributes zero, while each nonzero block is irreducible.
Thus \cref{sp:circular-regularity} implies $V\in C^1(\C\setminus E)$.

Consider a nonzero diagonal block and let $\rho=r(\Phi)$ be the
Perron eigenvalue of its covariance map. By
\cref{sp:perron-eigenmatrices}, $\rho$ is simple and depends
real-analytically on $z$. It is algebraic over $\R(x,y)$, with
$z=x+iy$, since the matrix of $\Phi$ on the real vector space of
Hermitian matrices has algebraic entries in $x,y$. If $\rho\equiv1$,
\cref{sp:circular-regularity} shows that the block potential vanishes.
Otherwise, the minimal polynomial of $\rho$ over $\R(x,y)$, after
clearing denominators, gives
\[
 q(x,y,\rho)=0,\qquad q(x,y,1)\not\equiv0.
\]
Consequently, $\{\rho=1\}$ is contained in the zero set of the
nonzero polynomial $q(x,y,1)$. By \cref{sp:circular-regularity},
the block potential is real-analytic away from this set. Multiplying
the finitely many resulting polynomials (with the empty product
understood as $1$) proves the claim.
\end{proof}

\begin{lemma}\label[lemma]{sp:no-curve-mass}
If the logarithmic potential $u$ of a compactly supported finite
positive measure $\mu$ is $C^1$ on $\Omega$, then
$\mu(B(z,r))=o(r)$ uniformly on compact subsets of $\Omega$.
In particular every compact $C^1$ arc in $\Omega$ has zero $\mu$-mass.
\end{lemma}
\begin{proof}
Fix a compact $K\subset\Omega$. The circular mean identity used in
\cref{lem:disk} gives, uniformly for $z\in K$,
\[
 (\log 2)\mu(B(z,r))
 \le \frac1{2\pi}\int_0^{2\pi}u(z+2re^{i\theta})\dd\theta-u(z)
 =o(r).
\]
Indeed, the linear term of the first-order expansion has mean zero,
and $\nabla u$ is uniformly continuous on a compact neighborhood of
$K$. A compact $C^1$ arc is covered by $O(1/r)$ disks of radius $r$
centered on the arc. Summing the bound proves the claim.
\end{proof}

\begin{proof}[Proof of \cref{sp:circular-main}]
By \eqref{sp:normalized-pencil}, $U_a(z)$ is $\log|P(0)-z|$ plus
$m$ times the potential of a circular pencil with rational coefficient
matrices in $z$. \Cref{sp:circular-rational} gives the stated regularity
after adjoining $P(0)$ to $E$. Fix such a disk and let $p$ denote the corresponding polynomial.

The zero set of $p$ in the disk is covered by countably many
compact $C^1$ arcs. Indeed, let $(x_0,y_0)$ be a zero of $p$,
and choose a multi-index $\alpha$ of minimal order such that
$\partial^\alpha p(x_0,y_0)\ne0$. There is a multi-index
$\beta$ with $|\beta|=|\alpha|-1$ such that
\[
 \partial^\beta p(x_0,y_0)=0,\qquad
 \nabla(\partial^\beta p)(x_0,y_0)\ne0.
\]
Thus every zero of $p$ lies in the regular zero set of some
partial derivative of $p$. These regular zero sets are locally
$C^1$ curves and can be covered by countably many compact
$C^1$ arcs.

By \cref{sp:no-curve-mass}, these arcs have zero $\mu_a$-measure.
On $\{p\ne0\}$, the measure has density $(2\pi)^{-1}\Lap U_a$.
A countable covering argument gives absolute continuity on
$\C\setminus E$. Finally, \cref{sp:finite-potential} implies that
$U_a$ is finite everywhere and hence that $\mu_a$ has no atoms.
In particular, $\mu_a(E)=0$.
\end{proof}

The exceptional set need not be empty. For example, let $P(c)=c^2$.
For $w^2=z$, multiplicativity of the Fuglede--Kadison determinant
gives $U_{c^2}(z)=U_c(w)+U_c(-w)$. By the circular law \cite{HL},
\[
 U_{c^2}(z)=
 \begin{cases}
 |z|-1,&|z|\le1,\\
 \log|z|,&|z|\ge1.
 \end{cases}
\]
Thus $U_{c^2}$ is finite everywhere but is not differentiable at
$z=0$.

\subsection{A mixed-degree polynomial with disk support}
\label{alg:quadratic-example}
Let $c,d$ be free standard circular elements and put $P=c^2+cd+d$.
The linear and quadratic terms transform differently under a common
rotation of $c,d$. Nevertheless the spectrum is a disk. To see why,
factor $P=c(c+d)+d$ and take the Schur complement at the lower-right
entry of the pencil
\begin{equation}\label{alg:example-pencil}
 \begin{gathered}
 L_z=\begin{pmatrix}d-z&c\\-c-d&1\end{pmatrix}
 =K_z+B_c c+B_d d,\\
 K_z=\operatorname{diag}(-z,1),\qquad
 B_c=\begin{pmatrix}0&1\\-1&0\end{pmatrix},\qquad
 B_d=\begin{pmatrix}1&0\\-1&0\end{pmatrix}.
 \end{gathered}
\end{equation}
For $z\ne0$ the covariance of the normalized pencil sends the diagonal
entries $(h,\ell)$ of a Hermitian matrix to
$((h+\ell)/|z|^2,2h)$, independently of its off-diagonal entries.
Its spectral radius is therefore the Perron eigenvalue of
\[
 \begin{pmatrix}|z|^{-2}&|z|^{-2}\\2&0\end{pmatrix}.
\]
Indeed the iterates of the covariance on $I$ are positive, so their
norms lie between half their trace and their trace, and these traces
are determined by the displayed matrix. Its Perron eigenvalue is less
than one precisely when $|z|^2>3$. Applying
\eqref{sp:circular-invertibility} with $K=K_z$ to this $2\times2$ pencil
gives invertibility for $z\ne0$ exactly when $|z|>\sqrt3$.
At the origin the same criterion rules out invertibility because $K_0$
is singular. Taking the Schur complement and applying \cref{sp:main} gives
\[
 \spec P=\supp\mu_P=\{z:|z|\le\sqrt3\}.
\]

\begin{figure}[htbp]
 \centering
 \includegraphics[width=\textwidth]{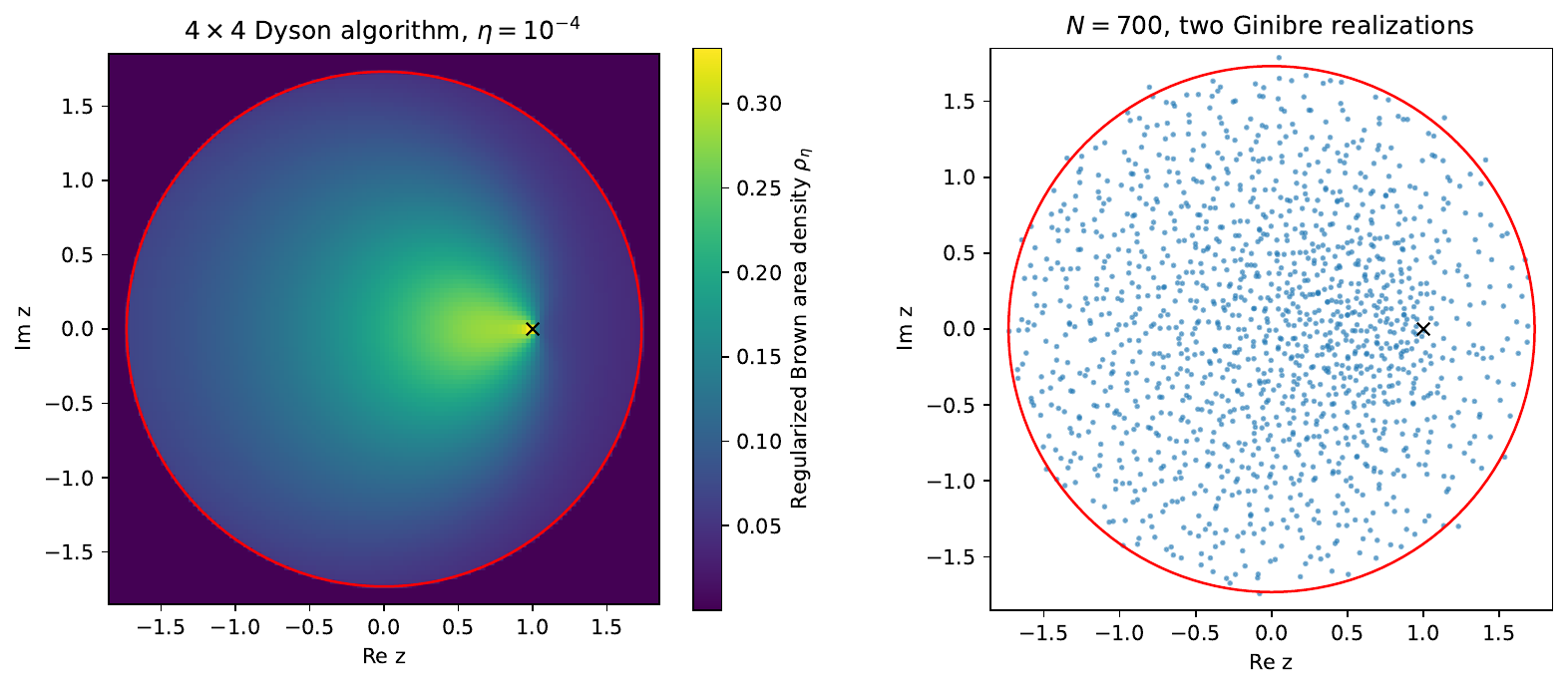}
 \caption{The polynomial $c^2+cd+d$. Left: the regularized density at
 $\eta=10^{-4}$ from the Hermitization of \eqref{alg:example-pencil}.
 Right: the eigenvalues of $C_N^2+C_ND_N+D_N$ for two independent samples
 at $N=700$, with $C_N,D_N$ independent complex Ginibre matrices.
 The circle has radius $\sqrt3$ and the cross marks $z=1$.}
 \label{fig:quadratic-density}
\end{figure}

The numerically computed density in \cref{fig:quadratic-density}
varies with direction. To compute the limiting measure, write
\[
 b=1+c,\quad C_0=bb^*,\quad w_z=z-c^2,\quad
 Q(z,k)=w_zw_z^*+kC_0,\qquad F(z,k)=\tau\log Q(z,k).
\]
Our \cref{thm:main,thm:hausdorff}
give almost-sure convergence of both the measures and the spectra (see also \cite{CGH, Han}),
and the determinant formula identifies the limit explicitly.

\begin{proposition}\label[proposition]{alg:example-measure}
For $|z|<\sqrt3$, $z\ne1$, there is a unique $k=k(z)\in(0,1)$ such that
$F_k(z,k)=\tau(Q(z,k)^{-1}C_0)=1$. There
\begin{equation}\label{alg:example-density}
 U_P(z)=\tfrac12(F(z,k)-k),\qquad
 \rho_P(z)=\frac1{4\pi}\left(F_{xx}+F_{yy}
             -\frac{F_{xk}^2+F_{yk}^2}{F_{kk}}\right)(z,k)>0 ,
\end{equation}
where $x=\Re z$ and $y=\Im z$. The density is real-analytic on this
punctured disk. For $|z|\ge\sqrt3$, $U_P(z)=\log|z|$, and
$U_P(1)=\tfrac12\log2-\tfrac14$.
\end{proposition}

\begin{proof}
\textit{Step 1: reduction to a scalar equation.}
Write $H_z$ for the Hermitization of $L_z$ and $H_c$ for its value at
$d=0$. The two circular covariances add, with
$\mathcal S_a(N)=\operatorname{diag}(B_aN_{22}B_a^*,B_a^*N_{11}B_a)$
for $a=c,d$. Put $E=E_{11}\in M_2(\C)$,
$J=\operatorname{diag}(E,E)$, and $\Lambda=\varepsilon J+\delta I_4$
for $\varepsilon,\delta>0$. In the Dyson equation for the expected
resolvent $M=(\id\otimes\tau)((i\Lambda-H_z)^{-1})$, move
$\mathcal S_d(M)$ into the deterministic part. The resulting
$c$-pencil has the same expected resolvent by uniqueness for the Dyson
equation \cite{HRS}. Put $\Lambda_c=\Lambda+i\mathcal S_d(M)$.
Both $\Lambda$ and $\Lambda_c$ are positive definite
and block diagonal, since the diagonal blocks of $M$
are skew-adjoint with negative imaginary part.
Apply \eqref{alg:action} to both pencils after
normalizing their respective regularization matrices
to $I_4$ by congruence.
The determinant factors arising from these congruences
cancel the corresponding changes in $\log|\det M|$.
Subtraction gives
\[
 2\tau_4\log|i\Lambda-H_z|
 =2\tau_4\log|i\Lambda-H_c-\mathcal S_d(M)|
       +\tfrac14\Re\Tr\big(M\mathcal S_d(M)\big).
\]
Here $\tau_4=4^{-1}(\Tr\otimes\tau)$. At $\delta=0$ the auxiliary
block of the full pencil has determinant of modulus one and its
Schur complement is the Hermitization of $P-z$ at height $\varepsilon$.
Thus the full resolvent has a norm limit and the left side tends to
$V_\varepsilon=\tfrac12\tau\log(|P-z|^2+\varepsilon^2)$.

For the limiting $M$ set $u=(1,-1)^{\mathsf T}$,
$\alpha=-\Im(M_{22})_{11}$, $\beta=-\Im(u^*M_{11}u)$, and
$s=\varepsilon+\beta$. Then
$-\mathcal S_d(M)=i\operatorname{diag}(\alpha uu^*,\beta E)$,
and the primary resolvent bound gives $0\le\varepsilon\alpha\le1$.
The auxiliary block of the shifted $c$-pencil is
$\left(\begin{smallmatrix}i\alpha&-1\\-1&0\end{smallmatrix}\right)$.
Its Schur complement and the further complement at $is$ are
\[
 \begin{pmatrix}i(\varepsilon+\alpha C_0)&w_z\\w_z^*&is\end{pmatrix},
 \qquad iQ_\varepsilon/s,\qquad
 Q_\varepsilon=Q(z,k_\varepsilon)+\varepsilon sI,\quad
 k_\varepsilon=s\alpha .
\]
Both are invertible. Hence the shifted determinant also passes to
$\delta=0$ and equals $\tfrac12\tau\log Q_\varepsilon$ with the factor
$2\tau_4$ above. The covariance term is $-\alpha\beta/2$, giving
\[
 V_\varepsilon=\tfrac12\big(\tau\log Q_\varepsilon
                       -k_\varepsilon+\varepsilon\alpha\big).
\]
Vary the two shifts $\alpha,\beta$ independently in the shifted
determinant before evaluating at $M$. The resolvent gives the
derivatives $\beta/2,\alpha/2$, respectively. Differentiating its
Schur complement expression gives
\begin{equation}\label{alg:disk-scalars}
 \tau(Q_\varepsilon^{-1}C_0)=1-\varepsilon/s,\qquad
 \alpha=s\tau(Q_\varepsilon^{-1}),\qquad
 k_\varepsilon=s^2\tau(Q_\varepsilon^{-1}).
\end{equation}

\textit{Step 2: the root.}
For $z\ne1$, the identity $w_z+b(c-1)=(z-1)1$ gives
\[
 Q(z,k)\ge |z-1|^2(1+\|c-1\|^2/k)^{-1}I,\qquad k>0,
\]
by the right inverse of the row $(w_z,\sqrt{k}b)$. In particular
$F$ is real-analytic for $k>0$, and trace differentiation gives
$F_{kk}=-\tau((Q^{-1/2}C_0Q^{-1/2})^2)<0$. Also $F_k(z,1)<1$,
since $Q(z,1)-C_0=w_zw_z^*\ne0$. For $|z|>1$ we have
$r(c^2/z)<1$. Expanding $w_z^{-1}b$ in norm and using
$\tau((c^m)^*c^n)=\delta_{mn}$ gives
\[
 F_k(z,0)=\|w_z^{-1}b\|_2^2
 =\left\|\sum_{j\ge0}z^{-j-1}(c^{2j}+c^{2j+1})\right\|_2^2
 =\frac{2}{|z|^2-1}.
\]
For $|z|\le1$, $z\ne1$, take $\zeta^2=z$. The scalar circular Dyson
equation for $q=\eta\tau(((\zeta-c)(\zeta-c)^*+\eta^2)^{-1})$ is
$\eta/q=|\zeta|^2+(\eta+q)^2-1$. Thus $q/\eta\to\infty$ as
$\eta\downarrow0$: this ratio is monotone, and a bounded limit would
give a left side bounded away from zero and a nonpositive limit on the right. Since
\[
 Q(z,k)\le\|\zeta+c\|^2(\zeta-c)(\zeta-c)^*+k\|C_0\|I,
\]
also $\tau(Q^{-1})\to\infty$ as $k\downarrow0$. Multiplying the row
identity by $Q^{-1/2}$ gives
$|z-1|\sqrt{\tau(Q^{-1})}\le1+\|c-1\|\sqrt{F_k}$.
This proves $F_k(z,0+)=\infty$ here. Since $F_k$ is strictly decreasing
in $k$, these endpoint values give exactly one positive root for
$|z|<\sqrt3$, $z\ne1$, and none for $|z|\ge\sqrt3$.

\textit{Step 3: removal of the regularization.}
Since $Q_\varepsilon\ge k_\varepsilon C_0$,
\eqref{alg:disk-scalars} gives $k_\varepsilon\le1+\varepsilon\alpha\le2$.
It also gives $\beta=s\tau(Q_\varepsilon^{-1}C_0)\le\|b\|^2\alpha$.
Thus $\varepsilon s$ and $\|Q_\varepsilon\|$ stay bounded, and
$s^2\le k_\varepsilon\|Q_\varepsilon\|$ gives $\varepsilon s\to0$.
At a positive cluster value $k_*$ the gap for $Q(z,k_*)$ bounds
$\alpha$ and gives $\varepsilon/s=\varepsilon\alpha/k_\varepsilon\to0$.
Hence $F_k(z,k_*)=1$. At a zero cluster value, inverse order in
\eqref{alg:disk-scalars} gives $F_k(z,k)\le1$ for every $k>0$:
eventually $Q_\varepsilon\le Q(z,k)+\varepsilon sI$ along that
subsequence. This excludes zero inside the punctured disk. Thus
$k_\varepsilon\to k(z)$ there and $k_\varepsilon\to0$ for
$|z|\ge\sqrt3$. In both cases the limiting $Q$ is invertible, so the
determinant identity for $V_\varepsilon$ gives the claimed potential.
In the exterior and on the circle,
$\tfrac12 F(z,0)=\log\Delta(z-c^2)=\log|z|$ by the circular law.
At $z=1$, $1-P=(1+c)(1-c-d)$. Since $c+d$ is circular of variance
two and $\log\Delta(1+c)=0$, multiplicativity gives the stated value.

Finally the implicit function theorem applies to $F_k=1$ since
$F_{kk}<0$. Differentiating $F(z,k(z))-k(z)$ twice gives
\eqref{alg:example-density}. For strict positivity let $T=Q^{-1/2}w_z$.
Trace differentiation gives
\[
 F_{xx}+F_{yy}=4\tau\big(Q^{-1}(I-T^*T)\big)>0 .
\]
Indeed $TT^*\le I$. Equality in the trace would give $T^*T=I$
and hence $TT^*=I$ by finiteness, contrary to $Q-w_zw_z^*=kC_0\ne0$.
Absolute continuity from \cref{sp:circular-main} rules out any
remaining mass at $1$ or on the boundary circle.
\end{proof}

The circular-addition formula \cite[Proposition 3.5 and Example 5.4]{BercoviciZhong}
gives an alternative derivation by factoring $z-P$ in the affiliated
determinant class as $(1+c)((1+c)^{-1}(z-c^2)-d)$.

\section*{Acknowledgements}
B.~Collins acknowledges support from the Ministry of Education, Culture,
Sports, Science and Technology (MEXT) through SPReaD (Supporting
Pioneering Research through AI for 1,000 Discovery challenges).
He acknowledges inspiring discussions with Charles Bordenave and F\'elix Parraud.

B.~Collins and Y.~Kitagawa are grateful for discussions with Amir Dembo and Dan Voiculescu around non-normal random matrices during their visit to Kyoto in the fall of 2025, at the early stages of the project.

P.~Zhong is supported in part by NSF grant LEAPS-MPS-2516951 and NSF CAREER Award DMS-2516950. He is grateful to Serban Belinschi, Hari Bercovici, Brian Hall, Yi Han, Ching-Wei Ho, and Roland Speicher for useful and inspiring discussions on various occasions.

AI disclosure: several coauthors had independently obtained partial results on projects related to this manuscript. Discussions among the coauthors, followed by interactions with frontier LLMs based on these partial results, helped connect them into the more general convergence results presented here.
The coauthors reviewed, improved, and substantially simplified all results. Frontier LLMs were also used to check proofs and improve mathematical consistency. The authors take responsibility for the contents of the paper.

\printbibliography
\end{document}